\documentclass[11pt]{article}

\usepackage{color, amsmath, amssymb, amsthm, graphicx, bbm, footnote, mathtools}

\usepackage[top=1in, bottom=1in, left=1in, right=1in]{geometry}

\usepackage[toc,page,header]{appendix}
\usepackage{titletoc}
\usepackage{tikz}
\usepackage{titling}
\usepackage{enumitem}
\allowdisplaybreaks

\usepackage{setspace}
\usepackage{adjustbox}
\usepackage[colorlinks,allcolors=blue]{hyperref}
\usepackage[round]{natbib}  
  \hypersetup{
	breaklinks=true,   
	colorlinks=true,   
	pdfusetitle=true,  
}

\usepackage{apptools}
\AtAppendix{\counterwithin{theorem}{section}}
\AtAppendix{\counterwithin{remark}{section}}

\newtheorem{theorem}{Theorem}
\newtheorem{definition}{Definition}
\newtheorem{proposition}[theorem]{Proposition}
\newtheorem{corollary}[theorem]{Corollary}
\newtheorem{lemma}[theorem]{Lemma}
\newtheorem{remark}{Remark}
\newtheorem{assumption}{Assumption}

\newtheorem{example}{Example}
\newcommand{\possessivecite}[1]{\citeauthor{#1}'s (\citeyear{#1})}

\newcommand{\ind}{{\perp\!\!\!\perp}}
\newcommand{\esssup}{{\mathrm{ess}\sup}}

\title{\textsc{A Bootstrap Hypothesis Test for High-Dimensional Mean Vectors} }
\author{Alexander Giessing\thanks{Department of Statistics, University of Washington, Seattle, WA. E-mail: giessing@uw.edu.} \and Jianqing Fan\thanks{Department of ORFE, Princeton University, Princeton, NJ. E-mail: jqfan@princeton.edu.}}
\date{\today}
\begin{document}
\maketitle

\begin{abstract}
	This paper is concerned with testing global null hypotheses about population mean vectors of high-dimensional data. Current tests require either strong mixing (independence) conditions on the individual components of the high-dimensional data or high-order moment conditions. In this paper, we propose a novel class of bootstrap hypothesis tests based on $\ell_p$-statistics with $p \in [1, \infty]$ which requires neither of these assumptions. We study asymptotic size, unbiasedness, consistency, and Bahadur slope of these tests. 
	Capitalizing on these theoretical insights, we develop a modified bootstrap test with improved power properties and a self-normalized bootstrap test for elliptically distributed data. We then propose two novel bias correction procedures to improve the accuracy of the bootstrap test in finite samples, which leverage measure concentration and hypercontractivity properties of $\ell_p$-norms in high dimensions. Numerical experiments support our theoretical results in finite samples.
\end{abstract}

\section{Introduction}\label{sec:Intro}

\subsection{Testing high-dimensional global null hypotheses}\label{subsec:GlobalNull}
Let $X_1, \ldots, X_n \in \mathbb{R}^d$ be i.i.d. random vectors with unknown mean $\mu_0 \in \mathbb{R}^d$ and unknown positive semi-definite covariance matrix $\Sigma \neq \mathbf{0} \in \mathbb{R}^{d \times d}$. In this paper, we consider testing high-dimensional linear restrictions
\begin{align}\label{eq:sec:Intro-1}
	H_0: \:\: R\mu_0 = r \quad{}\quad{} \mathrm{vs.} \quad{}\quad{} H_1: \:\: R \mu_0 \neq r,
\end{align}
for some deterministic matrix $R \in \mathbb{R}^{t \times d}$ and vector $r \in \mathbb{R}^t$ when dimension $d$ and number of restrictions $t$ are potentially much larger than the sample size $n$, i.e. $d, t \gg n$. 

Global null hypothesis testing problems such as~\eqref{eq:sec:Intro-1} arise frequently in scientific applications. For instance, many biological processes involve regulation of multiple genes with small effect sizes on individual variants~\citep{manolio2009finding}. In such cases, analyzing genes grouped according to their biological functions or chromosomal location can increase the power of statistical tests, reduce the complexity of the analysis, and lead to a better understanding of the underlying genetic mechanisms~\citep{huang2022overview}. The problem of assessing whether a group of many genes is differentially expressed from another group of genes ultimately reduces to testing a high-dimensional population mean vector as in~\eqref{eq:sec:Intro-1}.
Another example are cross-sectional and longitudinal studies in financial econometrics~\citep{fan2015power, gagliardini2016time-varying}. In these studies the parameter of interest is usually not a simple population mean but a high-dimensional regression vector. The time series character of such data sets adds additional complications. However, often, the testing problems can still be recast as in~\eqref{eq:sec:Intro-1}.

There are two common approaches to testing~\eqref{eq:sec:Intro-1} in high dimensions: one based on $\ell_2$-norms, i.e. sum-of-squares of the sample mean~\citep[e.g.][]{bai1996effect, chen2010two, wang2015high-dimensional, xu2017adaptive, he2021Asymptotic, huang2022overview}, and the other one based on $\ell_\infty$-norms, i.e. component-wise maxima  of the sample mean~\citep[e.g.][]{chernozhukov2013GaussianApproxVec, cai2014two, chen2018GaussianApproxUStat, xue2020Distribution, lopes2020bootstrapping}. Typically, tests based on the $\ell_2$-norm are valid only under mixing (independence) conditions on the individual components of the high-dimensional random vectors, whereas tests based on the $\ell_\infty$-norm often require higher-order moment conditions.

In this paper, we develop a new class of bootstrap hypothesis tests based on $\ell_p$-norms, $p \in [1, \infty]$, which is valid without mixing conditions and under mild lower-order moment conditions. Along with a thorough analysis of the theoretical properties of the test (Sections~\ref{subsec:Assumptions}-\ref{subsec:BahadurSlope}) we address the following questions:
\begin{itemize}
	\item It is statistical folklore that $\ell_2$-type statistics have high power against dense alternatives whose signals are spread out over a large number of coordinates while $\ell_\infty$-type statistics have high power against sparse alternatives with only a few strong signals~\citep[][]{cai2014two, fan2015power, wang2015high-dimensional, he2021Asymptotic}. \emph{Does this apply to the bootstrap test as well and what is the use of other $\ell_p$-norms, $p \geq 1$? Can we combine $\ell_p$-norms to obtain a test that has high power against, both, sparse and dense alternatives?} (Section~\ref{subsec:PowerEnhancement})
	\item The validity of bootstrapping procedures (in, both, classical as well as modern high-dimensional settings) hinges on moment conditions~\citep{gine1990bootstrapping, giessing2023Gaussian}.  \emph{Can we improve the theoretical and practical performance of the bootstrap test by leveraging distributional assumptions commonly used in high-dimensional statistics?} (Section~\ref{subsec:EllipticallyDistributedData})
	\item The proposed bootstrap test depends on the data via an estimate of the high-dimensional covariance matrix. Estimating high-dimensional covariance matrices can be challenging~\citep{cai2010optimal, cai2011adaptive, avella-medina2018robust}. \emph{How does the estimated covariance matrix affect the validity of the bootstrap test?} (Section~\ref{subsec:HeuristicsBias})
	\item The proposed test statistic is non-pivotal; in particular, the sampling distribution of the test statistic depends on the unknown population covariance matrix. In practice, when bootstrapping non-pivotal test statistics the actual level of the test often differs substantially from the nominal one. Several (computationally expensive) bias correction schemes have been proposed in classical settings~\citep{davison1986bootstrap, beran1987prepivoting, hall1988resampling, shi1992accurate}. \emph{Can we draw on insights from high-dimensional probability theory to develop new (computationally more efficient) correction schemes?} (Sections~\ref{subsec:CovarianceEstimation}-\ref{subsec:SphericalBootstrap})
\end{itemize}
The reader interested in only the answers to these questions and related open problems, may directly jump to the discussion in Section~\ref{sec:Discussion}.

\subsection{A bootstrap hypothesis test based on $\ell_p$-norms}\label{subsec:BootstrapHypothesisTest}
We propose to test hypothesis~\eqref{eq:sec:Intro-1} on the basis of the $\ell_p$-statistic
\begin{align*} 
	T_{n,p}: = \|RS_n - \sqrt{n} r\|_p, \quad{}\quad{} S_n = \frac{1}{\sqrt{n}} \sum_{i=1}^n X_i,
	\quad{}\quad{} p \geq 1,
\end{align*}
and, given a nominal level $\alpha \in (0,1)$, reject the null hypothesis if and only if
\begin{align*} 
	T_{n,p} \geq c^*_p(1- \alpha;\widehat{\Omega}_n),
\end{align*}
where $c^*_p(\alpha;\widehat{\Omega}_n ) = \inf\left\{t \in \mathbb{R} : \mathrm{P}(T_{n,p}^* \leq t \mid X_1, \ldots, X_n)  \geq \alpha \right\}$
is the conditional $\alpha$-quantile of the Gaussian proxy statistic
\begin{align*}
	T_{n,p}^* := \|Z_n\|_p, \quad{}\quad{} Z_n \mid X_1, \ldots, X_n \sim N(0, \widehat{\Omega}_n),
\end{align*}
and $\widehat{\Omega}_n$ is a positive semi-definite estimate of $\Omega = R \Sigma R'$. We call above testing procedure the \emph{bootstrap test based on $T_{n,p}$ at level $\alpha$} and define
\begin{align*}
	\varphi_\alpha(T_{n,p},\widehat{\Omega}_n) := \mathbf{1}\left\{T_{n,p} \geq c^*_p(1- \alpha;\widehat{\Omega}_n)\right\}.
\end{align*}
The bootstrap test based on $T_{n,p}$ rejects $H_0$ at level $\alpha$ if and only if $\varphi_\alpha(T_{n,p},\widehat{\Omega}_n) = 1$.

To understand the rationale for the test, consider the low-dimensional case first: If dimensions $d,t$ are fixed and the data are i.i.d. with finite second moments, the CLT and the continuous mapping theorem imply that $T_{n,p}\overset{d}{\rightarrow} \|Z\|_p$, where $Z \sim N\left(0, \Omega\right)$. Hence, in this case, the Gaussian proxy statistic $T_{n,p}^*$ is just the parametric bootstrap estimate of the limiting random variable $\|Z\|_p$. Of course, if $d, t \geq \sqrt{n}$, then the CLT does not apply and the limiting random variable $Z$ needs not to exist. The reason why the bootstrap is nonetheless valid in high dimensions is that a CLT is not needed. Indeed, validity of the bootstrap test follows already from convergence of the laws of $T_{n,p}$ and $T_{n,p}^*$ with respect to the Kolmogorov distance (see Appendix~\ref{sec:TheoryBootstrap}). Since two sequences of laws are close in Kolmogorov distance if they have the same cluster points, this does not imply (or necessitate) existence of a limit law.

In the remainder of the paper we will study the asymptotic properties of this test, develop modifications for special scenarios, and provide answers to the questions raised in Section~\ref{subsec:GlobalNull}. 



\subsection{Contributions and outline of the paper}
In detail, our \emph{contributions to statistical theory} are as follows: First, we show that the bootstrap tests $\varphi_\alpha(T_{n,p},\widehat{\Omega}_n)$ have asymptotic correct size for all $p \in [1, \infty]$. By inverting the test statistics we then obtain Scheff{\'e}-type simultaneous confidence regions for multiple testing problems 
(Section~\ref{subsec:AsympSizeLevel}). Second, we establish that the bootstrap tests are asymptotically unbiased and consistent against certain alternatives (Section~\ref{subsec:AsympCorrectness}). In particular, we show that for each exponent $p \in [1, \infty]$ there exist at least three different classes of alternatives characterized by their signal strength relative to noise level and expected size of the Gaussian proxy statistic $T_{n,p}^*$. 
Third, to assess the power of the bootstrap tests (against sparse and dense alternatives and for different exponents $p \geq 1$) we adapt the classical concept of the Bahadur slope to the high-dimensional setting. 
We then show that $\varphi_\alpha(T_{n,2},\widehat{\Omega}_n)$ is most powerful among all bootstrap tests when testing against dense alternatives whose signals are spread out over a large number of coordinates whereas $\varphi_\alpha(T_{n,\log t},\widehat{\Omega}_n)$ is most powerful among all bootstrap tests when testing against sparse alternatives with only a few strong signals
(Section~\ref{subsec:BahadurSlope}). Fourth, we discuss three modifications of the basic bootstrap tests: We construct a bootstrap test with power against, both, sparse and dense alternatives (Section~\ref{subsec:PowerEnhancement}), propose a self-normalized test for elliptically (heavy-tailed) data (Section~\ref{subsec:EllipticallyDistributedData}), and briefly discuss bootstrap tests for population parameters other than the mean (Section~\ref{subsec:ApproximateSampleAverages}). 

The key \emph{methodological innovations} of this paper are two new bias correction schemes to improve the accuracy of the bootstrap test. Typically, when implementing bootstrap methods via repeated sampling from the empirical distribution or, as in our case, repeated sampling from the Gaussian proxy statistic $T_{n,p}^*$ the actual level of the test differs from the nominal level $\alpha$ even if the number of samples $B$ is large. We develop a three-fold bias decomposition (Section~\ref{subsec:HeuristicsBias}) and propose methods to mitigate two of the biases: First, we show how to incorporate structural information of the high-dimensional covariance matrix of the data into the bootstrap test (Section~\ref{subsec:CovarianceEstimation}). Second, we propose a novel spherical multiplier bootstrap test which exploits hypercontractivity properties of $\ell_p$-norms of high-dimensional random vectors (Section~\ref{subsec:SphericalBootstrap}). We provide supporting results from Monte Carlo experiments in Section~\ref{sec:NumericalStudies}. 

The main \emph{probabilistic contribution} of this paper is that the bootstrap test $\varphi_\alpha(T_{n,p},\widehat{\Omega}_n)$ is indeed valid \emph{for all $p \geq 1$} even if dimensions $d,t$ far exceed the sample size $n$. We develop the necessary technical tools in Appendix~\ref{sec:TheoryBootstrap}. These include, among other things, Gaussian approximation and anti-concentration inequalities for high-dimensional $\ell_p$-statistics (Appendix~\ref{subsec:GeneralResults}) and lower bounds on the variances of $\ell_p$-norms of (an-isotropic) Gaussian random vectors (Appendix~\ref{subsec:LowerBoundVariance}). For the reader's convenience, we also provide a brief exposition of the relevant results from the two companion papers~\cite{giessing2022anticoncentration, giessing2023Gaussian} (Appendix~\ref{subsec:ResultsEPT}).

\subsection{Notation}
We denote the $\ell_p$-norm, $p \in [1, \infty)$, and the $\ell_\infty$-norm of a vector $x \in \mathbb{R}^d$ by $\|x\|_p = (\sum_{k=1}^d |x_k|^p)^{1/p}$ and $\|x\|_\infty = \max_{1 \leq k \leq d}|x_k|$, respectively. For a symmetric matrix $A$, we denote its operator, Frobenius, and $q,p$-matrix norm by $\|A\|_{op}$, $\|A\|_F$, and $\|A\|_{q \rightarrow p} := \sup_{\|u\|_q=1}\|Au\|_p$, $1 \leq p, q \leq \infty$. For two deterministic sequences $(a_n)_{n\geq 1}$ and $(b_n)_{n\geq 1}$, we write $a_n \lesssim b_n$ if $a_n = o(b_n)$. Also, $a_n \asymp b_n$ if there exist absolute constants $C_1, C_2 > 0$ such that $C_1 b_n \leq a_n \leq C_2 b_n$ for all $n \geq 1$. We denote the set of mean vectors which are equivalent under the null hypothesis by $\mathcal{H}_0 = \{\mu \in \mathbb{R}^d : R\mu = r\}$ and the set of those belonging to the alternative hypothesis by $\mathcal{H}_1 = \mathcal{H}_0^c$.

\section{Properties of the bootstrap hypothesis test}\label{sec:BootstrapTest}


\subsection{High-level assumptions and simple sufficient conditions}\label{subsec:Assumptions}
To carry out the theoretical analysis of the bootstrap tests $\varphi_\alpha(T_{n,p},\widehat{\Omega}_n)$ we need a set of minimal assumption which we present in this section (Assumptions~\ref{assumption:ControlThirdMoments}--\ref{assumption:ConsistentCovariance}). We include three examples (Lemmas~\ref{lemma:example:SubGaussianData}--\ref{lemma:example:HeavyTailedData}) which illustrate these minimal assumptions for concrete data generating processes.

Throughout, $X, X_1, \ldots, X_n \in \mathbb{R}^d$ are i.i.d. random vectors with mean $\mu_n \in \mathbb{R}^d$ and positive semi-definite covariance matrix $\Sigma \neq \mathbf{0} \in \mathbb{R}^{d \times d}$, $Z \sim N(0, \Omega)$ with $\Omega \equiv R \Sigma R' \neq \mathbf{0} \in \mathbb{R}^{t \times t}$ and $R \in \mathbb{R}^{t \times d}$. We allow the dimension $d$ and the number of restrictions $t$ to vary as a function of the sample size $n$; in particular, it is possible that $d, t \rightarrow \infty$ as $n \rightarrow \infty$.

The high-level assumptions that we impose are conditions on the moments of the data. These moment conditions originate from the technical results concerning the approximation of the sampling distribution of $T_{n,p}$ (see Appendix~\ref{sec:TheoryBootstrap}) and implicitly restrict the growth rates of dimensions $d, t$ relative to the sample size $n$. 

\begin{assumption}[Control of third moments]\label{assumption:ControlThirdMoments}
The $\ell_p$- and Gaussian proxy statistic satisfy
\begin{align*}
	\left(\mathrm{E}[\|R(X-\mu_n)\|_p^3]\right)^{1/3} \vee \left(\mathrm{E}[\|Z\|_p^3]\right)^{1/3}   \lesssim n^{1/6}\mathrm{Var}(\|Z\|_p)^{1/2}.
\end{align*}
\end{assumption}

\begin{assumption}[Control of ratio of moments]\label{assumption:RatioMoments}
There exists $s > 3$ such that 
\begin{align*}
	\frac{\left(\mathrm{E}[\|R(X-\mu_n)\|_p^s]\right)^{1/s} }{\left(\mathrm{E}[\|R(X-\mu_n)\|_p^3]\right)^{1/3}} \lesssim n^{1/3-1/s}.
\end{align*}
\end{assumption}

\begin{assumption}[Consistent estimate of covariance matrix]\label{assumption:ConsistentCovariance} 
	There exists a nonparametric statistic $\widehat{\Omega}_n$ based on the $X_i$'s only such that
	\begin{align*}
		\|\widehat{\Omega}_n- \Omega\|_{q \rightarrow p} = o_p\big( \mathrm{Var}(\|Z\|_p) \big),  \quad{}\quad{} 1/p + 1/q = 1,
	\end{align*}
	where $\|M\|_{q \rightarrow p} := \sup_{\|u\|_q=1}\|Mu\|_p$ for matrices $M$.
\end{assumption}

\begin{remark}[On Assumption~\ref{assumption:ConsistentCovariance}]\label{remark:assumption:ConsistentCovariance}
	We require $\widehat{\Omega}_n$ to be a nonparametric estimator to guarantee that the Gaussian proxy statistic $T_{n,p}^*$ is ancillary. Ancilliarity of $T_{n,p}^*$ is necessary for computing critical values for composite hypotheses and confidence regions. If we were interested in testing only one null and alternative hypothesis at a time, we could do without it, e.g. we could use $\widehat{\Omega}_n = n^{-1} \sum_{i=1}^n R(X_i - \mu)(X_i - \mu)'R'$, where $\mu \in \mathbb{R}^d$ depends on the (null) hypothesis under consideration.
\end{remark}

These assumptions are mild and can easily be verified to hold in a range of scenarios. In the following we give three examples with simple sufficient conditions. First, we consider data with sub-Gaussian tails which is a common assumption in high dimensional statistics. In the second example we introduce data with log-concave distribution which is a useful assumption for proving several subtle properties of the bootstrap test in Sections~\ref{subsec:AsympCorrectness} and~\ref{subsec:BahadurSlope}. The third example is about heavy-tailed data; it is important in applications and becomes relevant in Section~\ref{subsec:EllipticallyDistributedData} where we introduce a modified test statistic for elliptically distributed data which allows us to further relax Assumptions~\ref{assumption:ControlThirdMoments}--\ref{assumption:ConsistentCovariance}. To keep the exposition simple, we only consider exponents $p \in \{2, \infty\}$. However, in principle, all three examples can be extended to all $p \in [1, \infty]$. In the following, $\widehat{\Omega}_n = n^{-1}\sum_{i=1}^n R(X_i - \bar{X}_n)(X_i - \bar{X}_n)'R'$ denotes the (nonparametric) sample covariance matrix of the transformed data $RX_1, \ldots, RX_n \in \mathbb{R}^t$.

\begin{lemma}[Sub-Gaussian design]\label{lemma:example:SubGaussianData}
	Let $X, X_1, \ldots X_n \in \mathbb{R}^d$ be a simple random sample of sub-Gaussian random vectors with mean $\mu_n \in \mathbb{R}^d$ and positive semi-definite covariance matrix $\Sigma \in \mathrm{R}^{d \times d}$ such that $\| (X - \mu_n)'u\|_{\psi_2} \lesssim \|\Sigma^{1/2}u\|_2$ for all $u \in \mathbb{R}^d$. Then, Assumptions~\ref{assumption:ControlThirdMoments}--\ref{assumption:ConsistentCovariance} are satisfied
	\begin{itemize}
		\item[(i)] for $p=2$ if 
		\begin{align*}
			\frac{r(\Omega)}{\sqrt{r(\Omega^2)}} = o(n^{1/6}) \quad{} \quad{} \mathrm{and}\quad{}\quad{} \frac{r(\Omega)}{r(\Omega^2)}  \left( \sqrt{\frac{r(\Omega)}{n}} \vee \frac{r(\Omega)}{n} \right) = o(1),
		\end{align*}
		where $r(M) := \mathrm{tr}(M)/\|M\|_{op}$ is the effective rank of $M \in \{\Omega, \Omega^2\}$;
		\item[(ii)] for $p = \infty$ if
		\begin{align*}
		\frac{\omega_{(t)}}{\omega_{(1)}} \frac{\sqrt{\log t} }{n^{1/6}} \vee  \frac{\omega_{(t)}^2}{\omega_{(1)}^2} \frac{ \log t}{n^{1/6}} = o(1) \quad{}\quad{} \mathrm{and}\quad{}\quad{}\frac{\omega_{(t)}^2 }{\omega_{(1)}^2} \sqrt{\frac{\log t }{n}} \vee \frac{\omega_{(t)}^4 }{\omega_{(1)}^4} \sqrt{\frac{(\log t)^3}{n}} = o(1),
		\end{align*}
	\end{itemize}
	where $\omega_{(1)}^2 \leq \ldots \leq \omega_{(t)}^2$ are the ordered diagonal elements of $\Omega$.	
\end{lemma}

\begin{remark}[On the exponent $p$ and covariance matrix $\Omega$]
	The conditions in Lemma~\ref{lemma:example:SubGaussianData} (i) on $\Omega$ are easily satisfied for matrices with variance decay or bounded effective rank such that $\mathrm{r}(\Omega), \mathrm{r}(\Omega^2) < \infty$ independent of or slowly increasing in the dimensions $d,t$. We refer to~\cite{lopes2020bootstrapping} for more explicit examples. In Lemma~\ref{lemma:example:SubGaussianData} (ii) the ratio $\omega_{(t)}^2/ \omega_{(1)}^2$ can be upper bounded by the conditioning number of $\Omega$, i.e. ratio of its largest to smallest eigenvalue.
\end{remark}

\begin{lemma}[Log-concave design]\label{lemma:example:LogConcaveData}
	Let $X, X_1, \ldots X_n \in \mathbb{R}^d$ be a simple random sample of random vectors with log-concave density, i.e. the density $f$ satisfies $f = e^{-\varphi}$ with $\varphi$ convex. Then, Assumptions~\ref{assumption:ControlThirdMoments}--\ref{assumption:ConsistentCovariance} are satisfied
	\begin{itemize}
		\item[(i)] for $p=2$ if
		\begin{align*}
			\frac{r(\Omega)}{\sqrt{r(\Omega^2)}} = o(n^{1/6}) \quad{} \quad{} \mathrm{and}\quad{}\quad{} \frac{r(\Omega)}{r(\Omega^2)} \left(\left(\frac{r(\Omega)}{n}\right)^{1/4} \vee \frac{r(\Omega)}{n} \right) = o(1),
		\end{align*}
		where $r(M) := \mathrm{tr}(M)/\|M\|_{op}$ is the effective rank of $M \in \{\Omega, \Omega^2\}$;
		\item[(ii)] for $p = \infty$ if 
		\begin{align*}
			\frac{\omega_{(t)}}{\omega_{(1)} } \frac{\log (t K_\varphi)}{n^{1/6}} \vee  \frac{\omega_{(t)}^2}{\omega_{(1)}^2} \frac{\sqrt{\log (tK_\varphi)^3}}{n^{1/6}} = o(1) \quad{} \mathrm{and}\quad{}\frac{\omega_{(t)}^2 }{\omega_{(1)}^2} \left(\frac{(\log t)^5}{n}\right)^{1/4} \vee \frac{\omega_{(t)}^4}{\omega_{(1)}^4}\sqrt{\frac{(\log t)^3}{n}} = o(1),
		\end{align*}
	\end{itemize}
	where $\omega_{(1)}^2 \leq \ldots \leq \omega_{(t)}^2$ are the ordered diagonal elements of $\Omega$ and $K_\varphi > 1$ is a constant depending on (the marginals of) the density $f = e^{-\varphi}$.
\end{lemma}
\begin{remark}[On the convexity of $\varphi$]
If $\varphi$ is strictly convex in the sense that $\varphi'' \geq \lambda I_d$ for some $\lambda > 0$, then the $X_i$'s are sub-Gaussian and the bounds of Lemma~\ref{lemma:example:SubGaussianData} apply. For details (in particular, how the sub-Gaussian tail bounds depend on $\lambda >0$) we refer to Theorems 5.2 and 5.3 in~\cite{ledoux2001concentration}.
\end{remark}

\begin{lemma}[Heavy-tailed design]\label{lemma:example:HeavyTailedData}
	Let $X, X_1, \ldots X_n \in \mathbb{R}^d$ be a simple random sample of random vectors. Then, Assumptions~\ref{assumption:ControlThirdMoments}--\ref{assumption:ConsistentCovariance} are satisfied
\begin{itemize}
	\item[(i)] for $p=2$ if there exists $s > 3$,
	\begin{align*}
		m_{2, 3} \sqrt{\frac{r(\Omega)}{r(\Omega^2)}} = o(n^{1/6}), \quad{} \frac{m_{2, s}}{m_{2,3}} = o(n^{1/3-1/s}), \quad{} \frac{r(\Omega)}{r(\Omega^2)} \left(\sqrt{ \frac{ M_n^2 (\log n)}{n}}  \vee  \frac{ M_n^2 (\log n)}{n} \right) = o(1),
	\end{align*}
	where $r(M) := \mathrm{tr}(M)/\|M\|_{op}$ is the effective rank of $M \in \{\Omega, \Omega^2\}$, $M_n^2:= \mathbb{E}[\max_{1 \leq i \leq n}\|X_i-\mu_n\|_2^2] / \|\Omega\|_{op}$, and $m_{2,s}^s := \mathbb{E}[\|X-\mu_n\|_2^s] / \|\Omega\|_{op}^{s/2}$;
	\item[(ii)] for $p = \infty$ if there exists $s > 3$, 
	\begin{align*}
	&m_{\infty,3}\left(\frac{\omega_{(t)}}{\omega_{(1)}} \vee \frac{\omega_{(t)}^2}{\omega_{(1)}^2}\sqrt{\frac{\log t }{n}} \right)  = o(n^{1/6}), \quad{} \frac{m_{\infty, s}}{m_{\infty,3}} = o(n^{1/3-1/s}),\\
	&\frac{\omega_{(t)}^2 }{\omega_{(1)}^2} \left(\frac{(\log t)^5}{n}\right)^{1/4} \vee \frac{\omega_{(t)}^4}{\omega_{(1)}^4}\sqrt{\frac{(\log t)^3}{n}} = o(1),
	\end{align*}
\end{itemize}
where $\omega_{(1)}^2 \leq \ldots \leq \omega_{(t)}^2$ are the ordered diagonal elements of $\Omega$ and $m_{\infty, s}^s := \mathbb{E}[\|X-\mu_n\|_\infty^s] / \omega_{(t)}^s$.
\end{lemma}

\begin{remark}[On the $(\log n)$-factor]
 If the $\|X_i\|_2 \leq M_n$ almost surely, we can remove the $(\log n)$-factors in Lemma~\ref{lemma:example:HeavyTailedData} (i). Using an alternative argument based on the proof of Lemma~\ref{lemma:example:LogConcaveData}, we can replace the third relation in Lemma~\ref{lemma:example:HeavyTailedData} (i) by $\frac{r(\Omega)}{r(\Omega^2)}  \Big( \big(\frac{r(\Omega)}{n}\big)^{1/4} \vee \frac{r(\Omega)}{n} \Big) = o(1)$.
\end{remark}

\subsection{Asymptotic size}\label{subsec:AsympSizeLevel}
As a first result we show that the bootstrap test $\varphi_\alpha(T_{n,p},\widehat{\Omega}_n)$ has asymptotic correct size, i.e. asymptotically the bootstrap test controls the type 1 error. By inverting the hypothesis tests we then obtain asymptotic $1-\alpha$ confidence sets and Scheff{\'e}-type simultaneous confidence intervals for multiple comparison problems.

\begin{theorem}[Asymptotic size $\alpha$ test]\label{theorem:SizeAlphaTest}
	Suppose that Assumptions~\ref{assumption:ControlThirdMoments}--\ref{assumption:ConsistentCovariance} hold. Then, for $1 \leq p \leq \infty$,
	\begin{align*}
		\lim_{n \rightarrow \infty}\sup_{\alpha \in (0,1)} \sup_{\mu \in \mathcal{H}_0}\Big|\mathrm{E}_\mu[\varphi_\alpha(T_{n,p}, \widehat{\Omega}_n) ] - \alpha \Big| =0.
	\end{align*}
\end{theorem}
Since $T_{n,p} = \|RS_n - \sqrt{n}r\|_p = \sup_{\|u\|_q = 1} u'(RS_n - \sqrt{n}r)$ for $1/q + 1/p = 1$, we immediately obtain the following corollary:

\begin{corollary}\label{corollary:theorem:SizeAlphaTest-0}
	Recall the conditions of Theorem~\ref{theorem:SizeAlphaTest}. Let $u \in \mathbb{R}^t$ be a (possibly) random vector not necessarily independent of $X_1, \ldots, X_n$.
	Then, for $U_{n,p} :=  u'(RS_n - \sqrt{n}r)/\|u\|_q$ with $1/p + 1/q = 1$,
	\begin{align*}
		\lim_{n \rightarrow \infty} \sup_{\mu \in \mathcal{H}_0} \mathrm{E}_\mu[	\varphi_\alpha(U_{n,p}, \widehat{\Omega}_n) ] \leq \alpha.
	\end{align*}
\end{corollary}
Corollary~\ref{corollary:theorem:SizeAlphaTest-0} is of some interest because it implies that we can use the bootstrap test $\varphi_\alpha(U_{n,p}, \widehat{\Omega}_n)$  to test (simple) data dependent hypotheses while still controlling the type 1 error:

\begin{example}[Post-selective test of subsets of a population mean vector]
	Let $S_{n, (1)} \leq \ldots \leq S_{n, (d)}$ be the absolute values of the entries of $S_n \in \mathbb{R}^d$ arranged in ascending order. Let $\widehat{J} = \{i : S_{n, (i)} \geq S_{n, (B)}\} \subset \{1, \ldots, d\}$ be the random subset of the $B \geq 1$ largest entries in $S_n$ and consider the data dependent hypotheses $H_0: \mu_{0,\widehat{J}} = 0$ vs. $H_1: \mu_{0,\widehat{J}} \neq 0$. Since $\|S_{n, \widehat{J}}\|_p \leq \sup_{\|u\|_q =1} u'S_n = \|S_n\|_p$ we can test this hypothesis using above Corollary with $R= I$ and $r = 0$~\citep[such tests are relevant in the context of exploratory gene expression/ gene set/ pathway analyses, e.g.][]{matur2018gene}. Since this test has a high threshold for rejecting the null hypothesis (i.e. the $(1-\alpha)$ quantile of the Gaussian proxy statistic $T_{n,p}^* = \|Z_n\|_p = \sup_{\|u\|_q =1} u'Z_n$), it is conservative and under-powered. It is not a panacea to all post-selective inference problems.
	
	A test with higher power is the one based on the $(1-\alpha)$ quantile of the Gaussian proxy statistic $\sup_{J \subset \{1, \ldots, d\}, \: |J| = B}\|Z_{n,J}\|_p$. Intuitively, this proxy statistic captures the maximum spurious signal of $B$ variables. The set $\widehat{J}$ is considered significant only if its signal is stronger than the  maximum spurious signal of $B$ variables. Under properly adjusted Assumptions~\ref{assumption:ControlThirdMoments}--\ref{assumption:ConsistentCovariance}, it is trivial to modify the proofs of Theorem~\ref{theorem:SizeAlphaTest} and Corollary~\ref{corollary:theorem:SizeAlphaTest-0} to verify that this test has size at most $\alpha$, too.
\end{example}

Another easy consequence of Theorem~\ref{theorem:SizeAlphaTest} is the following:
\begin{corollary}\label{corollary:theorem:SizeAlphaTest}
	Recall the conditions of Theorem~\ref{theorem:SizeAlphaTest}. Suppose that in addition there exists a random matrix $\widehat{R}_n \in \mathbb{R}^{t \times d}$, not necessarily independent of $X_1, \ldots, X_n$, such that
	\begin{align}\label{eq:theorem:SizeAlphaTesttheorem:SizeAlphaTest-0}
		\|\widehat{R}_n - R\|_{q \rightarrow p} \|S_n - \sqrt{n} \mu_0\|_p = o_p\left(\sqrt{\mathrm{Var}(\|Z\|_p)}\right), \quad{}\quad{} 1/p + 1/q = 1.
	\end{align}
	Then, for $\widetilde{T}_{n,p} := \|\widehat{R}_n (S_n - \sqrt{n} \mu_0)\|_p$ and all $1 \leq p \leq \infty$,
	\begin{align*}
		\lim_{n \rightarrow \infty}\sup_{\alpha \in (0,1)} \sup_{\mu \in \mathcal{H}_0}\Big|\mathrm{E}_\mu[	\varphi_\alpha(\widetilde{T}_{n,p}, \widehat{\Omega}_n) ]  - \alpha \Big| = 0.
	\end{align*}
\end{corollary}
Corollary~\ref{corollary:theorem:SizeAlphaTest} shows that (simple) preprocessing of the data does not affect the size of the bootstrap test. A common preprocessing step is to studentize the data. Studentization is particularly appropriate if the data set has been obtained by aggregating data from various sources with different noise levels or is known to be heteroscedastic. Moreover, in the classical setting with fixed dimensions, (bootstrap) tests based on studentized statistics are known to be more accurate~\citep{hall1986bootstrapCI}. We have evidence from simulations that the same holds true in high dimensions.

\begin{example}[Tests based on studentized statistics]\label{example:StudentizedStatistics}
	Suppose that Assumptions~\ref{assumption:ControlThirdMoments}--\ref{assumption:ConsistentCovariance} hold with $R^2 = \mathrm{diag}(\Sigma)^{-1}$ and that eq.~\eqref{eq:theorem:SizeAlphaTesttheorem:SizeAlphaTest-0} is satisfied with $\widehat{R}_n^2 = \mathrm{diag}(\widehat{\Sigma}_n)^{-1}$, where $\widehat{\Sigma}_n = (\hat{\sigma}_{kj})_{k,j=1}^d$ is the sample covariance matrix based on the $X_i$'s. By Corollary~\ref{corollary:theorem:SizeAlphaTest} we can test $H_0 : \mu =\mu_0$ versus $H_1: \mu \neq \mu_0$ at level $\alpha$ via $\varphi_\alpha(\widetilde{T}_{n,p}, \widehat{\Omega}_n)$, where $\widetilde{T}_{n,p} := \|\widetilde{S}_n - \tilde{\mu}_0\|_p$ is the studentized $\ell_p$-statistic with $\tilde{\mu}_0 = (\tilde{\mu}_{0,k}/ \hat{\sigma}_{kk})_{k =1}^d$, $\widetilde{X}_i = (X_{ik}/\hat{\sigma}_{kk})_{k=1}^d$, and $\widetilde{S}_n = n^{-1/2}\sum_{i=1}^n \widetilde{X}_i$. Notice that $\|\widehat{R}_n - R\|_{q \rightarrow p} \asymp  \max_{1 \leq k \leq d}|\hat{\sigma}_{kk}^2 - \sigma_{kk}^2|$ for all $p \in \{2\} \cup [\log d, \infty]$. Thus, for these exponents, eq.~\eqref{eq:theorem:SizeAlphaTesttheorem:SizeAlphaTest-0} is easily verified by combining the results in Section~\ref{subsec:LowerBoundVariance} with the commonly used moment conditions in high-dimensional statistics. 
\end{example}

By inverting the test $\varphi_\alpha(T_{n,p}, \widehat{\Omega}_n)$ we obtain asymptotic $1- \alpha$ bootstrap confidence regions and Scheff{\'e}-type simultaneous bootstrap confidence intervals. The following theorem is essentially a restatement of Theorem~\ref{theorem:SizeAlphaTest} and therefore does not need a proof.

\begin{proposition}[Asymptotic confidence regions and simultaneous confidence intervals]\label{theorem:ConfidenceSets} 
	Suppose that Assumptions~\ref{assumption:ControlThirdMoments}--\ref{assumption:ConsistentCovariance} hold. Set $\hat{\mu}_n = n^{-1/2} S_n = n^{-1} \sum_{i=1}^n X_i$. Then, for all $\alpha \in (0,1)$ and $1 \leq p, q \leq \infty$ such that $1/p + 1/q =1$ the following holds.
	\begin{itemize}
		\item[(i)] The $\ell_p$-norm ellipsoid 
		\begin{align*}
			\mathcal{E}_{n,p}(R, \widehat{\Omega}_n) : = \left\{ \mu \in \mathbb{R}^d : \sqrt{n}\|R(\hat{\mu}_n - \mu)\|_p \leq c^*_p(1-\alpha;\widehat{\Omega}_n)\right\}
		\end{align*}
		is an asymptotic $(1-\alpha)$-confidence region for $\mu_0$, $\lim_{n \rightarrow \infty} \mathrm{P}\left(\mu_0 \in \mathcal{E}_{n,p}(R, \widehat{\Omega}_n) \right) = 1- \alpha$.
		\item[(ii)]	The intervals
		\begin{align*}
			\mathcal{I}_{n,p}(h;R, \widehat{\Omega}_n) : = \left[ h' R\hat{\mu}_n -\frac{\|h\|_q}{\sqrt{n}} c^*_p(1-\alpha;\widehat{\Omega}_n), \: h'R\hat{\mu}_n+\frac{\|h\|_q}{\sqrt{n}} c^*_p(1-\alpha;R,\widehat{\Omega}_n)  \right],\quad{} h \in \mathbb{R}^d,
		\end{align*}
		are asymptotic simultaneous $(1-\alpha)$-confidence intervals for all linear combinations $h'R\mu_0$, $h \in \mathbb{R}^d$, i.e. $\lim_{n \rightarrow \infty} \mathrm{P} \left(h'R\mu_0 \in \mathcal{I}_{n,p}(h;R, \widehat{\Omega}_n), \: \forall h \in \mathbb{R}^d \right) = 1- \alpha$.
	\end{itemize}
\end{proposition}
\begin{remark}[On the matrix of linear restrictions $R$]
	Under the conditions of Corollary~\ref{corollary:theorem:SizeAlphaTest} the deterministic matrix $R \in \mathbb{R}^{t \times d}$ can be replaced by a random matrix $\widehat{R}_n$. Thus, analogous to Corollary~\ref{corollary:theorem:SizeAlphaTest}, we can construct confidence regions nd simultaneous confidence intervals for the studentized mean vector.
\end{remark}

Apart from asymptotic exact $1- \alpha$ confidence sets for the mean $\mu_0$ and linear combinations $h'\mu_0$ Proposition~\ref{theorem:ConfidenceSets} also yields asymptotic level $1- \alpha$ confidence intervals for the norm $\|\mu_0\|_p$:

\begin{example}[Conservative confidence intervals for $\|\mu_0\|_p$]\label{example-1}
	Suppose that Assumptions~\ref{assumption:ControlThirdMoments}--\ref{assumption:ConsistentCovariance} hold with $R = I_d  \in \mathbb{R}^{d\times d}$. Let $\widehat{\Sigma}_n = (\hat{\sigma}_{kj})_{k,j=1}^d$ be the sample covariance matrix based on the $X_i$'s. By the triangle inequality we have $\big| \|\hat{\mu}_n\|_p - \|\mu_0\|_p \big| \leq \|\hat{\mu}_n - \mu_0\|_p$.  Hence, by Proposition~\ref{theorem:ConfidenceSets} (i)
	\begin{align*}
		\mathcal{C}_{n,p} : = \left[ \|\hat{\mu}_n\|_p - n^{-1/2} c^*_p(1-\alpha; \widehat{\Sigma}_n), \:\|\hat{\mu}_n\|_p + n^{-1/2} c^*_p(1-\alpha; \widehat{\Sigma}_n ) \right], \quad{} p \geq 1,
	\end{align*}
	is an asymptotic level $1-\alpha$ interval for $\|\mu_0\|_p$, i.e. $\lim_{n \rightarrow \infty} \mathrm{P}\left(\|\mu_0\|_p \in \mathcal{C}_{n,p} \right) \geq 1- \alpha$.
\end{example}

If one is interested in inference about a single linear combination of a high-dimensional mean vector, say $x'\mu_0$, it is advisable to use Proposition~\ref{theorem:ConfidenceSets} (i) with $R = x'$ (not Proposition~\ref{theorem:ConfidenceSets} (ii) with $h =x$ and $R = I_d$) in order to obtain intervals with short width:

\begin{example}[Confidence interval for a single linear combination of a high-dimensional mean vector]\label{example-3}
	Suppose that Assumptions~\ref{assumption:ControlThirdMoments}--\ref{assumption:ConsistentCovariance} hold with $R = x'$. Let $\widehat{\Omega}_x$ be an estimate of the scalar variance $x'\Sigma x$. Then, by Proposition~\ref{theorem:ConfidenceSets} (i),
	\begin{align*}
		\left[x'\hat{\mu}_n - n^{-1/2}  c^*_p(1-\alpha; \widehat{\Omega}_x ), \: x'\hat{\mu}_n + n^{-1/2}c^*_p(1-\alpha; \widehat{\Omega}_x )  \right]
	\end{align*}
	is an asymptotic $(1-\alpha)$-confidence interval for $x'\mu_0$. Notice that this result holds without sparsity assumptions of $x$ or $\mu_0$.
\end{example}

\subsection{Asymptotic unbiasedness and consistency}\label{subsec:AsympCorrectness}
Next, we establish unbiasedness and consistency of the bootstrap test $\varphi_\alpha(T_{n,p}, \widehat{\Omega}_n)$. Both properties turn out to be more nuanced in high dimensions than the classical setting with fixed dimensions. We conclude the section with a discussion of the effect of the exponent $p \in [1, \infty]$ on the consistency of $\varphi_\alpha(T_{n,p}, \widehat{\Omega}_n)$ against sparse and dense alternatives.

The following theorem shows that (asymptotically) no alternative in $\mathcal{H}_1$ has probability of rejection less than the size of the bootstrap test. This means that we accept the alternative with higher probability when it is correct than when it is false.

\begin{theorem}[Unbiasedness]\label{theorem:UnbiasedTest} 
	Suppose that Assumptions~\ref{assumption:ControlThirdMoments}--\ref{assumption:ConsistentCovariance} hold and, in addition, the $X_i$'s are symmetric around their mean $\mu_n$ and have a log-concave density, i.e. the density satisfies $f = e^{-\varphi}$ with $\varphi$ convex and $\varphi(\mu_n - x) = \varphi(\mu_n + x)$ for all $x \in \mathbb{R}^d$. Then, for all $1 \leq p \leq \infty$ and all $\alpha \in (0,1)$,
	\begin{align*}
		\lim_{n \rightarrow \infty} \mathrm{E}_{\mu}[\varphi_\alpha(T_{n,p}, \widehat{\Omega}_n)] \geq \alpha \quad \quad \forall \mu \in \mathcal{H}_1.
	\end{align*}
\end{theorem}
\begin{remark}[On the assumptions on the density]
	In low dimensions, when $S_n = n^{-1/2} \sum_{i=1}^n R(X_i-\mu_n)$ converges weakly under the alternative hypothesis to a Gaussian random vector, the additional condition on the density of the $X_i$'s can be dropped since the density of the limiting distribution $N(0, \Omega)$ is automatically symmetric and log-concave. The situation is more complex in high dimensions when $S_n$ does not converge to a Gaussian random vector and is not even Gaussian approximable under the alternative hypothesis. In this case, symmetry and log-concavity of the $X_i$'s allow us to invoke Anderson's lemma and reduce the problem to Gaussian approximation under the null hypothesis which falls within the scope of Theorem~\ref{theorem:Bootstrap-Lp-Norm-Quantiles}.
\end{remark}

It is straightforward to verify that in the classical setting with fixed dimensions $d,t < \infty$ the bootstrap test $\varphi_\alpha(T_{n,p}, \widehat{\Omega}_n)$ is asymptotically consistent against any alternative $\mu \in \mathcal{H}_1$. In high dimensions with $d, t \rightarrow \infty$ as $n \rightarrow \infty$ asymptotic consistency is a more delicate property because the critical value $c^*_{n,p}(1-\alpha; \widehat{\Omega}_n)$  depends on the dimensions $d, t$ and diverges as $n \rightarrow \infty$. Thus, intuitively, in high dimensions the bootstrap test is consistent only if the signal of the alternative diverges faster than the critical value as $n \rightarrow \infty$. This is the content of the next theorem.

\begin{theorem}[Consistency under high-dimensional alternatives]\label{theorem:ConsistencyLocalAlternatives} 
	Suppose that Assumptions~\ref{assumption:ControlThirdMoments}--\ref{assumption:ConsistentCovariance} hold. Then, for all $1 \leq p \leq \infty$ the following is true.
	\begin{itemize}
		\item[(i)] For any $\alpha \in (0,1)$,
		\begin{align*}
			\lim_{n \rightarrow \infty}\sup_{(\mu_n)_{n \geq 1} \in \mathcal{A}_p} \mathrm{E}_{\mu_n}[\varphi_\alpha(T_{n,p}, \widehat{\Omega}_n)] = \alpha,
		\end{align*}
		where $	\mathcal{A}_p = \left\{\left(\mu_n\right)_{n \geq 1} :  \sqrt{n}\|R\mu_n -r\|_p \lesssim  \mathrm{Var}(\|Z\|_p)^{1/2} \right\}$.
		\item[(ii)] For any $\alpha \in (0,1/2)$,
		\begin{align*}
			\lim_{n \rightarrow \infty}\sup_{(\mu_n)_{n \geq 1} \in \mathcal{B}_p} \mathrm{E}_{\mu_n}[\varphi_\alpha(T_{n,p}, \widehat{\Omega}_n)]  < 1,
		\end{align*}
		where $\mathcal{B}_p = \left\{\left(\mu_n\right)_{n \geq 1} : \sqrt{n}\|R\mu_n -r\|_p \lesssim \mathrm{E}\|Z\|_p \right\}$.
		\item[(iii)] For any $\alpha \in (0,1)$,
		\begin{align*}
			\lim_{n \rightarrow \infty}\inf_{(\mu_n)_{n \geq 1} \in \mathcal{C}_p}\mathrm{E}_{\mu_n}[\varphi_\alpha(T_{n,p}, \widehat{\Omega}_n)] = 1,
		\end{align*}
		where $\mathcal{C}_p = \left\{ \left(\mu_n\right)_{n \geq 1} :  \mathrm{E}\|Z\|_p \lesssim \sqrt{n}\|R\mu_n -r\|_p \right\}$.
	\end{itemize}
\end{theorem}
In above theorem, $\mathcal{A}_p$ is the set of alternatives (i.e. sequences of mean vectors $(\mu_n)_{n\geq 1}$) whose signals $\sqrt{n}\|R\mu_n - r\|_p$ are asymptotically negligible compared to the noise level of the Gaussian proxy statistic $\|Z\|_p$ with $Z \sim N(0, \Omega)$. As one would expect, alternatives in $\mathcal{A}_p$ are undetectable in the sense that the (asymptotic) probability of rejection is not larger than the (asymptotic) size of the test $\varphi_\alpha(T_{n,p}, \widehat{\Omega}_n)$. The set $\mathcal{B}_p$ contains those alternatives whose signals are asymptotically smaller than the expected value of $\|Z\|_p$. The bootstrap test is inconsistent for such alternatives, but unlike in the case of $\mathcal{A}_p$ the power function is not necessarily tacked to size $\alpha$. Lastly, the set $\mathcal{C}_p$ contains alternatives whose signals asymptotically dominate noise level and expected value of the Gaussian proxy statistic $\|Z\|_p$. The bootstrap test is asymptotically consistent for these alternatives. 

From Theorem~\ref{theorem:ConsistencyLocalAlternatives} it is conceivable that there exist alternatives $(\mu_n)_{n\geq 1} \in \mathcal{H}_1$ and exponents $p, q \geq 1$ such that the bootstrap test based on $T_{n,p}$ is consistent whereas the test based on $T_{n,q}$ is inconsistent. In the following we show that this intuition is indeed correct.

\begin{definition}[Sparse or dense alternatives]\label{definition:SparseDenseAlternatives}
	Given a sparsity level $s \in \{1, \ldots, t\}$ and a signal strength $\delta > 0$, define
	\begin{align*}
		\mathcal{D}_{\delta,s} = \left\{ \left(\mu_n\right)_{n \geq 1} : \|R\mu_n -r\|_0 = s, \: |(R\mu_n -r)_k| \asymp \delta \sqrt{\mathrm{E}[Z_k^2]}, \:  k \in \{j : (R\mu_n -r)_j \neq 0\} \right\},
	\end{align*}
	where $Z \sim N(0, \Omega)$. The set $\mathcal{D}_{\delta,s}$ contains ``sparse'' alternatives if $s \ll t$ and ``dense'' alternatives if $s \asymp t$. 
\end{definition}
\begin{remark}[On the uniform signal-to-noise ratio in Definition~\ref{definition:SparseDenseAlternatives}]
	In the definition of sparse and dense alternatives we require that asymptotically the signal-to-noise ratios of all non-zero entries in the alternative are of the same order $\delta > 0$. We believe that this is a realistic assumption, since in applications the high-dimensional data is either collected from a single data generating mechanism or, if aggregated from different sources, it is standardized to guarantee that the measurements have comparable levels of variability.
\end{remark}

\begin{proposition}\label{lemma:example:SparseDenseAlternatives} Consider the setup in Definition~\ref{definition:SparseDenseAlternatives}. Recall the notation of Theorem~\ref{theorem:ConsistencyLocalAlternatives}. Suppose that Assumptions~\ref{assumption:ControlThirdMoments}--\ref{assumption:ConsistentCovariance} hold and let $1 \leq p < \log t \leq q \leq \infty$. Denote by $\omega_{(1)}^2 \leq \ldots \leq \omega_{(t)}^2$ the  ordered values of $\omega_k^2 = \mathrm{E}[|Z_k|^2]$ where $Z \sim N(0, \Omega)$ and define $a_t^2 = \omega_{(t)}/\omega_{(1)}$ and $b_t^p = \big(t^{-1}\sum_{k=1}^t \omega_{(k)}^p\big)/ \big(s^{-1}\sum_{k=1}^s\omega_{(k)}^p\big)$ for $s \geq 1$.
	\begin{itemize}
		\item[(i)] If $s \ll t$ and $\sqrt{\frac{a_t^2\log t}{n}} \lesssim \delta \lesssim \sqrt{\frac{b_t^2(t/s)^{2/p}}{n}}$, then
		\begin{align*}
			\mathcal{D}_{\delta,s} \subseteq \mathcal{B}_p \cap \mathcal{C}_q.
		\end{align*}
		Thus, there exist signal strengths $\delta > 0$ such that $\varphi_\alpha(T_{n,q}, \widehat{\Omega}_n)$ with $q \geq \log t$ is consistent and $\varphi_\alpha(T_{n,p}, \widehat{\Omega}_n)$ with $p < \log t$ is inconsistent against sparse $(\mu_n)_{n\geq 1} \in \mathcal{D}_{\delta,s}(\Omega)$.
		
		\item[(ii)] If $s \asymp t $ and 
		$\sqrt{\frac{b_t^2 p}{n}} \lesssim \delta \lesssim \sqrt{\frac{a_t^2 \log t}{n} }$, then
		\begin{align*}
			\mathcal{D}_{\delta,s} \subseteq \mathcal{C}_p \cap \mathcal{B}_q.
		\end{align*}
		Thus, there exist signal strengths $\delta > 0$ such that $\varphi_\alpha(T_{n,p}, \widehat{\Omega}_n)$ with $p < \log t$ is consistent and $\varphi_\alpha(T_{n,q}, \widehat{\Omega}_n)$ with $q \geq \log t$ is inconsistent against dense $(\mu_n)_{n\geq 1} \in \mathcal{D}_{\delta,s}$.
	\end{itemize}
\end{proposition}

Above proposition substantiates the statistical folklore that tests based on sum-of-squares type statistics (i.e. $T_{n,p}$ with $p < \log t$) have good power against dense alternatives whose signals $R\mu_n - r$ are spread out over a large number of coordinates, whereas tests based on maximum type statistics (i.e. $T_{n,p}$ with $p \geq \log t$) are more powerful against ``sparse" alternatives with only a few strong signals $R\mu_n - r$.

Moreover, part (ii) of the proposition shows that the bootstrap tests $\varphi_\alpha(T_{n,p}, \widehat{\Omega}_n)$ with $p < \log t$ are consistent against dense alternatives with coordinate-wise effects as small as $O(1/\sqrt{n})$ (provided that $b_t = O(1)$). Thus, these tests can detect deviations from the global null hypothesis when the coordinate-wise effects are well below the detection threshold once traditional multiple comparison adjustments such as Bonferoni-based family-wise error rate control are taken into account. This is useful for establishing the existence of significant aggregate effects composed of otherwise negligible individual effects.

\subsection{Worst case Bahadur slope}\label{subsec:BahadurSlope}
In the preceding section we showed that there exist alternatives and test statistics $T_{n,p}$ and $T_{n,q}$ such the bootstrap test based on $T_{n,p}$ is consistent and the one based on $T_{n,q}$ is inconsistent. In such a situation, we will always choose the consistent over the inconsistent test. But how can we choose between two statistics $T_{n,p}$ and $T_{n,q}$ if the associated tests $\varphi_\alpha(T_{n,p}, \widehat{\Omega}_n)$ and  $\varphi_\alpha(T_{n,q}, \widehat{\Omega}_n)$ are both consistent against the alternative under consideration?

In the classical setting with fixed dimensions $d,t < \infty$ and in certain high-dimensional settings when the test statistic has a (Gaussian) limit distribution~\cite[e.g.][]{chen2010two, cai2014two, wang2015high-dimensional, he2021Asymptotic} this question can be easily answered by computing the asymptotic relative efficiency of the tests based on Pitman alternatives. Since our test statistics do not have a limiting distribution and the Gaussian approximation results hold only under the null hypothesis, we cannot pursue this approach. Instead, we suggest comparing the tests $\varphi_\alpha(T_{n,p}, \widehat{\Omega}_n)$ based on their \emph{observed significance levels} as this requires only knowledge of the sampling distribution under the null hypothesis. Observed significance levels can be interpreted as a measure of strength of the observed sample as evidence against the null hypothesis. In particular, the smaller an observed significance level, the greater the evidence against the null hypothesis. The disadvantage of using observed significance levels is that they depend on the (unknown) distribution of the data. To obtain concrete results, we will therefore impose additional distributional assumptions.

To formalize this approach, let $F_{n, \mu}$ be the sampling distribution of $T_{n,p}$ if the $X_i$'s have mean $\mu$. Since the bootstrap test $\varphi_\alpha(T_{n,p}, \widehat{\Omega}_n)$ is significant for large values of $T_{n,p}$, the observed significance level is
\begin{align}\label{eq:subsec:BahadurSlope-1}
	\sup_{ \mu \in \mathcal{H}_0}\mathrm{P}_{\mu} \left( T_{n,p} \geq u \right) \mid_{u = T_{n,p}} = \sup_{ \mu \in \mathcal{H}_0}\big( 1 - F_{n, \mu}(T_{n,p}) \big).
\end{align}
If the $X_i$'s were drawn from $P_{\mu}$ with $\mu \in \mathcal{H}_0 = \{\mu_0\}$ the observed significance level would be uniformly distributed on $(0,1)$, whereas if the $X_i$'s were drawn from $P_{\mu}$ with $\mu \notin \mathcal{H}_0$ (and dimensions $d, t < \infty$ fixed) it would converge to zero at an exponential rate by Cram{\'e}r's large deviation principle. It is therefore of interest to understand the behavior of 
\begin{align}\label{eq:subsec:BahadurSlope-2}
	-\frac{1}{n^\gamma} \log\Big( \sup_{ \mu \in \mathcal{H}_0}\big(1 - F_{n, \mu}(T_{n,p}) \big)\Big) \quad{} \quad{} \mathrm{as}\quad{} \quad{} n \rightarrow \infty \quad{}\quad{} \mathrm{for\: some\:\:}\gamma \in (0,1].
\end{align}
The limit, if it exists, is called the Bahadur slope~\citep[][Chapter 10.4]{serfling1980approximation}. The larger the Bahadur slope of a test, the ``faster'' the rate at which the test gathers evidence against the null hypothesis at any given sample size $n$. In high dimensions with $d, t \rightarrow \infty$ as $n \rightarrow \infty$ obtaining the exact limit is extremely challenging as it requires a large deviation principle for high-dimensional random vectors. We believe that developing such a large deviation principle should be addressed in a separate, future paper. Here, we only derive an asymptotic lower bound on the quantity in eq.~\eqref{eq:subsec:BahadurSlope-2} which affords ranking tests based on their ``worst case'' Bahadur slopes. 

\begin{theorem}[Asymptotic bounds on Bahadur slope]\label{theorem:BahadurARE}
	Suppose that Assumptions~\ref{assumption:ControlThirdMoments}--\ref{assumption:ConsistentCovariance} hold and $\limsup_{n \rightarrow \infty} \mathrm{E}[\|Z\|_p]/(\sqrt{n}\|R\mu_n - r\|_p) < 1$ for $1 \leq p \leq \infty$.
	\begin{itemize}
		\item[(i)] If the $X_i$'s have a log-concave density, i.e. the density satisfies $f = e^{-\varphi}$ with $\varphi$ convex, then
		\begin{align*}
			-\frac{1}{\sqrt{n}} \log\Big( \sup_{ \mu \in \mathcal{H}_0}\big(1 - F_{n, \mu}(T_{n,p}) \big)\Big) \gtrsim \frac{\|R\mu_n -r\|_p}{\|\Omega^{1/2}\|_{2 \rightarrow p}} + o_p(1).
		\end{align*}
		\item[(ii)] If the $X_i$'s have a strictly log-concave density, i.e. the density satisfies $f = e^{-\varphi}$ where $\varphi'' \geq \lambda I$ with $\lambda > 0$, then
		\begin{align*}
			-\frac{1}{n} \log\Big( \sup_{ \mu \in \mathcal{H}_0}\big(1 - F_{n, \mu}(T_{n,p}) \big)\Big) \gtrsim \left(\frac{\|R\mu_n -r\|_p}{\|\Omega^{1/2}\|_{2 \rightarrow p}/\sqrt{\lambda}}\right)^2 + o_p(1).
		\end{align*}
		\item[(iii)] If the $X_i$'s are Gaussian random vectors and, in addition, $\log \left(\|R\mu_n -r\|_p/\|\Omega^{1/2}\|_{2 \rightarrow p}\right) = o(n)$, then
		\begin{align*}
				-\frac{1}{n} \log\Big( \sup_{ \mu \in \mathcal{H}_0}\big(1 - F_{n, \mu}(T_{n,p}) \big)\Big) \asymp  \left(\frac{\|R\mu_n -r\|_p}{\|\Omega^{1/2}\|_{2 \rightarrow p}}\right)^2  +  o_p(1).
		\end{align*}
	\end{itemize}
\end{theorem}

\begin{remark}[On the distributional assumptions]
	As in Theorem~\ref{theorem:UnbiasedTest} we impose log-concavity to navigate around certain shortcomings of the available technical tools. If there existed a large deviation principle for random vectors whose dimension $d$ varies with the sample size $n$, we would not need this assumption. However, in the absence of such a result, log-concavity conveniently guarantees via Borell's lemma that $\ell_p$-norms of the averages of random vectors have sub-exponential tails for all $d$ and $n$. In the case of Gaussian random vectors we are able to further strengthen the result and derive a matching upper bound. Notice that the condition $\log \left(\|R\mu_n -r\|_p/\|\Omega^{1/2}\|_{2 \rightarrow p}\right) = o(n)$ imposes an implicit constraint on the growth rate of the dimension $d$ and the Poincar{\'e} constant $\|\Omega^{1/2}\|_{2 \rightarrow p}$. This hints at some of the challenges of obtaining large deviation principles for (arbitrary) random vectors whose dimension is a function of the sample size.
\end{remark}

\begin{remark}[On the Poincar{\'e} constant $\|\Omega^{1/2}\|_{2 \rightarrow p}$]
	Since $\sqrt{\mathrm{Var}(\|Z\|_p)}  \leq \|\Omega^{1/2}\|_{2 \rightarrow p}$, one might wonder whether the lower bounds in Theorem~\ref{theorem:BahadurARE} should in fact depend on the variance of the Gaussian proxy statistic $\|Z\|_p$, just as the results in Section~\ref{subsec:AsympCorrectness}. Statement (iii) of Theorem~\ref{theorem:BahadurARE} shows that this is not the case and that the Poincar{\'e} constant $\|\Omega^{1/2}\|_{2 \rightarrow \infty}$ is indeed the correct quantity. While the results in Section~\ref{subsec:AsympCorrectness} rely on Gaussian approximation arguments and thus ``local'' deviations characterized by the standard deviation $\sqrt{\mathrm{Var}(\|Z\|_p)}$, Theorem~\ref{theorem:BahadurARE} depends on large deviations and concentration of measure type arguments involving the Poincar{\'e} constant.
\end{remark}

Theorem~\ref{theorem:BahadurARE} complements Theorem~\ref{theorem:ConsistencyLocalAlternatives} insofar as it introduces another qualitative measure by which we can compare tests: While Theorem~\ref{theorem:ConsistencyLocalAlternatives} enables us to compare tests based on whether they are consistent or not, Theorem~\ref{theorem:BahadurARE} allows to rank consistent tests based on their Bahadur slopes. We formalize this idea in the following definition:

\begin{definition}[Bahadur preferred test]\label{definition:BahadurPreferredTest} Consider two tests $\varphi_\alpha(T_{n, p}, \widehat{\Omega}_n)$ and $\varphi_\alpha(T_{n,p'}, \widehat{\Omega}_n)$, both consistent against alternatives in $\mathcal{H} \subseteq \mathcal{H}_1$. We say that $\varphi_\alpha(T_{n,p}, \widehat{\Omega}_n)$ is Bahadur preferred over $\varphi_\alpha(T_{n,p'}, \widehat{\Omega}_n)$ when testing against alternatives in $\mathcal{H}$ if, for all alternatives in $\mathcal{H}$, the lower bound on the Bahadur slope of $\varphi_\alpha(T_{n,p}, \widehat{\Omega}_n)$ is at least as large as the lower bound on the Bahadur slope of $\varphi_\alpha(T_{n,p'}, \widehat{\Omega}_n)$.
\end{definition}

To illustrate this approach of ranking tests, we continue the discussion of sparse and dense alternatives started in Proposition~\ref{lemma:example:SparseDenseAlternatives}. 
	
\begin{proposition}\label{lemma:example:BahadurSlope} Consider the setup in Definition~\ref{definition:SparseDenseAlternatives}. Recall the notation of Theorem~\ref{theorem:ConsistencyLocalAlternatives}. Suppose that Assumptions~\ref{assumption:ControlThirdMoments}--\ref{assumption:ConsistentCovariance} and any of the sufficient conditions of Theorem~\ref{theorem:BahadurARE} hold. In addition, suppose that $\Omega$ is a diagonal matrix.
	\begin{itemize}
		\item[(i)] If $p \in [\log t, \infty]$, $s \ll t$, and $\sqrt{(\log t)/n} \lesssim \delta$, then
		\begin{align*}
			\mathcal{D}_{\delta, s} \subseteq \mathcal{C}_{\log t} \cap \mathcal{C}_p \quad{}\quad{} \mathrm{and} \quad{}\quad{} \frac{\|R\mu_n-r\|_{\log t}}{\|\Omega^{1/2}\|_{2 \rightarrow \log t}} \geq  \frac{\|R\mu_n-r\|_p}{\|\Omega^{1/2}\|_{2 \rightarrow p}} \quad{} \forall \: (\mu_n)_{n \geq 1} \in \mathcal{D}_{\delta,s}.
		\end{align*}
		Thus, when testing against sparse alternatives $\mathcal{D}_{\delta, s}$ there exist signal strengths $\delta > 0$ such that $\varphi_\alpha(T_{n,\log t},\widehat{\Omega}_n)$ is Bahadur preferred over all tests  $\varphi_\alpha(T_{n,p},\widehat{\Omega}_n)$ with $p \in [\log t, \infty]$.
		
		\item[(ii)] If $p, q \in [1, \log t]$, $p \leq q$, $s \asymp t$, and $\sqrt{q/n}\lesssim \delta$, then
		\begin{align*}
			\mathcal{D}_{\delta, s} \subseteq \mathcal{C}_2 \cap \mathcal{C}_p \quad{}\quad{} \mathrm{and} \quad{}\quad{} \frac{\|R\mu_n-r\|_2}{\|\Omega^{1/2}\|_{2 \rightarrow 2}} \geq  \frac{\|R\mu_n-r\|_p}{\|\Omega^{1/2}\|_{2 \rightarrow p}} \quad{} \forall \: (\mu_n)_{n \geq 1} \in \mathcal{D}_{\delta,s}.
		\end{align*}
		Thus, when testing against dense alternatives $\mathcal{D}_{\delta, s}$ there exist signal strengths $\delta > 0$ such that $\varphi_\alpha(T_{n,2},\widehat{\Omega}_n)$ is Bahadur preferred over all tests $\varphi_\alpha(T_{n,p},\widehat{\Omega}_n)$ with $p \in [1, q]$.
	\end{itemize}
\end{proposition}
\begin{remark}[On the assumption on $\Omega$]
	Definition~\ref{definition:BahadurPreferredTest} is only useful if we can compute the Poincar{\'e} constant $\|\Omega^{1/2}\|_{2 \rightarrow p}$. Evaluating $\|\Omega^{1/2}\|_{2 \rightarrow p}$ for arbitrary positive semi-definite $\Omega$ is difficult and can even be NP-hard~\citep[][]{foucart2013mathematical}. Therefore, the assumption that $\Omega$ is diagonal is a compromise: on the one hand, it is specific enough to simplify the computation of $\|\Omega^{1/2}\|_{2 \rightarrow p}$; on the other hand, it is broad enough to cover a large class of tests.
\end{remark}
\begin{remark}
	Obviously, case (i) also holds for $s \asymp t$. But then case (ii) simply implies that $\varphi_\alpha(T_{n,2},\widehat{\Omega}_n)$ is Bahadur preferred over all tests $\varphi_\alpha(T_{n,p},\widehat{\Omega}_n)$, $p \in [1, \infty]$.
\end{remark}

Proposition~\ref{lemma:example:BahadurSlope} is a surprising result: Unlike one could have expected, neither $\varphi_\alpha(T_{n,\infty},\widehat{\Omega}_n)$ nor $\varphi_\alpha(T_{n,1},\widehat{\Omega}_n)$ are Bahadur preferred tests. Instead, under the setup of Definition~\ref{definition:SparseDenseAlternatives} the bootstrap tests based on $T_{n, \log t}$ and $T_{n,2}$ are ``best'' against sparse and dense alternatives, respectively. The simulation study in Section~\ref{sec:NumericalStudies} further corroborates this conclusion even in finite samples.

\begin{remark}[On the relation of observed significance level and $p$-value] 
	Observed significance level and $p$-value are two different concepts that yield the same numbers only in a few special situations. Given an observed test statistic $T_{n,p}$ the $p$-value of $\varphi_\alpha(T_{n,p},\widehat{\Omega}_n)$ is defined as 
	\begin{align*}
		\inf \{ \alpha \in (0,1) : \varphi_\alpha(T_{n,p},\widehat{\Omega}_n) = 1\} = \mathrm{P} \left( T^*_{n,p} \geq u \mid X_1, \ldots, X_n\right) \mid_{u = T_{n,p}} =  1 - \widehat{G}_n(T_{n,p}),
	\end{align*}
	where $\widehat{G}_n$ is the conditional distribution of the Gaussian proxy statistic
	\begin{align*}
		T_{n,p}^* := \|Z_n\|_p, \quad{}\quad{} Z_n \mid X_1, \ldots, X_n \sim N(0, \widehat{\Omega}_n).
	\end{align*}
	Typically, $\widehat{G}_n \neq F_{n,\mu}$, and hence $p$-value and observed significance level of $\varphi_\alpha(T_{n,p},\widehat{\Omega}_n)$ are different.
	
	From Theorem~\ref{theorem:Bootstrap-Lp-Norm-Quantiles} we know that $\widehat{G}_n$ and $F_{n, \mu}$ are close in Kolomgorov-Smirnov distance for all $\mu \in \mathcal{H}_0$. This implies that $p$-values and observed significance levels are close as well. However, this does not imply that the $p$-value and observed significance level accumulate evidence against the null hypothesis at the same rate. Indeed, from Theorem~\ref{theorem:BahadurARE} (iii) we infer that
	\begin{align*}
		-\frac{1}{n}\log(1 - \widehat{G}_n(T_{n,p})) \gtrsim  \left(\frac{\|R\mu_n -r\|_p}{\|\widehat{\Omega}_n^{1/2}\|_{2 \rightarrow p}}\right)^2  +  o_p(1),
	\end{align*}
	independent of the sampling distribution of $T_{n,p}$. Incidentally, under the distributional assumptions of Theorem~\ref{theorem:BahadurARE}, ranking tests $\varphi_\alpha(T_{n,p},\widehat{\Omega}_n)$ by the rates at which their $p$-values vanish yields the same order as ranking them by their worst case Bahadur slopes.
\end{remark}

\section{Three modifications of the basic bootstrap hypothesis test}\label{sec:Modifications}

\subsection{A test with power against sparse and dense alternatives}\label{subsec:PowerEnhancement}
The first modification of the basic bootstrap hypothesis test that we consider is motivated by the theoretical analysis in Section~\ref{subsec:BahadurSlope} and aims at enhancing the power of the test.

Consider the implications of Proposition~\ref{lemma:example:BahadurSlope} when the signal strength satisfies $\delta \asymp \sqrt{(\log t)/n}$: In the case of dense alternatives with $s \asymp t$ we easily verify (by Theorem~\ref{theorem:ConsistencyLocalAlternatives}) that all bootstrap tests based on $\ell_p$-statistics are consistent. Hence, combining Proposition~\ref{lemma:example:BahadurSlope} (i) and (ii) (with $q = \log t$) we deduce that the test based on $T_{n,2}$ is ``best''  among all tests based on $\ell_p$-statistics with any $p \in [1, \infty]$. Moreover, in the case of sparse alternatives with $s \ll t$, Theorem~\ref{theorem:ConsistencyLocalAlternatives} implies that tests based on $T_{n,p}$ with $p \in [1, \log t)$ are inconsistent. Hence, by Proposition~\ref{lemma:example:BahadurSlope} (i) the bootstrap test based on $T_{n,\log t}$ is ``best'' not just among all tests based on $\ell_p$-statistics with $p \in [\log t, \infty]$ but with any $p \in [1, \infty]$.

We therefore propose to combine test statistics $T_{n,2}$ and $T_{n,\log t}$ into the single test statistic
\begin{align*} 
	W_n := T_{n,2} + T_{n,\log t} = \|RS_n - \sqrt{n} r\|_2  + \|RS_n - \sqrt{n} r\|_{\log t}, \quad{}\quad{} S_n = \frac{1}{\sqrt{n}} \sum_{i=1}^n X_i.
\end{align*}
Given a nominal level $\alpha \in (0,1)$, we reject the null hypothesis if and only if
\begin{align*} 
	W_n \geq c^*_W(1- \alpha;\widehat{\Omega}_n),
\end{align*}
where $c^*_W(\alpha; \widehat{\Omega}_n) = \inf\left\{s \in \mathbb{R} : \mathrm{P}(W_n^* \leq s \mid X_1, \ldots, X_n)  \geq \alpha \right\}$
is the conditional $\alpha$-quantile of the Gaussian proxy statistic
\begin{align*} 
	W_n^* := \|Z_n\|_2 + \|Z_n\|_{\log t}, \quad{}\quad{} Z_n \mid X_1, \ldots, X_n \sim N(0, \widehat{\Omega}_n),
\end{align*}
and $\widehat{\Omega}_n$ is a positive semi-definite estimate of $\Omega = R \Sigma R'$. We define the \emph{modified bootstrap test based on $W_n$ at level $\alpha$} as
\begin{align*}
	\varphi_\alpha(	W_n,\widehat{\Omega}_n) := \mathbf{1}\left\{W_n \geq c^*_W(1- \alpha;\widehat{\Omega}_n)\right\}.
\end{align*}
Given the theoretical results in Appendix~\ref{sec:TheoryBootstrap} it is straightforward to show that this test has asymptotically the correct size $\alpha$:

\begin{proposition}[Asymptotic size $\alpha$ test]\label{theorem:SizeAlphaTestL1LInfty}
	Suppose that Assumptions~\ref{assumption:ControlThirdMoments}--\ref{assumption:ConsistentCovariance} hold for $p \in \{2, \log t\}$. Then,
	\begin{align*}
		\lim_{n \rightarrow \infty}\sup_{\alpha \in (0,1)} \sup_{\mu \in \mathcal{H}_0}\Big|\mathrm{E}_\mu[\varphi_\alpha(W_n, \widehat{\Omega}_n) ] - \alpha \Big| =0.
	\end{align*}
\end{proposition}

Moreover, by simple modification of Theorem~\ref{theorem:ConsistencyLocalAlternatives} we conclude that the modified bootstrap test is consistent against both, dense and sparse alternatives. 

\begin{proposition}\label{theorem:ConsistencySizeAlphaTestL1LInfty} 
Suppose that Assumptions~\ref{assumption:ControlThirdMoments}--\ref{assumption:ConsistentCovariance} hold for $p \in \{2, \log t\}$. Then, for $\alpha \in (0,1)$ arbitrary,
\begin{align*}
	\lim_{n \rightarrow \infty}\inf_{(\mu_n)_{n \geq 1} \in \mathcal{W}}\mathrm{E}_{\mu_n}[\varphi_\alpha(W_n, \widehat{\Omega}_n)] = 1,
\end{align*}
where $\mathcal{W} = \left\{ \left(\mu_n\right)_{n \geq 1} :  \mathrm{E}\|Z\|_2 + \mathrm{E}\|Z\|_{\log t} \lesssim \sqrt{n}\|R\mu_n -r\|_2 + \sqrt{n}\|R\mu_n -r\|_{\log t} \right\}$.	
\end{proposition}
We easily verify that $\mathcal{W} \supset \mathcal{D}_{\delta, s}$ for sparse alternatives characterized by sparsity level $s \lesssim t$ and signal strength $\delta \gtrsim \sqrt{(\log t)/n}$ as well as for dense alternatives characterized by $s \asymp t$ and $\delta \gtrsim 1/\sqrt{n}$ (recall notation of Definition~\ref{definition:SparseDenseAlternatives}). Thus, the modified bootstrap test $\varphi_\alpha(	W_n,\widehat{\Omega}_n)$ is consistent against dense and sparse alternatives simultaneously.

\begin{remark}[Comparison with the \emph{power enhancement method} by~\cite{fan2015power}]
	The starting point for~\possessivecite{fan2015power} power enhancement method is the observation that sum-of-square type statistics based on (such as $T_{n,2}^2$) tend to have low power against sparse alternatives. To rectify this, they suggest to add a \emph{power enhancement component} $J_0$ to the initial sum-of-square type statistics. Such $J_0$ should satisfy three intuitive properties: First, non-negativity, i.e. $J_0 \geq 0$ almost surely; second, no size distortion, i.e. $\inf_{\mu \in \mathcal{H}_0}\mathrm{P}_\mu(J_0 = 0) \rightarrow 1$; third, power enhancement, i.e. $J_0 \rightarrow \infty$ in probability under some alternative $\mathcal{H}_1$.
	The power enhancement component typically depends on a tuning parameter $\delta > 0$ that determines the minimal (asymptotic) signal strength of an alternative $\mathcal{H}_1$ for which $J_0 \rightarrow \infty$ in probability. Adapted to our setup,~\cite{fan2015power} propose to set $\delta \asymp (\log \log n) \sqrt{(\log t)/n}$.
	
	The philosophy of our bootstrap test $\varphi_\alpha(W_n,\widehat{\Omega}_n)$ and the power enhancement by~\cite{fan2015power} are thus very different. Our bootstrap test is not based on the idea of enhancing the power the test statistic $T_{n,2}$ without distorting the size of the test (with the critical value being the $\alpha$-quantile of the sampling distribution of $T_{n,2}$). Instead, our bootstrap test is based on $T_{n,2} + T_{n,\log t}$ and the critical value is the $\alpha$-quantile of the sampling distribution of $T_{n,2} + T_{n,\log t}$. Hence, size distortion is no issue. Moreover, the bootstrap test does not depend on a tuning parameter and the minimal signal strength that the bootstrap test can detect can be as small as $\delta \gtrsim \sqrt{(\log t)/n}$ (sparse alternative) and $\delta \gtrsim 1/\sqrt{n}$ (dense alternative).
\end{remark}

\subsection{A test statistic for elliptically distributed data}\label{subsec:EllipticallyDistributedData}
We have so far kept the assumptions on the data generating process minimal and have mostly relied on the three moment conditions Assumptions~\ref{assumption:ControlThirdMoments}--\ref{assumption:ConsistentCovariance}. These moment conditions implicitly constrain the growth rates of the dimensions $d,t$ relative to the sample size $n$. In this section, we propose a test statistic based on self-normalized random vectors which has asymptotic correct size under stronger assumptions on the data generating process but considerably weaker implicit constraints on $d,t, n$. Consequently, the self-normalized test performs well in even higher dimensions.

The results in this section hold under only the following assumption, i.e. in the following  Assumptions~\ref{assumption:ControlThirdMoments}--\ref{assumption:ConsistentCovariance} are not needed.

\begin{assumption}[Elliptically distributed data]\label{assumption:EllipticallyDistributed}
	Let $X, X_1, \ldots, X_n \in \mathbb{R}^d$ be a simple random sample of elliptically distributed random vectors, i.e.
\begin{align*}
	X = \mu_n + \Gamma \theta U^{(s)},
\end{align*}
where $U^{(s)} \in \mathbb{R}^s$ is a random vector uniformly distributed on the $\ell_2$-norm unit sphere $\mathbb{S}^{s-1}$, $\theta$ is a non-negative random variable independent of $U$, $\mu_n \in \mathbb{R}^d$, and $\Gamma \in \mathbb{R}^{d \times s}$. In particular, $X$ has mean $\mu_n$ and covariance matrix $\Sigma = s^{-1} \mathrm{E}[\theta^2]\Gamma \Gamma' \in \mathbb{R}^{d \times d}$ with $\mathrm{rank}(\Sigma) = s \wedge d$.
\end{assumption}
\begin{remark}[On elliptically distributed data]
	The class of elliptical distributions includes many classical non-Gaussian multivariate distributions such as multivariate t-distribution, multivariate logistic distribution, Kotz-type multivariate distribution, and Pearson II type multivariate distribution~\citep{fang1990symmetric}. Elliptical distributions have proved useful for modeling tail dependence (i.e. clustering of extremes) in financial data~\citep{mcneil2015quantitative}. Assumption~\ref{assumption:EllipticallyDistributed} does not (!) imply Assumptions~\ref{assumption:ControlThirdMoments}--\ref{assumption:ConsistentCovariance}.
\end{remark}

Let $X_1, \ldots, X_n \in \mathbb{R}^d$ be a simple random sample of elliptically distributed random vectors and denote by $\widetilde{X}, \widetilde{X}_1, \ldots, \widetilde{X}_n \in \mathbb{R}^t$ their projections on the Euclidean unit sphere with center $r \in \mathbb{R}^t$,
\begin{align*}
	\widetilde{X} := \begin{cases}
	\frac{RX - r}{\|RX - r\|_2},  &\|RX - r\|_2 \neq 0,\\
	0, & o/w.
	\end{cases}
\end{align*}
From the characterization in Assumption~\ref{assumption:EllipticallyDistributed} it follows that $\widetilde{X} \sim R\Gamma U^{(s)}$, where $U^{(s)} \in \mathbb{R}^s$ is a random vector uniformly distributed on the $\ell_2$-norm unit sphere $\mathbb{S}^{s-1}$ for all $\mu \in \mathcal{H}_0$. In particular, $\mathrm{E}_\mu[\widetilde{X}] = 0$ for all $\mu \in \mathcal{H}_0$. Moreover, in the proofs to below results we will show that the projection $\widetilde{X}$ satisfies the simple sufficient conditions of Lemma~\ref{lemma:example:HeavyTailedData}. We therefore use these projections $\widetilde{X}_1, \ldots, \widetilde{X}_n$ as building blocks for a \emph{self-normalized test statistic}: Define
\begin{align*} 
	V_n := \left\|\frac{1}{\sqrt{n}} \sum_{i=1}^n \widetilde{X}_i \right\|_2, 
\end{align*}
and, given a nominal level $\alpha \in (0,1)$, we reject the null hypothesis if and only if
\begin{align*} 
	V_n \geq c^*_V(1- \alpha;\widetilde{\Omega}_n),
\end{align*}
where $c^*_V(\alpha; \widetilde{\Omega}_n) = \inf\left\{s \in \mathbb{R} : \mathrm{P}(V_n^* \leq s \mid X_1, \ldots, X_n)  \geq \alpha \right\}$
is the conditional $\alpha$-quantile of the Gaussian proxy statistic
\begin{align*} 
	V_n^* := \|Z_n\|_2, \quad{}\quad{} Z_n \mid X_1, \ldots, X_n \sim N(0, \widetilde{\Omega}_n),
\end{align*}
and $\widetilde{\Omega}_n$ is a positive semi-definite estimate of $\Omega = \mathrm{E}[\widetilde{X}\widetilde{X}']$. We define the \emph{self-normalized bootstrap test based on $V_n$ at level $\alpha$} as
\begin{align*}
	\varphi_\alpha(	V_n,\widetilde{\Omega}_n) := \mathbf{1}\left\{V_n \geq c^*_V(1- \alpha;\widetilde{\Omega}_n)\right\}.
\end{align*}

We have the following result:
\begin{proposition}[Asymptotic size $\alpha$ test]\label{theorem:SizeAlphaTestEllipticallyDistributed}
	Let $\widetilde{\Omega}_n = n^{-1} \sum_{i=1}^n\widetilde{X}_i \widetilde{X}_i'$. If Assumption~\ref{assumption:EllipticallyDistributed} holds and $\|R\Gamma\Gamma'R'\|_{op} = o(n^{1/6} s^{-1}\|R\Gamma\Gamma' R'\|_F )$, then
	\begin{align*}
		\lim_{n \rightarrow \infty}\sup_{\alpha \in (0,1)} \sup_{\mu \in \mathcal{H}_0}\Big|\mathrm{E}_\mu[\varphi_\alpha(V_n, \widetilde{\Omega}_n) ] - \alpha \Big| =0.
	\end{align*}
\end{proposition}
Since $s \geq \mathrm{rank}(R\Gamma\Gamma'R')$, we typically have $s^{-1}\|R\Gamma\Gamma'R'\|_F = O(1)$. Therefore, the requirement $\|R\Gamma\Gamma'R'\|_{op} = o(n^{1/6} s^{-1} \|R\Gamma\Gamma' R'\|_F )$ is very mild and easily satisfied even when the covariance matrix $\Sigma$ has full rank and/ or is dense. In particular, the dimensions $d,t, s$ can be arbitrarily large relative to the sample size $n$. However, these mild requirements come at a slight disadvantage: The estimate $\widetilde{\Omega}_n = n^{-1} \sum_{i=1}^n\widetilde{X}_i \widetilde{X}_i'$ depends on the hypothesis $R\mu_n = r$. Thus, the Gaussian proxy statistic $V^*_n$ is no longer ancillary. Consequently, we can establish consistency of the self-normalized bootstrap test only for one alternative hypothesis at a time, not uniformly over all alternatives. We make this explicit in the next proposition by indexing the estimate $\widetilde{\Omega}_n$ by the alternative $\mu_n$.

\begin{proposition}[Consistency under high-dimensional alternatives]\label{theorem:ConsistencySizeAlphaTestEllipticallyDistributed} 
	For arbitrary $\mu \in \mathbb{R}^d$ define $\widetilde{\Omega}_{\mu} := n^{-1} \sum_{i=1}^n\widetilde{X}_{i, \mu} \widetilde{X}_{i, \mu}'$, where $\widetilde{X}_{i, \mu} = R(X_i - \mu)/\|R(X_i - \mu)\|_2$. If Assumption~\ref{assumption:EllipticallyDistributed} holds and $\|R\Gamma\Gamma'R'\|_{op} = o(n^{1/6} s^{-1}\|R\Gamma\Gamma' R'\|_F )$, then for $\alpha \in (0,1)$ arbitrary,
		\begin{align*}
			\lim_{n \rightarrow \infty}\inf_{(\mu_n)_{n \geq 1} \in \mathcal{E}}\mathrm{E}_{\mu_n}[\varphi_\alpha(V_n, \widetilde{\Omega}_{\mu_n})] = 1,
		\end{align*}
		where $\mathcal{E} = \left\{ \left(\mu_n\right)_{n \geq 1} :  \lim_{n \rightarrow \infty} \sqrt{n} \left\| \mathrm{E}_{\mu_n}\left[\frac{RX - r}{\|RX - r\|_2}\right] \right\|_2 = \infty \right\}$.
\end{proposition}
Notice that $\mathrm{E}[\|Z\|_2] \leq 1$ for $Z \sim N(0, \widetilde{\Omega}_{\mu_n})$ (see proof of Proposition~\ref{theorem:SizeAlphaTestEllipticallyDistributed}). Hence, we can re-write the set of alternatives as $\mathcal{E} = \{ \left(\mu_n\right)_{n \geq 1} :  \mathrm{E}[\|Z\|_2] = o\big(\sqrt{n}\| \mathrm{E}_{\mu_n}[(RX - r)/\|RX - r\|_2]\|_2 \big) \}$. In this form, $\mathcal{E}$ matches the intuition developed in Section~\ref{subsec:AsympSizeLevel} and Theorem~\ref{theorem:ConsistencyLocalAlternatives}: An alternative $\mu_n \in \mathcal{E}$ is detectable because its signal $\sqrt{n} \|\mathrm{E}_{\mu_n}[(RX - r)/\|RX - r\|_2]\|_2 $ dominates the expected value of the Gaussian proxy statistic $\|Z\|_2$.

\begin{remark}[Comparison with the \emph{high-dimensional nonparametric test} by~\cite{wang2015high-dimensional}]
	Our self-normalized bootstrap test is inspired by~\cite{wang2015high-dimensional} who, to the best of our knowledge, were the first to propose a sum-of-square type statistic for elliptically distributed data based on the normalized random vectors $X_1/ \|X_1\|_2, \ldots, X_n/\|X_n\|_2$.
	
	The most interesting and consequential difference between their test and ours lies in the theoretical analysis. The asymptotic theory in~\cite{wang2015high-dimensional} is driven by ``asymptotics over the dimension''. By this we mean that the test statistic is constructed in such a way that increasing the dimension reduces the variability of each individual summand of the test statistic. Then, provided that the entries in the high-dimensional random vectors are not be too correlated, a Martingale CLT applies. Earlier attempts in this direction by~\cite{bai1996effect} and~\cite{chen2010two} required data generated by a factor model and strong moment conditions. The assumptions in~\cite{wang2015high-dimensional} are much weaker, but they are violated if $\mathrm{tr}(\Sigma) = O(1)$, i.e. if the data lie in a (approximate) low-dimensional subspace (a.k.a. if the data are highly correlated). In contrast, the asymptotic theory of the bootstrap test relies on Gaussian approximation which is related to the classical ``asymptotics over the sample''. As seen in Propositions~\ref{theorem:SizeAlphaTestEllipticallyDistributed} and~\ref{theorem:ConsistencySizeAlphaTestEllipticallyDistributed} this results in very mild assumptions on the dimension and the covariance matrix $\Sigma$.
\end{remark}

\subsection{Bootstrap test statistics based on approximate sample averages}\label{subsec:ApproximateSampleAverages}
The basic bootstrap test (including above two modifications) is based on the average of i.i.d. random vectors. However, in many applications, the ``natural'' statistic to base an $\ell_p$-statistic on is a nonlinear transformation of the data, i.e. $T_{n,p} = \|F(X_1, \ldots, X_n)\|_p$ with $F: \bigtimes_{i=1}^n \mathbb{R}^d \rightarrow \mathbb{R}^t$ and $F$ nonlinear. For example, $F$ could be a function that captures certain pre-processing steps such as removal of biases or handling of missingness and heteroscedasticity in a data set. Or $F$ could be a procedure that returns an estimate of a high-dimensional regression vector, e.g. solution to a generalized estimating equation.

In classical statistics with fixed dimension, asymptotic theory can be developed whenever $F$ is sufficiently smooth via continuous mapping theorem, delta method, or linearization. An analogous result holds true for the high-dimensional bootstrap test:

\begin{proposition}\label{theorem:ApproximateSampleAverages}
	Consider a random sample $Y_1, \ldots, Y_n$. Let $F: \bigtimes_{i=1}^n \mathbb{R}^d \rightarrow \mathbb{R}^t$. If, under the null hypothesis,
	\begin{align}\label{eq:theorem:ApproximateSampleAverages-1}
		F(Y_1, \ldots, Y_n) -\sqrt{n} r = \frac{1}{\sqrt{n}}\sum_{i=1}^n R(X_i - \mu) + W_n \quad{}\quad{} \mathrm{and} \quad{}\quad{} \|W_n\|_p = o_p\left(\sqrt{\mathrm{Var}(\|Z\|_p)}\right),
	\end{align}
	where  $X_1, \ldots X_n \in \mathbb{R}^d$ are i.i.d. random vectors (possibly changing with $n$) with mean $\mu \in \mathbb{R}^d$ and positive semi-definite covariance matrix $\Sigma \neq \mathbf{0} \in \mathbb{R}^{d \times d}$, $Z \sim N(0, \Omega)$ with $\Omega \equiv R \Sigma R' \neq \mathbf{0} \in \mathbb{R}^{t \times t}$ and $R \in \mathbb{R}^{t \times d}$ (possibly changing with $n$), then Theorems~\ref{theorem:SizeAlphaTest} and~\ref{theorem:ConsistencyLocalAlternatives} and Propositions~\ref{theorem:ConfidenceSets},~\ref{theorem:SizeAlphaTestL1LInfty}, and~\ref{theorem:ConsistencySizeAlphaTestL1LInfty} hold with the test statistic $\|F(Y_1, \ldots, Y_n)- \sqrt{n}r\|_p$ substituted for $T_{n,p}= \|RS_n - \sqrt{n} r\|_p$.
\end{proposition}
Intuitively, this result says that the bootstrap hypothesis test is valid whenever $F$ affords a linear expansion with a remainder term $W_n$ that is asymptotically negligible compared to the noise level the leading term of the linear expansion (which is captured by $\sqrt{\mathrm{Var}(\|Z\|_p)}$). Theorems~\ref{theorem:UnbiasedTest} and~\ref{theorem:BahadurARE} do not necessarily hold under the condition in eq.~\eqref{eq:theorem:ApproximateSampleAverages-1} because they require stronger non-asymptotic lower bounds on tail probabilities of log-concave measures.

\section{A practical guide to the bootstrap hypothesis test}\label{sec:Implementation}

\subsection{Three sources of bias}\label{subsec:HeuristicsBias}
To test a global null hypothesis such as~\eqref{eq:sec:Intro-1} via the bootstrap test $\varphi_\alpha(T_{n,p},\widehat{\Omega}_n)$ we need to compute the conditional $(1-\alpha)$-quantile $c^*_p(1- \alpha;\widehat{\Omega}_n)$ of the Gaussian proxy statistic $T_{n,p}^*$. Since there does not exist an analytic expression of the distribution of $T_{n,p}^*$, we suggest approximating $c^*_p(1- \alpha;\widehat{\Omega}_n)$ via Monte Carlo simulations. Thus, in practice, we test hypothesis~\eqref{eq:sec:Intro-1} using the \emph{Monte Carlo bootstrap hypothesis test} defined as
\begin{align*} 
	\varphi_{\alpha,B}(T_{n,p},\widehat{\Omega}_n) := \mathbf{1}\left\{T_{n,p} \geq c^*_{p,B}(1- \alpha;\widehat{\Omega}_n)\right\},
\end{align*}
where $c^*_{p, B}(1-\alpha; \widehat{\Omega}_n)$ denotes the $\lfloor (1-\alpha) B\rfloor$th order statistic of the random sample $T_{n,p,1}^*, \ldots, T_{n,p,B}^*$ with $T_{n,p,b}^* := \|Z_{n,b}\|_p$ and $Z_{n,b} \mid X_1, \ldots, X_n \sim_{iid} N(0, \widehat{\Omega})$ for all $b = 1, \ldots, B$. 

Since the test statistic is non-pivotal (i.e. the sampling distribution of the test statistic depends on the unknown population covariance $\Omega$) theoretical and empirical results from the classical, low-dimensional bootstrap literature suggest that the Monte Carlo bootstrap test might be severely biased~\citep[e.g.][]{davison1986bootstrap, beran1987prepivoting, hall1988resampling, shi1992accurate, shao1995jackknife}. In the following we therefore discuss three sources of potential biases in our high-dimensional setting. This discussion lays the foundation for developing strategies to mitigate two of these biases (Sections~\ref{subsec:CovarianceEstimation} and~\ref{subsec:SphericalBootstrap}). 
 
Consider the following decomposition of the bias of the Gaussian Monte Carlo bootstrap test:
\begin{align}\label{eq:subsec:HeuristicBias-2}
	\begin{split}
	\mathrm{E}_\mu\left[ \varphi_{\alpha,B}(T_{n,p},\widehat{\Omega}_n) \right] - \alpha &= \underbrace{\mathrm{E}_\mu\left[ \varphi_{\alpha}(T_{n,p},\Omega) \right] - \alpha}_{\mathrm{Gaussian\:approximation\:error}} + \:\:\:\underbrace{\mathrm{E}_\mu\left[ \varphi_{\alpha}(T_{n,p},\widehat{\Omega}_n) \right]  - \mathrm{E}_\mu\left[ \varphi_{\alpha}(T_{n,p},\Omega) \right]}_{\mathrm{Gaussian\:comparison\:error}} \\
	&\quad{}  + \underbrace{\mathrm{E}_\mu\left[ \varphi_{\alpha,B}(T_{n,p},\widehat{\Omega}_n) \right] -\mathrm{E}_\mu\left[ \varphi_{\alpha}(T_{n,p},\widehat{\Omega}_n) \right]}_{\mathrm{Monte\:Carlo\:error} }.
	\end{split}
\end{align}

The ``Gaussian approximation error'' in above display is the error of approximating the $(1-\alpha)$th quantile of the test statistic $T_{n,p}$ with the $(1-\alpha)$th quantile of the Gaussian proxy statistic $\|Z\|_p$, where $Z \sim N(0, \Omega)$. In other words, this error measures the distance between the distributions of $T_{n,p}$ and $\|Z\|_p$. Since in the classical setting with fixed dimensions bootstrapping procedures are known to be consistent if and only if the bootstrapped statistic satisfies a CLT~\citep[][Theorem 3.6.1]{mammen1993bootstrap, vandervaart1996weak}, it is intuitive that the Gaussian approximation error is part of the bias decomposition~\eqref{eq:subsec:HeuristicBias-2}. Theorem~\ref{theorem:CLT-Lp-Norm} (in Appendix~\ref{subsec:GeneralResults}) provides an upper bound on this error uniformly in $\alpha \in (0,1)$. Under Assumptions~\ref{assumption:ControlThirdMoments} and~\ref{assumption:RatioMoments} this upper bound vanishes as $n, d, t \rightarrow \infty$. The rate at which the Gaussian approximation error vanishes is only of theoretical interest as it is independent of the Monte Carlo procedure and the data $X_1, \ldots, X_n$. For all practical purposes, the only way to reduce this error is by increasing the sample size $n$.

The ``Gaussian comparison error'' in bias decomposition~\eqref{eq:subsec:HeuristicBias-2} arises from comparing the $(1-\alpha)$th quantiles of the Gaussian proxy statistics $\|Z\|_p$ and $\|Z_n\|_p$, where $Z \sim N(0,\Omega)$ and $Z_n \mid X_1, \ldots, X_n \sim N(0, \widehat{\Omega}_n)$, respectively. From Theorem~\ref{theorem:Bootstrap-Lp-Norm} (in Appendix~\ref{subsec:GeneralResults}) we deduce that the rate at which this error vanishes depends on how fast $\|\widehat{\Omega}_n- \Omega\|_{q \rightarrow p}/\mathrm{Var}(\|Z\|_p) \rightarrow 0$ in probability, where $1/p+1/q = 1$. 
Consequently, the choice of the estimator $\widehat{\Omega}_n$ has significant impact on the Gaussian comparison error. We discuss in Section~\ref{subsec:CovarianceEstimation} how to exploit this simple observation to improve the accuracy of the bootstrap hypothesis test.  

The third term in bias decomposition~\eqref{eq:subsec:HeuristicBias-2} is the ``Monte Carlo error'' of approximating the $(1-\alpha)$th quantile of the Gaussian proxy statistic $T_{n,p}^*$. Under the mild assumption that the test statistic $T_{n,p}$ has a density (or probability mass function) we have the following result:
\begin{theorem}[Gaussian Monte Carlo Error]\label{theorem:BoundMonteCarloError}
	Let $M_{n,p} = \esssup_{z, \mu} f_{n,p, \mu}(z)$, where $f_{n,p, \mu}$ is the density of the test statistics $T_{n,p}$ when the $X_i$'s have mean $\mu$. Then, for all $1 \leq p \leq \infty$,
		\begin{align*}
			\sup_{\alpha \in (0,1)} \sup_{\mu} \left| \mathrm{E}_\mu\left[ \varphi_{\alpha,B}(T_{n,p},\widehat{\Omega}_n) \right] -\mathrm{E}_\mu\left[ \varphi_{\alpha}(T_{n,p},\widehat{\Omega}_n) \right] \right| &\lesssim B^{-1/2} (\log B) \sqrt{  M_{n,p}^2 \mathrm{E}\left[ \mathrm{Var}(T_{n,p}^* \mid X_1, \ldots, X_n)\right]},
		\end{align*}
	where $\lesssim$ hides an absolute constant independent of $p, n, d, t, B$, and the distribution of the $X_i$'s.
\end{theorem}
\begin{remark}[On size, consistency, and power of the Gaussian Monte Carlo bootstrap test]\label{remark:SizeGaussianMonteCarloTest}
	The upper bound on the Monte Carlo error increases in the conditional variance of the Gaussian proxy statistic $T_{n,p}^*$ and the mode $M_{n,p}$ of the density of the statistic $T_{n,p}$, but decreases in the number of Monte Carlo samples $B$. Thus, Theorems~\ref{theorem:SizeAlphaTest} and~\ref{theorem:ConsistencyLocalAlternatives} and Propositions~\ref{theorem:ConfidenceSets},~\ref{theorem:SizeAlphaTestL1LInfty}--\ref{theorem:ApproximateSampleAverages} continue to hold if we substitute the Monte Carlo Bootstrap test $\varphi_{\alpha,B}(T_{n,p},\widehat{\Omega}_n)$ for the infeasible bootstrap test $\varphi_\alpha(T_{n,p},\widehat{\Omega}_n)$ and choose $B \gtrsim n \vee \left( M_{n,p}^2 \mathrm{E}\left[ \mathrm{Var}(T_{n,p}^* \mid X_1, \ldots, X_n)\right]\right)^{1 + \gamma}$, $\gamma > 0$ arbitrary.
\end{remark}


The dependence of the upper bound in Theorem~\ref{theorem:BoundMonteCarloError} on the number of Monte Carlo samples $B$ and the conditional variance $\mathrm{Var}(T_{n,p}^* \mid X_1, \ldots, X_n)$ is intuitively obvious. It is, however, noteworthy, because its proof 
applies to any Monte Carlo sampling scheme and does not rely on the specific Gaussian proxy statistic $T_{n,p}^*\equiv \|Z_n\|_p$, $Z_n \mid X_1, \ldots, X_n \sim N(0, \widehat{\Omega}_n)$. This raises the prospect that it might be possible to reduce the Monte Carlo error by designing an alternative (potentially non-Gaussian) Monte Carlo sampling scheme based on proxy statistics that have smaller variances than the Gaussian proxy statistics $T_{n,p}^*$.

Finding such an alternative Monte Carlo sampling scheme is a challenging task since for certain exponents $p$ and covariance matrices $\widehat{\Omega}_n$ the Gaussian proxy statistic $T_{n,p}^*$ is ``superconcentrated'', i.e. for certain exponents and covariance matrices the variance of $T_{n,p}^*$ is a decreasing function of the dimension of the Gaussian random vector $Z_n$ and hence (asymptotically) already very small (see Remark~\ref{remark:SuperconcentrationGaussianRV}). In Section~\ref{subsec:SphericalBootstrap} we succeed in developing a simple and practical alternative using hypercontractivity tools for Markov semigroups~\citep[i.e.][]{cordero-erausquin2012hypercontractive, tanguy2017quelques}. Numerical experiments in Section~\ref{sec:NumericalStudies} show that this alternative Monte Carlo sampling scheme substantially improves the accuracy of the bootstrap tests in finite samples.

\begin{remark}[On the superconcentration of $\ell_p$-norms of Gaussian random vectors in $\mathbb{R}^d$]\label{remark:SuperconcentrationGaussianRV}
	In the case of isotropic Gaussian random vectors, $Z \sim N(0, I_d)$, the variance of $\|Z\|_p$ has several phase transitions depending on the interplay between dimension $d$ and exponent $p$~\citep[][]{lytova2017variance}. In particular, $\mathrm{Var}(\|Z\|_2) \asymp 1$ (i.e. not superconcentrated) and $\mathrm{Var}(\|Z\|_\infty) \asymp (\log d)^{-1}$ (i.e. superconcentrated), which, by Theorem~\ref{theorem:BoundMonteCarloError}, lends theoretical support to the empirical observation that the bias of the Gaussian Monte Carlo hypothesis test based on the $\ell_\infty$-norm is significantly smaller than the one based on the $\ell_2$-norm.
	
	In the statistically more relevant case of an-isotropic Gaussian random vectors, $Z \sim N(0, \Sigma)$, the variance of $\|Z\|_p$ is less well studied. In principle, one can obtain upper bounds on $\mathrm{Var}(\|Z\|_p)$ via Talagrand's $L_1$-$L_2$ inequality~\citep{cordero-erausquin2012hypercontractive}. For $\ell_2$- and $\ell_\infty$-norms the following useful estimates have been derived via ad-hoc approaches: $\mathrm{Var}(\|Z\|_2) \asymp \mathrm{tr}(\Sigma^2)/\mathrm{tr}(\Sigma)$~\citep[][p. 344]{valettas2019tightness} and $\mathrm{Var}(\|Z\|_\infty) \lesssim \sigma^2 (\log \log 2d/ \log 2d)^{1/4} + \sigma^2 (\log \sum_{j,k} (2d)^{-2/(1 + r_{jk})}/ \log 2d)^{1/4}$, where $\Sigma = (\sigma_{jk})_{j,k=1}^d$, $\sigma_{jj} = \sigma^2$ and $r_{jk} = \sigma_{jk}/\sigma^2$ for all $1 \leq j,k \leq d$~\citep[][Theorem 1.11]{chatterjee2008chaos}. While less obvious than in the isotropic case, the takeaway of the an-isotropic case is similar: $\|Z\|_\infty$ is superconcentrated (provided that the correlations between the entries in the vector $Z$ are ``small''), whereas $\|Z\|_2$ is not superconcentrated (irrespective of the correlations). Thus, by way of Theorem~\ref{theorem:BoundMonteCarloError} this provides again a partial explanation for why the Gaussian Monte Carlo hypothesis test based on the $\ell_\infty$-norm is more accurate than the one based on the $\ell_2$-norm.
\end{remark}


\subsection{Bias reduction via structured covariance matrix estimation}\label{subsec:CovarianceEstimation}
In the preceding section we argued that the rate at which $\|\widehat{\Omega}_n- \Omega\|_{q \rightarrow p}/\mathrm{Var}(\|Z\|_p) \rightarrow 0$ impacts the accuracy of the (Monte Carlo) bootstrap test via the Gaussian approximation error. We now show how to reduce this error by leveraging additional structural information about the population covariance $\Omega$. Note that reducing the Gaussian approximation error is tantamount to relaxing the sufficient conditions of Lemmas~\ref{lemma:example:SubGaussianData}--\ref{lemma:example:HeavyTailedData}. The results in this section are far from comprehensive, they are meant to be illustrative.

Throughout, let $\Omega =  R\Sigma R' = (\omega_{jk})_{j,k=1}^t \in \mathbb{R}^{t \times t}$ and $\widehat{\Omega}_n^{\mathrm{naive}}= n^{-1} \sum_{i=1}^n R(X_i - \bar{X}_n)(X_i - \bar{X}_n)'R'$ be the population and (naive) sample covariance matrix of $RX_1, \ldots, RX_n \in \mathbb{R}^t$, respectively. Consider the following canonical structural assumption on $\Omega$:

\begin{assumption}[Approximately sparse covariance matrix]\label{assumption:ApproxSparseCovariance}
	There exist constants $\gamma \in [0,1]$ and $R_\gamma > 0$ such that $	\max_{1 \leq j\leq t} \sum_{k=1}^t |\omega_{jk}|^\gamma \leq R_\gamma$.
\end{assumption}
Under Assumption~\ref{assumption:ApproxSparseCovariance} it is natural to estimate the covariance matrix via (hard) thresholding of the sample covariance as $\mathcal{T}_{\lambda}(\widehat{\Omega}_n^{\mathrm{naive}})$, where
\begin{align*}
	\mathcal{T}_\lambda(\Omega) := \big(\omega_{jk}\mathbf{1}\{|\omega_{jk}| > \lambda\}\big)_{j,k=1}^t,
\end{align*}
and $\lambda > 0$ a suitably chosen thresholding parameter. Since the bootstrap test requires a positive semi-definite estimate of the covariance matrix, we project this estimate onto the cone of positive semi-definite matrices. The resulting positive semi-definite projection $\mathcal{T}_{\lambda}^+(\widehat{\Omega}_n^{\mathrm{naive}})$ maintains the same order of error as the original thresholding estimate~\citep[][p. 275]{avella-medina2018robust}.

We have the following result; it should be compared to the case $p=2$ in Lemmas~\ref{lemma:example:SubGaussianData}--\ref{lemma:example:HeavyTailedData}:

\begin{lemma}\label{lemma:example:ApproxSparseCovariance}
	Let $X_1, \ldots X_n \in \mathbb{R}^d$ be a simple random sample of random vectors with mean $\mu_n \in \mathbb{R}^d$ and positive semi-definite covariance matrix $\Sigma \in \mathrm{R}^{d \times d}$. Suppose that $\Omega = R \Sigma R'$ satisfies Assumption~\ref{assumption:ApproxSparseCovariance}. Let $\widehat{\Omega}_n := \mathcal{T}_{\lambda}^+(\widehat{\Omega}_n^{\mathrm{naive}})$ for $\lambda > 0$ to be specified below. Then, Assumptions~\ref{assumption:ControlThirdMoments}--\ref{assumption:ConsistentCovariance} hold for exponent $p=2$
	\begin{itemize}
		\item[(i)] if the $X_i$'s are sub-Gaussian, $\lambda \asymp \sqrt{ (\log t)/n} \vee (\log t)/n$, and
		\begin{align*}
				\frac{r(\Omega)}{\sqrt{r(\Omega^2)}} = o(n^{1/6}) \quad{} \quad{} \mathrm{and}\quad{}\quad{} \frac{r(\Omega)}{r(\Omega^2)} \frac{R_\gamma}{\omega_{(t)}^{2\gamma} } \lambda^{1- \gamma} = o(1);
		\end{align*}
		\item[(ii)] if the $X_i$'s have a log-concave density, $\lambda \asymp \big((\log t)/n\big)^{1/4} \vee \sqrt{(\log t)/n}$, and
		\begin{align*}
			\frac{r(\Omega)}{\sqrt{r(\Omega^2)}} = o(n^{1/6}) \quad{} \quad{} \mathrm{and}\quad{}\quad{} \frac{r(\Omega)}{r(\Omega^2)} \frac{R_\gamma}{\omega_{(t)}^{\gamma-1}} \lambda^{1- \gamma} = o(1);
		\end{align*}
		\item[(i)] if the $X_i$'s are heavy-tailed, $\lambda \asymp \big((\log t)/n\big)^{1/4} \vee \sqrt{(\log t)/n}$, and there exists $s > 3$ such that
		\begin{align*}
			m_{2, 3} \sqrt{\frac{r(\Omega)}{r(\Omega^2)}} = o(n^{1/6}), \quad{} \frac{m_{2, s}}{m_{2,3}} = o(n^{1/3-1/s}), \quad{} \frac{r(\Omega)}{r(\Omega^2)} \frac{R_\gamma}{\omega_{(t)}^{\gamma-1}} \lambda^{1- \gamma} = o(1),
		\end{align*}
	\end{itemize}
		where $r(M) := \mathrm{tr}(M)/\|M\|_{op}$ is the effective rank of $M \in \{\Omega, \Omega^2\}$ and $m_{2,s}^s := \mathbb{E}[\|X_i-\mu_n\|_2^s] / \|\Omega\|_{op}^{s/2}$.
\end{lemma}
Three comments are in order: First, the rates in above lemma are less stringent than those in Lemmas~\ref{lemma:example:SubGaussianData}--\ref{lemma:example:HeavyTailedData}, i.e. the thresholding estimator reduces the Gaussian approximation error in an asymptotic sense. Numerical evidence shows that the improvement is sizable even in finite samples (see Section~\ref{sec:NumericalStudies}). Second, we only provide results for exponent $p =2$ because structural assumptions on the population covariance $\Omega$ do not improve the rates for exponent $p= \infty$. In fact, the rates in above lemma for $p=2$ reduce approximately to the rates of $p = \infty$ in Lemmas~\ref{lemma:example:SubGaussianData}--\ref{lemma:example:HeavyTailedData} raised to the power $1- \gamma$. Third, other structural assumptions on $\Omega$ such as bandedness or low-rank result in similar improvements of the asymptotic rates. Above lemma is merely intended as a proof of concept.

\begin{remark}[Comparison with the \emph{Gaussian multiplier bootstrap} by~\cite{chernozhukov2013GaussianApproxVec}]\label{remark:GaussianMultiplierBootstrap}
	The bootstrap test is at least as accurate as the Gaussian multiplier test by~\cite{chernozhukov2013GaussianApproxVec}. Moreover, if additional information about the population covariance matrix is available, the bootstrap test can substantially outperform the Gaussian multiplier test. To see this, recall that the Gaussian multiplier bootstrap is defined as
	\begin{align*}
		\varphi_{\alpha}^m(T_{n,p}, X_1, \ldots, X_n) := \mathbf{1}\left\{T_{n,p} \geq c^m_p(1- \alpha;X_1, \ldots, X_n)\right\},
	\end{align*}
	where $c^m_p(\alpha;X_1, \ldots, X_n) = \inf\left\{u \in \mathbb{R} : \mathrm{P}(T_{n,p}^m \leq u \mid X_1, \ldots, X_n)  \geq \alpha \right\}$
	is the conditional $\alpha$-quantile of the Gaussian multiplier proxy statistic
	\begin{align*}
		T_{n,p}^m := \left\| \frac{1}{\sqrt{n}} \sum_{i=1}^n \xi_i R(X_i - \bar{X}_n)\right\|_p, \quad\quad \xi_1, \ldots, \xi_n \sim_{iid} N(0,1),
	\end{align*}
	and the Gaussian multipliers $\xi_1, \ldots, \xi_n$ are independent of the data $X_1, \ldots, X_n$. Since $n^{-1/2} \sum_{i=1}^n \xi_i\\ R(X_i - \bar{X}_n) \mid X_1, \ldots, X_n \sim N(0, \widehat{\Omega}_n^{\mathrm{naive}})$, the Gaussian multiplier test $\varphi_{\alpha}^m(T_{n,p}, X_1, \ldots, X_n)$ is equivalent to the bootstrap test $\varphi_{\alpha}(T_{n,p}, \widehat{\Omega}_n^{\mathrm{naive}})$. Thus, the Gaussian multiplier bootstrap does not offer a way to leverage structural information about the population covariance matrix.
\end{remark}

\subsection{Improved efficiency via spherical bootstrapping}\label{subsec:SphericalBootstrap}
As detailed in Section~\ref{subsec:HeuristicsBias}, we can, in theory, reduce the bias of the Monte Carlo bootstrap test by using a more efficient non-Gaussian sampling scheme based on highly superconcentrated proxy statistics. When implementing this idea, we have to trade off two conflicting goals: On the one hand, we need to find a non-Gaussian distribution such that the proxy statistics are ``more'' superconcentrated than their Gaussian counterparts. On the other hand, the non-Gaussian distribution must not be ``too far'' from a Gaussian distribution with mean zero and covariance $\widehat{\Omega}_n$ or else our proof strategy based on Gaussian proxy statistics fails (viz. Gaussian approximation and Gaussian comparison, see Section~\ref{subsec:HeuristicsBias}).

Formalizing this bias-variance trade-off as a constrained optimization problem and finding the optimal non-Gaussian distribution is a formidable task. 
Here, we only propose one possible solution to this problem; namely, the \emph{spherical bootstrap hypothesis test} defined as
\begin{align*} 
	\varphi_{\alpha}^s(T_{n,p},\widehat{\Gamma}_n) := \mathbf{1}\left\{T_{n,p} \geq c^s_p(1- \alpha;\widehat{\Gamma}_n)\right\},
\end{align*}
where $c^s_p(\alpha; \widehat{\Gamma}_n) = \inf\left\{u \in \mathbb{R} : \mathrm{P}(T_{n,p}^s \leq u \mid X_1, \ldots, X_n)  \geq \alpha \right\}$
is the conditional $\alpha$-quantile of the \emph{spherical proxy statistic}
\begin{align*} 
	T_{n,p}^s := \sqrt{s}\|\widehat{\Gamma}_nU^{(s)}\|_p, \quad{}\quad{} U^{(s)} \sim \mathrm{Unif}(\mathbb{S}^{s-1}),  \quad{}\quad{} \widehat{\Gamma}_n\widehat{\Gamma}_n' = \widehat{\Omega}_n  \in \mathbb{R}^{t \times t},
\end{align*}
where $\mathbb{S}^{s-1}$ is the $\ell_2$-norm unit sphere in $\mathbb{R}^s$, $2 \leq s \leq t$, and $\widehat{\Omega}_n$ is a positive semi-definite estimate of $\Omega$. The corresponding feasible \emph{spherical Monte Carlo bootstrap hypothesis test} is
\begin{align*} 
	\varphi_{\alpha,B}^s(T_{n,p},\widehat{\Gamma}_n) := \mathbf{1}\left\{T_{n,p} \geq c^s_{p,B}(1- \alpha;\widehat{\Gamma}_n)\right\},
\end{align*}
where $c^s_{p, B}(1-\alpha; \widehat{\Gamma}_n)$ denotes the $\lfloor (1-\alpha) B\rfloor$th order statistic of the random sample $T_{n,p,1}^s, \ldots, T_{n,p,B}^s$ with $T_{n,p,b}^s := \sqrt{s}\|\widehat{\Gamma}_nU_b^{(s)}\|_p$ and $U_b^{(s)} \sim \mathrm{Unif}(\mathbb{S}^{s-1})$ for all $b = 1, \ldots, B$.

We expect the spherical bootstrap procedure to perform well for two reasons: First, the first two moments of $\sqrt{s}\widehat{\Gamma}_nU^{(s)}$ and $Z_n \mid X_1, \ldots, X_n \sim N(0, \widehat{\Omega}_n)$ coincide. Thus, the spherical proxy statistic $T_{n,p}^s$ and the Gaussian proxy statistic $T_{n,p}^*$ are close in Kolmogorov distance (Theorem~\ref{theorem:CLT-Lp-Norm}, Appendix~\ref{subsec:GeneralResults}). Second, the entries in the vector $U^{(s)}$ are ``negatively associated'', i.e. we deduce from Theorem 3.3 in~\cite{fang1990symmetric} that all moments of entries in $U^{(s)}$ are non-positively correlated. This suggests that $\|\widehat{\Gamma}_nU^{(s)}\|_p$ is more superconcentrated than $\|Z_n\|_p \overset{d}{=} \|\widehat{\Gamma}_nZ^{(s)}\|_p$, where $Z^{(s)} \sim N(0, I_s)$ (for a more rigorous statement, see Lemma~\ref{lemma:BoundsVarianceUniformRV} in Appendix~\ref{subsec:Miscellaneous}). 

The following two theorems confirm that this heuristic reasoning is indeed correct. The main takeaway is that the spherical Monte Carlo bootstrap test has asymptotic correct size, is consistent against high-dimensional alternatives, has the same power properties as the (Gaussian) bootstrap test, and has smaller (asymptotic) Monte Carlo error than the (Gaussian) Monte Carlo bootstrap test. Thus, for all practical purposes, we recommend using the spherical bootstrapping procedure. 

\begin{theorem}\label{theorem:SphericalBootstrapTest}
	Suppose that $\|\widehat{\Gamma}_n\|_{2 \rightarrow p} = o_p\left(s^{1/4} \sqrt{\mathrm{Var}(\|Z\|_p) } \right)$. Then,
	Theorems~\ref{theorem:SizeAlphaTest} and~\ref{theorem:ConsistencyLocalAlternatives} and Propositions~\ref{theorem:ConfidenceSets},~\ref{theorem:SizeAlphaTestL1LInfty}--\ref{theorem:ApproximateSampleAverages} hold with the spherical bootstrap test $\varphi_{\alpha}^s(T_{n,p},\widehat{\Gamma}_n)$ substituted for the (Gaussian) bootstrap test $\varphi_{\alpha}(T_{n,p},\widehat{\Omega}_n)$.
\end{theorem}

\begin{theorem}[Spherical Monte Carlo Error]\label{theorem:BoundSphericalMonteCarloError}
	Let $M_{n,p} = \esssup_{z, \mu} f_{n,p,\mu}(z)$, where $f_{n,p, \mu}$ is the density of the test statistics $T_{n,p}$ when the $X_i$'s have mean $\mu$. Then, for all $1 \leq p \leq \infty$,
	\begin{align*}
		\sup_{\alpha \in (0,1)} \sup_{\mu} \left| \mathrm{E}_\mu\left[ \varphi_{\alpha,B}^s(T_{n,p},\widehat{\Gamma}_n) \right] -\mathrm{E}_\mu\left[ \varphi_{\alpha}^s(T_{n,p},\widehat{\Gamma}_n) \right] \right| &\lesssim  B^{-1/2} s^{-1} (\log B) \sqrt{M_{n,p}^2 \mathrm{E}\big[\|\widehat{\Gamma}_n\|_{2 \rightarrow p}^2\big]},
	\end{align*}
	where $\lesssim$ hides an absolute constant independent of $p, n, d, t, s, B$, and the distribution of the $X_i$'s.
\end{theorem}

\begin{remark}[On size, consistency, and power of the spherical Monte Carlo bootstrap test]\label{remark:SizeSphericalMonteCarloTest}
	Combined with Theorem~\ref{theorem:SphericalBootstrapTest} the bound on the spherical Monte Carlo error implies that Theorems~\ref{theorem:SizeAlphaTest} and~\ref{theorem:ConsistencyLocalAlternatives} and Propositions~\ref{theorem:ConfidenceSets},~\ref{theorem:SizeAlphaTestL1LInfty}--\ref{theorem:ApproximateSampleAverages} to hold if we substitute the spherical Monte Carlo Bootstrap test $\varphi_{\alpha,B}^s(T_{n,p},\widehat{\Gamma}_n)$ for the infeasible $\varphi_\alpha^s(T_{n,p},\widehat{\Gamma}_n)$ and choose $B \gtrsim n \vee \left( s^{-2} M_{n,p}^2 \mathrm{E}\big[\|\widehat{\Gamma}_n\|_{2 \rightarrow p}^2]\right)^{1 + \gamma}$, $\gamma > 0$ arbitrary.
\end{remark}

\begin{remark}[On the rank $s$ of $\widehat{\Gamma}_n \in \mathbb{R}^{t\times s}$]
	If we estimate the covariance matrix $\Omega \in \mathbb{R}^{t \times t}$ via the sample analogue $\widehat{\Omega}_n = n^{-1}\sum_{i=1}^n R(X_i - \bar{X}_n)(X_i - \bar{X}_n)'R'$, then there exists $\widehat{\Gamma}_n \in \mathbb{R}^{t \times (n-1)}$ such that $\widehat{\Gamma}_n\widehat{\Gamma}_n' = \widehat{\Omega}_n$. Thus, in this case rank $s = n-1$ and the condition in Theorem~\ref{theorem:SphericalBootstrapTest} is typically satisfied whenever Assumption~\ref{assumption:ConsistentCovariance} holds. Moreover, for $s = n-1$ the upper bound in Theorem~\ref{theorem:BoundSphericalMonteCarloError} is (asymptotically) significantly smaller than the corresponding upper bound on the (Gaussian) Monte Carlo error of Theorem~\ref{theorem:BoundMonteCarloError}. If we estimate $\Omega \in \mathbb{R}^{t \times t}$ using a more sophisticated estimator (viz. Section~\ref{subsec:CovarianceEstimation}), rank $s$ may be smaller or larger than $n-1$. The smaller rank $s$, the smaller the variance reduction and, thus, the smaller the gains of using the spherical Monte Carlo bootstrap test. Evidence from simulations supports this theoretical insight.
\end{remark}

\begin{remark}[Comparison with the \emph{non-Gaussian multiplier bootstrap} by~\cite{deng2020beyond}]\label{remark:NonGaussianMultiplierBootstrap}
	The rationale for the non-Gaussian (multiplier) bootstrap procedures in~\cite{deng2020beyond} is very different from the one for the spherical bootstrap:~\cite{deng2020beyond} do not consider the Monte Carlo error of implementing a bootstrap procedure. Instead, they are concerned with reducing the theoretical biases corresponding to Gaussian approximation and Gaussian comparison error. By matching up to five moments of the distributions of test and proxy statistic, they achieve some improvements in the rates at which these theoretical biases vanish~\citep[compared to ][]{chernozhukov2013GaussianApproxVec}. However, the most recent results in~\cite{chernozhukov2019ImprovedCLT, chernozhukov2021NearlyOptimalCLT} show that the Gaussian multiplier bootstrap can achieve the same (or even better) rates.
	
	To the best of our knowledge, we are the first to propose a method to mitigate the Monte Carlo error associated with a high-dimensional bootstrap procedure. The practical gains and the theoretical justification of the spherical bootstrap hold regardless of any (future) improvements in the rates at which Gaussian approximation and Gaussian comparison errors vanish.
\end{remark}


\section{Numerical Experiments}\label{sec:NumericalStudies}

\subsection{Monte Carlo setup}\label{subsec:MonteCarloSim}
We study the performance of the bootstrap test and its modifications in four different data generating processes (DGPs). To get our main points across as concise as possible, we consider testing problem~\ref{eq:sec:Intro-1} with $R= I_d \in \mathbb{R}^{d \times d}$ and $r = \mathbf{0} \in \mathbb{R}^d$ (i.e. we test $H_0 : \: \mu_0= \mathbf{0}$ versus $H_1: \: \mu_0 \neq \mathbf{0}$) and only report results for test statistics with exponents $p \in \{2, \infty\}$. 
In DGPs 1--3 we generate vectors $X_i = (X_{i1}, \ldots, X_{id})'$ via the Gaussian copula model
\begin{align*}
	X_{ij} = F^{-1}\left(\Phi(Y_{ij}) \right) - 1, \quad \quad  Y_i = (Y_{i1}, \ldots, Y_{id})' \sim N(0, \Sigma), \quad \quad  1 \leq j \leq d, \quad \quad 1 \leq i \leq n,
\end{align*}
where the random vectors $Y_1, \ldots, Y_n \in \mathbb{R}^d$ are independent copies of each other and $F$ denotes the cdf of the Gamma distribution with shape and rate parameters $\alpha = \beta = 1$. The parameters $\alpha, \beta$ are chosen such that the Gamma distribution has mean and variance equal to one, skew equal to two, and an excess kurtosis of six. This ensures that the $X_i$'s are distributed markedly different from symmetric Gaussian or spherical random vectors and hence makes the bootstrap procedures a non-trivial exercise. We consider three different covariance structures for $\Sigma = (\Sigma_{jk})_{j,k=1}^d$: equicorrelated entries with $\Sigma_{jk} = 0.8 + 0.2 \times \mathbf{1}\{j=k\}$ (non-mixing dependencies) in DGP 1, weakly dependent entries with $\Sigma_{jk} = 0.8^{|j-k|}$ (Toeplitz structure of an autoregressive process of order one) in DGP 2, and $m$-dependent entries with  $\Sigma_{jk} = \max\left\{ 1- |j-k| / \lceil d^{1/3}/2 \rceil, 0\right\}$ (banded structure of an moving average process with 2, 3, and 4 lags for dimension 100, 500, 1000, respectively) in DGP 3. In DGP 4 we generate vectors $X_i \in \mathbb{R}^d$ from a multivariate $t(4)$-distribution. We choose 4 degrees of freedom because this guarantees that the first $3 + \delta$ moments exist, which is just enough for our theoretical results to hold. For the covariance matrix $\Sigma$ we choose the same Toeplitz structure as in DGP 2.

Under DGP 1, 2, and 4 we estimate the covariance matrix via the naive sample covariance. Under DGP 3 we use a thresholding estimate of the covariance matrix. Simply thresholding the sample covariance matrix and projecting it onto the cone of positive semi-definite matrices as suggested in Section~\ref{subsec:CovarianceEstimation} results in large biases. We therefore use the Cholesky-based regularization approach by~\cite{rothman2010new} which is guaranteed to produce positive semi-definite estimates. To reduce the computational burden, we assume that the correct bandwidth is known.

For DGPs 1--3 we report results for dimensions $d \in \{100, 500, 1000\}$ and sample sizes $n \in \{20, 50, 100\}$; for DGP 4 we present results for dimensions $d \in \{500, 1000, 5000\}$ and $n \in \{20, 50, 100\}$. We implement all bootstrap procedures with $B= 5000$ bootstrap samples. The QQ-plots in Section~\ref{subsec:Results} compare `nominal' and `actual' sizes of different bootstrap tests. We choose 99 `nominal' sizes $\alpha \in \{0.01, 0.02, \ldots, 0.99\}$ and compute the corresponding `actual' sizes as the relative frequency of false positives in 1000 independent Monte Carlo samples.

\subsection{Results}\label{subsec:Results}
\begin{figure}[p]
	\centering
	\includegraphics[scale=0.08]{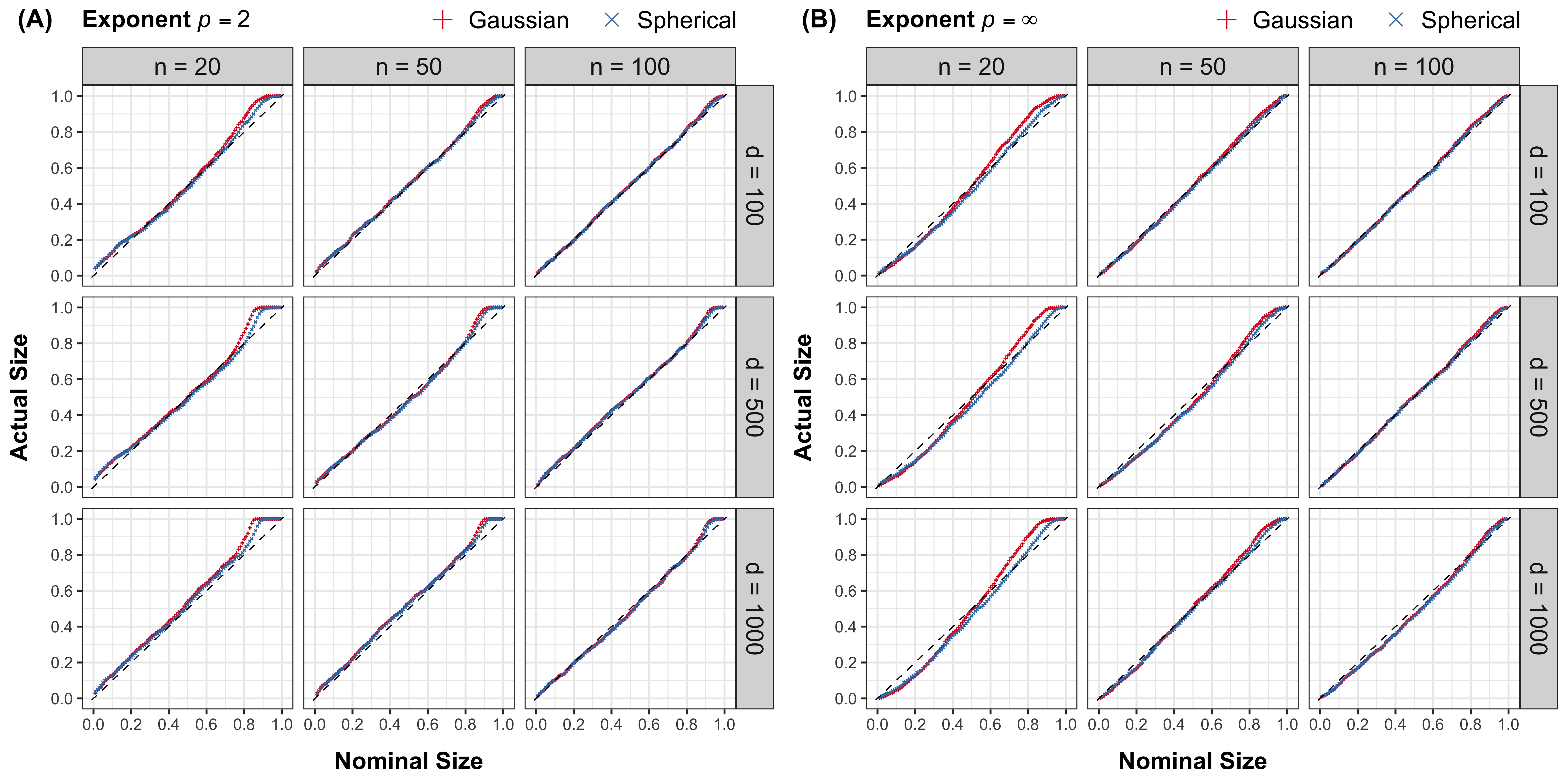}
	\caption{QQ-Plots of `nominal' and `actual' sizes of the Gaussian and spherical bootstrap test under equicorrelated design (DGP 1). Panel (A) shows simulation results for the bootstrap test with exponent $p=2$, panel (B) shows simulation results for the bootstrap test with exponent $p=\infty$. }
	\label{fig:1}
\end{figure}

\begin{figure}[p]
	\includegraphics[scale=0.08]{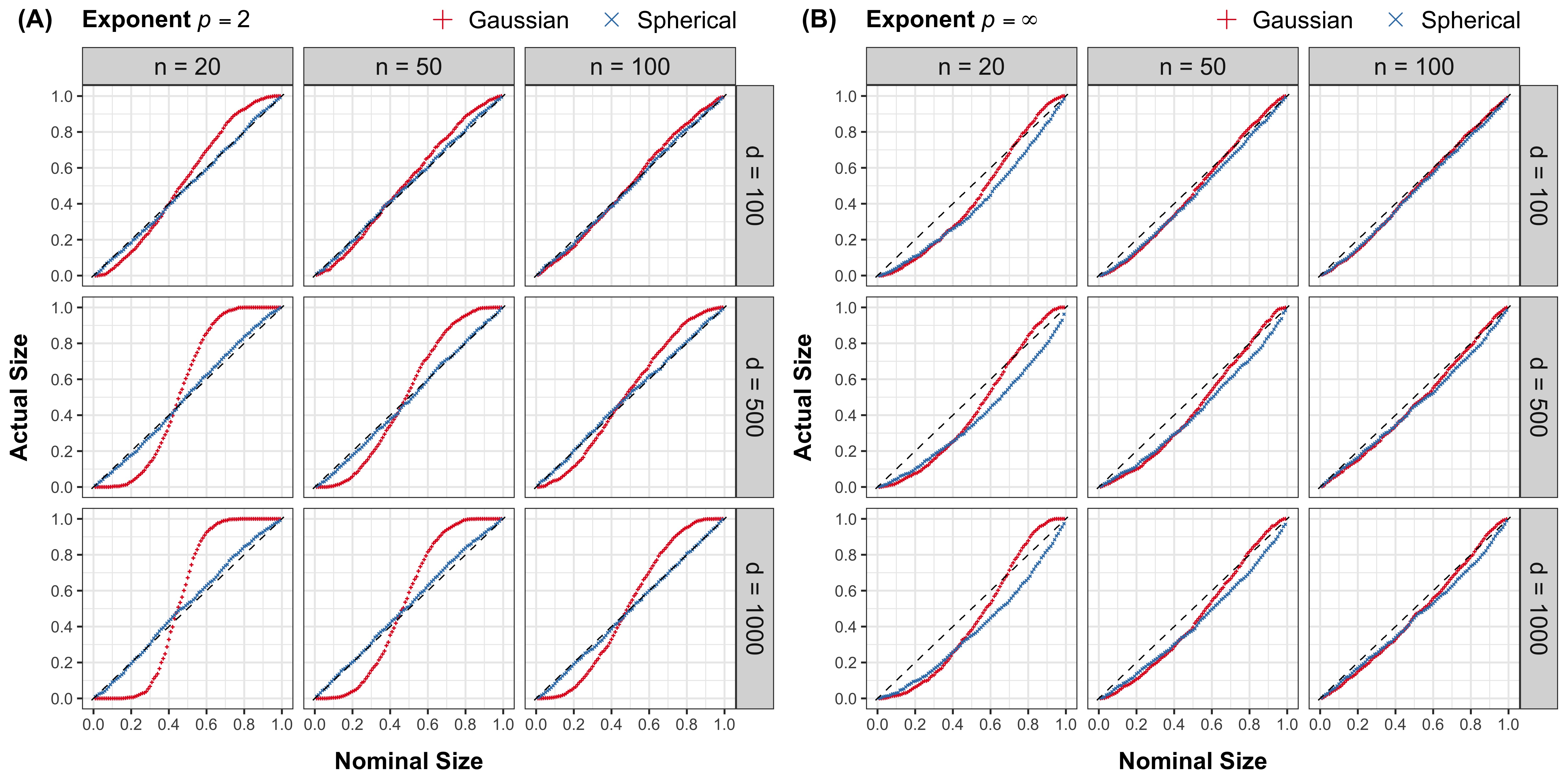}
	\caption{QQ-Plots of `nominal' and `actual' sizes of the Gaussian and spherical bootstrap test under Toeplitz design (DGP 2). Panel (A) shows simulation results for the bootstrap test with exponent $p=2$, panel (B) shows simulation results for the bootstrap test with exponent $p=\infty$. }
	\label{fig:2}
\end{figure}

In Figures~\ref{fig:1} and~\ref{fig:2} we compare the accuracy of the basic Gaussian bootstrap test from Section~\ref{subsec:BootstrapHypothesisTest} and the bias-corrected spherical bootstrap test from Section~\ref{subsec:SphericalBootstrap} for exponents $p \in \{2, \infty\}$. In the case of equicorrelated design (DGP 1) both tests are rather accurate, especially for the (important) small quantiles $\alpha \in \{0.01, \ldots, 0.1\}$. The test statistic based on the $\ell_2$-norm appears to be slightly more accurate that the test statistic based on the $\ell_\infty$-norm. For Toeplitz matrices (DGP 2) the spherical bootstrap substantially outperforms the Gaussian bootstrap across all quantile levels $\alpha \in \{0.01, 0.02, \ldots, 0.99\}$ for test statistic based on the $\ell_2$-norm. For test statistics based on the $\ell_\infty$-norm the spherical and Gaussian bootstrap test perform similarly for small quantiles $\alpha$, but for large quantiles $\alpha > 0.5$ the spherical test is significantly more conservative than the Gaussian test.

\begin{figure}[p]
	\centering
	\includegraphics[scale=0.08]{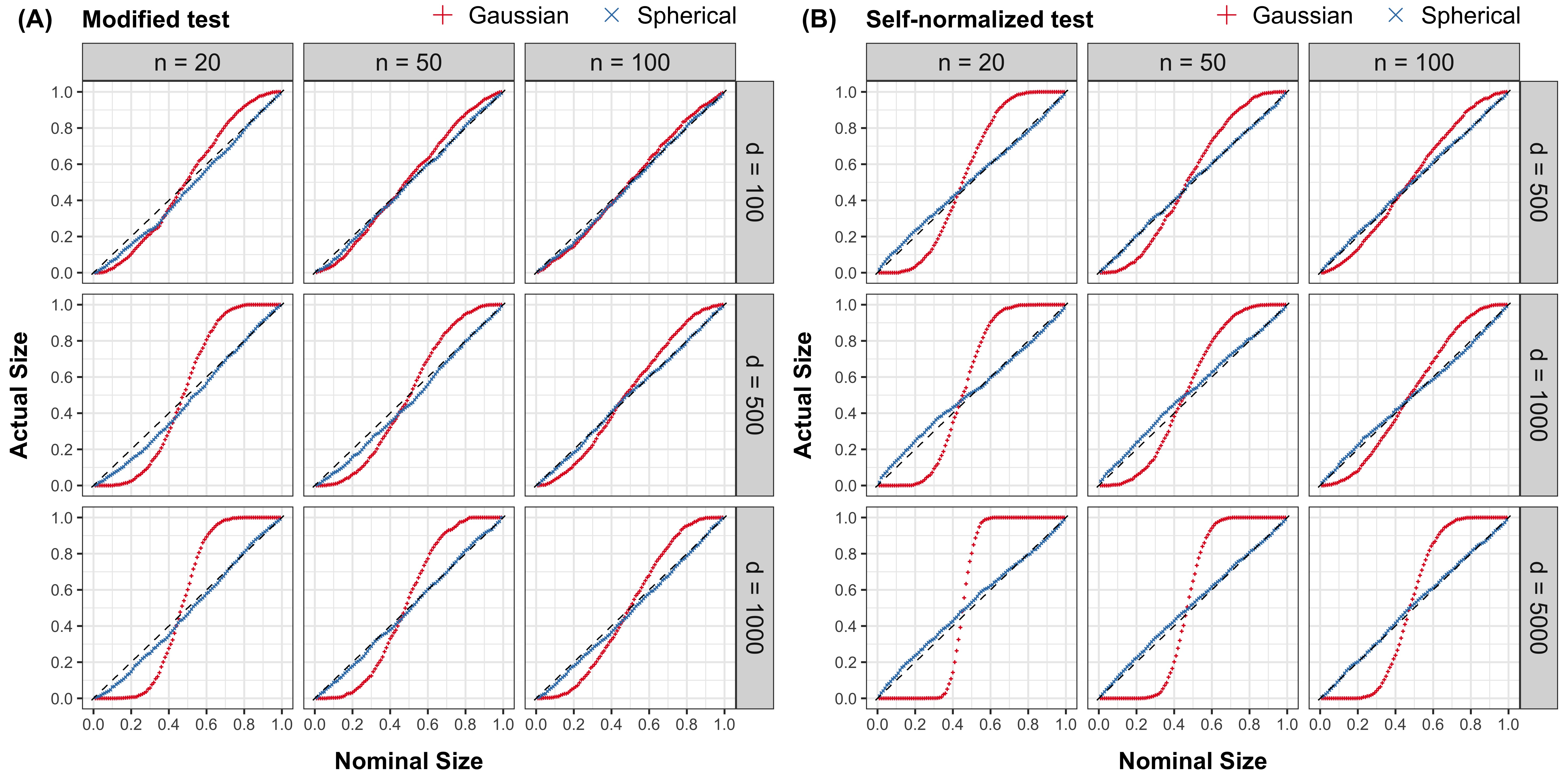}
	\caption{Panel (A) shows QQ-Plots of `nominal' and `actual' sizes of the Gaussian and spherical modified test bootstrap Section~\ref{subsec:PowerEnhancement} under Toeplitz design (DGP 2). Panel (B) shows shows QQ-Plots of `nominal' and `actual' sizes of the Gaussian and spherical self-normalized test from Section~\ref{subsec:EllipticallyDistributedData} for data with $t(4)$-distribution (DGP 4).}
	\label{fig:3}
\end{figure}

In Figure~\ref{fig:3} we show results for the Modified bootstrap test (Section~\ref{subsec:PowerEnhancement}) under DGP 2 and the self-normalized test (Section~\ref{subsec:EllipticallyDistributedData}) under DGP 4. In both cases, we observe that the spherical bootstrap test is more accurate than the Gaussian test. As predicted by our theory, the validity of the (spherical) self-normalized test is virtually unaffected by dimension. Notice that the fact that the spherical bootstrap performs so well for the self-normalized test is not entirely trivial: while each normalized vector $\widetilde{X}_i$ follows a spherical distribution, the (rescaled) sum of of spherical random vectors is not spherically distributed.

\begin{figure}[p]
	\includegraphics[scale=0.08]{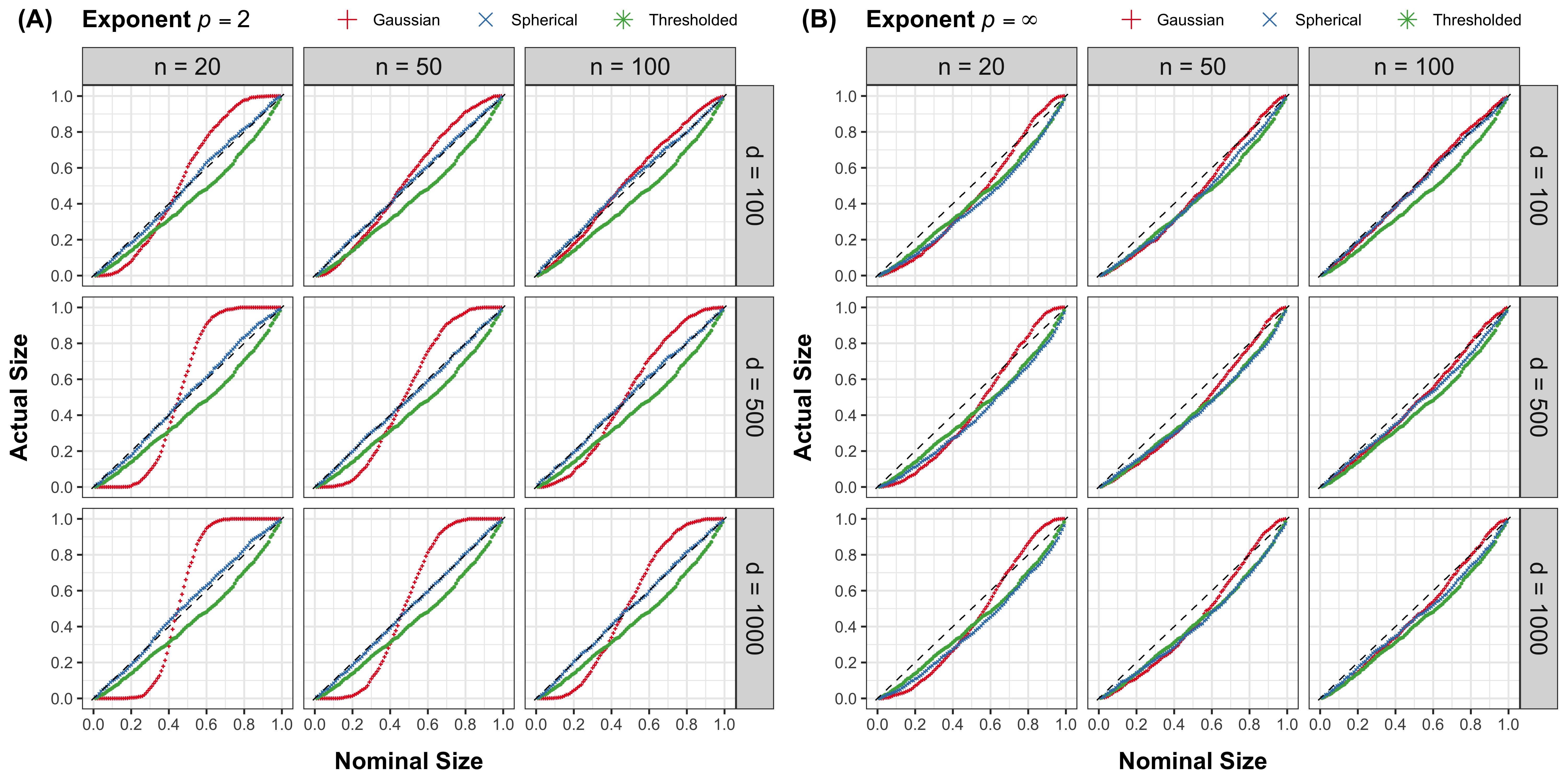}
	\caption{QQ-Plots of `nominal' and `actual' sizes of the Gaussian bootstrap test, the Gaussian bootstrap test with thresholded covariance matrix, and the spherical bootstrap test for banded covariance design (DGP 3). Panel (A) shows simulation results for the bootstrap test with exponent $p=2$, panel (B) shows simulation results for the bootstrap test with exponent $p=\infty$. }
	\label{fig:4}
\end{figure}

In Figure~\ref{fig:4} we compare the accuracy of the basic Gaussian bootstrap test, the bias-corrected spherical bootstrap test, and the Gaussian bootstrap test under DGP 3 with thresholded covariance matrix (Section~\ref{subsec:CovarianceEstimation}) for exponents $p \in \{2, \infty\}$. In accordance with our theory, the thresholded covariance matrix does indeed reduce the bias. However, it results in a conservative test whose actual level is consistently less than the targeted nominal level. Among the three tests, the spherical bootstrap test is still the most accurate one.

\section{Discussion}\label{sec:Discussion}
In Section~\ref{subsec:GlobalNull} we raised four questions about the bootstrap hypothesis test that motivated us to write this paper. Synthesizing our findings from Sections~\ref{sec:BootstrapTest}-\ref{sec:NumericalStudies} we can now answer these questions:
\begin{itemize}
	\item \emph{On the choice of exponent $p \geq 1$ and the power of the test:} Our analysis indicates that exponents $p \in \{2, \log t\}$ are the most relevant ones, with $\varphi_\alpha(T_{n,2},\widehat{\Omega}_n)$ having high power against dense alternatives and $\varphi_\alpha(T_{n,\log t},\widehat{\Omega}_n)$ high power against sparse alternatives. Since $\ell_\infty$- and $\ell_{\log t}$-norms are equivalent for $t$-dimensional vectors, this result is in line with the literature~\citep[][]{cai2014two, fan2015power, he2021Asymptotic}; however, for high-dimensional bootstrap tests this question was unsettled until now. Since the high-dimensional bootstrap tests do not satisfy a CLT (but only the weaker Gaussian approximability property, see Appendices~\ref{subsec:GeneralResults} and~\ref{subsec:ResultsEPT}), the power analysis is based on the Bahadur slope of the tests.
	\item \emph{Leveraging distributional assumptions:} Our discussion of the self-normalized bootstrap test for elliptically distributed data shows that it is indeed possible to exploit distributional assumptions in way that dramatically improves theoretical and empirical performance of the basic bootstrap test. If one accepts the stronger distributional assumption, then the additional sufficient condition for type 1 error control and consistency is extremely mild and, in many cases, holds for arbitrarily large dimensions $d,t$. It might well be possible to develop similar tests for other data generating processes~\citep[see][]{wang2015high-dimensional}.
	\item \emph{The effect of the estimated covariance matrix:} The effect of the estimated covariance matrix on the validity of the bootstrap is fully captured in Assumption~\ref{assumption:ConsistentCovariance}. Intuitively, this assumption says that the estimation error $\|\widehat{\Omega}_n - \Omega\|_{q \rightarrow p}$ has to be asymptotically negligible compared to the variance of the Gaussian proxy statistic $\|Z\|_p$. If the dimensions $d,t$ are fixed, $\mathrm{Var}(\|Z\|_p)$ plays no role. However, typically, $\mathrm{Var}(\|Z\|_p)\downarrow 0$ whenever $d , t \rightarrow \infty$ and $p > 2$ (see Appendix~\ref{subsec:LowerBoundVariance}); hence, Assumption~\ref{assumption:ConsistentCovariance} can be difficult to satisfy in high dimensions.
	
	\item \emph{Utilizing insights from high-dimensional probability theory to mitigate biases:} The bias decomposition of the (Monte Carlo) bootstrap test exposes three sources of potential bias: a Gaussian approximation error, a Gaussian comparison error, and a Monte Carlo error. The first error is unavoidable and can only be reduced by collecting more samples $n$. The second error can be mitigated by sampling from Gaussian proxy statistics $T_{n,p}^*$ that employ structured covariance matrix estimators to leverage low rank, bandedness, or approximate sparsity of the population covariance.
	The third error can be attenuated by improving the efficiency of the Monte Carlo sampling procedure. We propose the spherical bootstrap which is motivated by hypercontractivity properties of $\ell_p$-norms of high-dimensional spherically distributed random vectors. Simulation studies show that both bias correction schemes tend to produce bootstrap tests that have better control of type 1 error. 
\end{itemize}
We conclude with several comments and questions that arose while working on this paper:
\begin{itemize}
	\item \emph{How much of the theory carries over to statistics other than $\ell_p$-norms?} The practical power of the bootstrap principle is its ability to perform inference in complicated settings involving highly non-linear statistics. While $\ell_p$-norms are non-linear functions, they can also be expressed as simple suprema over their respective dual-norm unit balls. This duality relationship allows us to apply the general theory from~\cite{giessing2022anticoncentration,giessing2023Gaussian} and drives the majority of the theoretical results. This general theory can be applied to a virtually any statistic and function class. However, it remains a major challenge to derive sharp upper bounds on the third moments of arbitrary statistics and non-trivial lower bounds on the variances of the corresponding Gaussian proxy statistics.
	
	\item \emph{Why do all rates (in this paper) depend on $\mathrm{Var}(\|Z\|_p)$?} The fact that all rates in this paper (Appendices~\ref{subsec:GeneralResults} and~\ref{subsec:ResultsEPT}) feature the variance of the Gaussian proxy statistic $\|Z\|_p$ is intriguing and deserves further investigation. There appears to be a connection between ``super-concentration''~\citep[][]{chatterjee2014Superconcentration}, ``anti-concentration''~\citep[][]{lecam1986asymptotic}, and ``universality laws'' (e.g. classical CLTs, high-dimensional Gaussian approximation results, or~\possessivecite{mossel2010noise} multilinear invariance principle). For a brief discussion we refer to Appendix~\ref{subsec:LowerBoundVariance}.
	
	\item \emph{What are necessary conditions for bootstrap consistency in high dimensions?} We know that in the classical, low-dimensional setting the (empirical or wild) bootstrap is valid if and only if the involved function classes are Donsker~\citep{mammen1992bootstrap, vandervaart1996weak}. This and the bias decomposition in Section~\ref{subsec:HeuristicsBias} seem to suggest that Gaussian approximability of the function classes might not only be a sufficient but also a necessary condition for bootstrap consistency  in high-dimensional settings. However, this need not be true. In fact, the proofs in the theoretical companion papers~\cite{giessing2022anticoncentration, giessing2023Gaussian} use Gaussianity only twice, once when applying an anti-concentration argument and a second time when invoking Stein's lemma. Other distributions (might) satisfy similar properties that can be used to replace the original arguments.
		
	\item \emph{What are the connections between Efron's empirical, the wild, and the spherical bootstrap?} Numerical results from simulations (not included in this paper) show that Efron's empirical bootstrap is valid for $\ell_p$-statistics with large exponents $p \gtrsim \log d$ but fails for small exponents $p \lesssim \log d$. This is broadly in line with results by~\cite{chernozhukov2013GaussianApproxVec} and~\cite{elKaroui2018can}. We can contribute the following insight: Recall that Efron's empirical bootstrap can be interpreted as a multiplier bootstrap with multinomial multipliers~\citep[][p.345]{vandervaart1996weak}. The multipliers are thus negatively correlated and have correlation coefficients $-n^{-1}$. When bootstrapping the $\ell_p$-norm with $p \gtrsim \log d$ this results in an additional bias term of order $O(n^{-1})$ (modulo $(\log d)$-terms); whereas, when bootstrapping the $\ell_p$-norm with $p \lesssim \log d$ this yields a non-vanishing bias term of order $O(1)$.~\cite{elKaroui2018can} arrive at a similar conclusion and therefore advocate using uncorrelated multipliers in high dimensions. Our spherical bootstrap can now be viewed as just another multiplier bootstrap with uncorrelated multipliers. However, unlike the multipliers proposed by~\cite{elKaroui2018can} ours are dependent which, together with the hypercontractivity of $\ell_p$-norms, improves the efficiency and accuracy of the bootstrap procedure. Numerical results (not included in this paper) show that even in classical, low-dimensional settings the spherical bootstrap tends to outperform Efron's empirical bootstrap. Existing theory on exchangeable multiplier bootstrap procedures does not cover the spherical bootstrap because it requires non-negative multipliers~\citep[][p. 353ff]{vandervaart1996weak}. In this paper we have focused on validity and consistency of the spherical bootstrap procedure; a more comprehensive and comparative study of its efficiency and accuracy would be highly desirable.
\end{itemize}

\newpage
\section*{Acknowledgement}
Alexander Giessing's research is supported by NSF grant DMS-2310578, Jianqing Fan’s research by ONR grant N00014-19-1-2120 and the NSF grants DMS-2052926, DMS-2053832, and DMS-2210833.

\newpage
\normalsize
\setcounter{page}{1}

\bibliography{GBA_51_ref}

\newpage

\title{\textsc{Supplementary Materials for ``A bootstrap hypothesis test for high-dimensional mean vectors''}}
\author{Alexander Giessing\protect\footnotemark[1] \and Jianqing Fan\protect\footnotemark[2]}

\date{\today}

\maketitle

\appendix

\setcounter{page}{1}

\noindent{\bf\LARGE Contents}

\startcontents[sections]
\printcontents[sections]{l}{1}{\setcounter{tocdepth}{2}}

\section{Bootstrapping high-dimensional $\ell_p$-statistics}\label{sec:TheoryBootstrap}
In this section we present Gaussian and bootstrap approximation for high-dimensional $\ell_p$-statistics (Appendix~\ref{subsec:GeneralResults}), lower bounds on the variance of $\ell_p$-norms of (an-isotropic) Gaussian random vectors (Appendix~\ref{subsec:LowerBoundVariance}), and essential results from our companion papers~\cite{giessing2022anticoncentration, giessing2023Gaussian} (Appendix~\ref{subsec:ResultsEPT}).~ 

\subsection{Approximating sampling distribution and quantiles}\label{subsec:GeneralResults}
Unless otherwise stated, all parameters in this section may be thought of as indexed by the sample size $n$. Since we present non-asymptotic bounds that hold for all $n \geq 1$ we leave this dependence implicit. 
The results in this section are consequences of the more general results in our two companion papers~\cite{giessing2022anticoncentration, giessing2023Gaussian} summarized in Appendix~\ref{subsec:ResultsEPT}.

At the core of any bootstrap procedure is a CLT. Thus, the following theorem is key to all theoretical results in this paper:

\begin{theorem}[Gaussian approximation]\label{theorem:CLT-Lp-Norm}
	Let $X, X_1, \ldots X_n \in \mathbb{R}^d$ be i.i.d. random vectors with mean zero and positive semi-definite covariance matrix $\Sigma \neq \mathbf{0} \in \mathbb{R}^{d \times d}$. Set $S_n = n^{-1/2} \sum_{i=1}^n X_i$ and $Z \sim N(0, \Sigma)$. Then, for $1 \leq p \leq \infty$, $n \geq 1$,
	\begin{align*}
		&\sup_{s \geq 0} \Big|\mathrm{P}\left(\|S_n\|_p \leq s \right) -\mathrm{P}\left(\|Z\|_p \leq s \right) \Big|\\
		&\quad{}\quad{}\quad{} \lesssim \frac{(\mathrm{E}[\|X \|_p^3])^{1/3}}{n^{1/6}\sqrt{\mathrm{Var}(\|Z\|_p)}} + \frac{\mathrm{E} \left[ \|X\|_p^3 \mathbf{1}\{\|X\|_p^3 > n\: \mathrm{E}[\|X \|_p^3]\}\right]}{\mathrm{E} \left[ \|X\|_p^3\right]}  + \frac{ \mathrm{E}[\|Z\|_p]}{\sqrt{n\mathrm{Var}(\|Z\|_p)}},
	\end{align*}
	where $\lesssim$ hides an absolute constant independent of $p, n, d$, and the distribution of the $X_i$'s.
\end{theorem}
This Gaussian approximation differs from the ones in the literature~\citep[e.g.][]{chernozhukov2013GaussianApproxVec, chernozhukov2014AntiConfidenceBands, chernozhukov2014GaussianApproxEmp, chernozhukov2015ComparisonAnti, chernozhukov2017CLTHighDim, chernozhukov2019ImprovedCLT, chernozhukov2021NearlyOptimalCLT, zhang2017gaussian, xu2019L2Asymptotics, lopes2020bootstrapping, xue2020Distribution, fang2021HighDimCLT} in three ways: First, it holds for all $\ell_p$-norms with exponents $p \in [1, \infty]$ whereas existing results apply only to the $\ell_\infty$-norm or the square of the $\ell_2$-norm. Second, the non-asymptotic upper bound is dimension-free, i.e. it does not explicitly depend on the dimension $d$. Since the quantities $\mathrm{E}[\|X\|_p^3]$ and $\mathrm{E}[\|Z\|_p]$ can often be upper bounded in terms of the trace of the covariance matrix (e.g. see Lemmas~\ref{lemma:example:SubGaussianData}--\ref{lemma:example:HeavyTailedData}), this opens the possibility of leveraging the eigen-structure of the covariance matrix. Third, the result applies to degenerate distributions that do not have a strictly positive definite covariance matrix.

The most striking feature of this Gaussian approximation result is the dependence of the upper bound on the variance of the Gaussian proxy statistic $\mathrm{Var}(\|Z\|_p)$. 
The magnitude of this variance depends in a delicate way on the exponent $p$ and the dimensions $d$ and hints at an interesting connection between universality laws for high-dimensional data, anti-concentration, and super-concentration. We refer to Section~\ref{subsec:LowerBoundVariance} for details.

The dependence of the upper bounds on the sample size $n$ is most likely sub-optimal. For example, for $p =\infty$ and under additional assumptions on the moments of the data, sharper rates have been derived in~\cite{lopes2020bootstrapping} and~\cite{chernozhukov2021NearlyOptimalCLT}. 

Specializing to the $\ell_2$-norm and comparing the upper bound in Theorem~\ref{theorem:CLT-Lp-Norm} with the upper bound in Theorem 1.1 in~\cite{bentkus2003DependenceBerryEsseen} we observe that our bound is better by a factor $d^{1/4}$. The reason for this discrepancy is that Gaussian approximation result is in fact a weaker a statement than the~\possessivecite{bentkus2003DependenceBerryEsseen} Berry-Esseen-type result and the two results are incomparable: In Theorem~\ref{theorem:CLT-Lp-Norm} we consider the supremum over all $\ell_2$-balls with center at the origin, whereas in Theorem 1.1~\cite{bentkus2003DependenceBerryEsseen} considers the supremum over all $\ell_2$-balls.

Since the Gaussian distribution is fully characterized by its first two moments, Theorem~\ref{theorem:CLT-Lp-Norm} instantly suggests that it should be possible to approximate the sampling distribution of $\|S_n\|_p$ with the sampling distribution of $\|Z_n\|_p$, where $Z_n \mid X_1, \ldots, X_n \sim N(0, \widehat{\Sigma}_n)$ and $\widehat{\Sigma}_n$ is a positive semi-definite estimate of $\Sigma$. The next theorem formalizes this idea:

\begin{theorem}[Gaussian process bootstrap approximation]\label{theorem:Bootstrap-Lp-Norm}
	Let $X, X_1, \ldots X_n \in \mathbb{R}^d$ be i.i.d. random vectors with mean zero and positive semi-definite covariance matrix $\Sigma \neq \mathbf{0} \in \mathbb{R}^{d \times d}$. Let $\widehat{\Sigma}_n \equiv \widehat{\Sigma}_n(X_1, \ldots, X_n)$ be any positive semi-definite estimate of $\Sigma$. Set $S_n = n^{-1/2} \sum_{i=1}^n X_i$, $Z \sim N(0, \Sigma)$, and $Z_n \mid X_1, \ldots X_n \sim N(0,\widehat{\Sigma}_n )$. Then, for $1 \leq p \leq \infty$, $n \geq 1$,
	\begin{align*}
		&\sup_{s \geq 0} \Big|\mathrm{P}\left(\|S_n\|_p \leq s \right) -\mathrm{P}\left(\|Z_n\|_p \leq s \mid X_1, \ldots, X_n \right) \Big|\\
		&\quad{}\lesssim \frac{(\mathrm{E}[\|X \|_p^3])^{1/3}}{n^{1/6}\sqrt{\mathrm{Var}(\|Z\|_p)}} + \frac{\mathrm{E} \left[ \|X\|_p^3 \mathbf{1}\{\|X\|_p^3 > n\: \mathrm{E}[\|X \|_p^3]\}\right]}{\mathrm{E} \left[ \|X\|_p^3\right]} + \frac{ \mathrm{E}[\|Z\|_p]}{\sqrt{n\mathrm{Var}(\|Z\|_p)}} + \left( \frac{\|\widehat{\Sigma}_n- \Sigma\|_{q\rightarrow p}}{ \mathrm{Var}(\|Z\|_p)} \right)^{1/3},
	\end{align*}
	where $1/p+ 1/q = 1$ and $\lesssim$ hides an absolute constant independent of $p, n, d$, and the distribution of the $X_i$'s.
\end{theorem}

To test statistical hypotheses and construct confidence intervals, we do not need to estimate the entire sampling distribution but (only) specific quantiles. Since the covariance matrix $\Sigma$ is unknown, the quantiles of $\|Z\|_p$ with $Z \sim N(0, \Sigma)$ are infeasible. Instead, we use quantiles of $\|Z_n\|_p$ where $Z_n \mid X_1, \ldots, X_n \sim N(0, \widehat{\Sigma}_n)$: For $\alpha \in (0,1)$ arbitrary, we define 
\begin{align*}
	c_p^*(\alpha; \widehat{\Sigma}_n) := \inf \left\{s \geq 0: \mathrm{P}\left(\|Z_n\|_p\leq s \mid X_1, \ldots, X_n \right)  \geq \alpha \right\}.
\end{align*}
Since this quantity is random, it is not immediately obvious that it is a valid approximation of the $\alpha$-quantile of the sampling distribution of $\|S_n\|_p$. However, combing Theorem~\ref{theorem:Bootstrap-Lp-Norm} with (by now) standard arguments~\citep[e.g.][]{chernozhukov2013GaussianApproxVec}, we obtain the following result:

\begin{theorem}[Gaussian process bootstrap approximation of quantiles]\label{theorem:Bootstrap-Lp-Norm-Quantiles}
	Recall the conditions of Theorem~\ref{theorem:Bootstrap-Lp-Norm}. Let $(\Theta_n)_{n \geq 1} \in \mathbb{R}$ be a sequence of arbitrary random variables, not necessarily independent of $X_1, \ldots, X_n$. For $1 \leq p \leq \infty$, $n \geq 1$,
	\begin{align*}
		&\sup_{\alpha \in (0,1)} \Big|\mathrm{P}\left(\|S_n\|_p + \Theta_n \leq c_p^*(\alpha; \widehat{\Sigma}_n ) \right)  - \alpha \Big| \\
		&\quad{}\quad{}\lesssim \frac{(\mathrm{E}[\|X\|_p^3])^{1/3}}{n^{1/6}\sqrt{\mathrm{Var}(\| Z\|_p)}} + \frac{\mathrm{E} \left[ \|X\|_p^3 \mathbf{1}\{\|X\|_p^3 > n\: \mathrm{E}[\|X\|_p^3]\}\right]}{\mathrm{E}\left[ \|X\|_p^3\right]}+  \frac{ \mathrm{E}[\|Z\|_p]}{\sqrt{n\mathrm{Var}(\|Z\|_p)}}\\
		& \quad{}\quad{}\quad{}\quad{}+ \inf_{\delta > 0}\left\{ \left(\frac{\delta }{ \mathrm{Var}( \|Z\|_p)} \right)^{1/3}  + \mathrm{P}\left( \|\widehat{\Sigma}_n- \Sigma\|_{q \rightarrow p} > \delta\right)\right\} \\
		&\quad{}\quad{}\quad{} \quad{} \quad{}+ \inf_{\eta > 0} \left\{ \frac{\eta }{ \sqrt{\mathrm{Var}( \|Z\|_p)} } + \mathrm{P}\left(|\Theta_n| > \eta \right)\right\},
	\end{align*}
	where $1/p+ 1/q = 1$ and $\lesssim$ hides an absolute constant independent of $p, n, d$, and the distributions of the $X_i$'s and the $\Theta_n$'s.
\end{theorem}

\begin{remark}[On the name \emph{Gaussian process bootstrap}]
	Naming the approximations in Theorems~\ref{theorem:Bootstrap-Lp-Norm} and~\ref{theorem:Bootstrap-Lp-Norm-Quantiles} ``Gaussian process bootstrap'' may appear overdone since the approximations involve only a single finite-dimensional Gaussian random vector. The name originates from the more general bootstrap procedure presented in Section~\ref{subsec:ResultsEPT}. We therefore decided to keep.
\end{remark}

\subsection{Lower bounds on the variance of the Gaussian proxy statistic}\label{subsec:LowerBoundVariance}
Theorems~\ref{theorem:CLT-Lp-Norm}--\ref{theorem:Bootstrap-Lp-Norm-Quantiles} from the preceding section depend on the moments of $\|X\|_p$ and $\|Z\|_p$ and the variance of $\|Z\|_p$. Bounding moments of $\ell_p$-norms of (high-dimensional) random vectors is easy. The reader can find some useful results and pointers to the literature in Section~\ref{subsec:Miscellaneous} and (the proofs of) Lemmas~\ref{lemma:example:SubGaussianData}--\ref{lemma:example:HeavyTailedData}. In contrast, obtaining (lower) bounds on the variance of $\ell_p$-norms of (high-dimensional) Gaussian random vectors is extremely challenging. Yet, lower bounds on the variance are precisely what is needed for Theorems~\ref{theorem:CLT-Lp-Norm}--\ref{theorem:Bootstrap-Lp-Norm-Quantiles} to be useful.

The lower bounds that we present in this section are not necessarily optimal in all scenarios. However, we have tried hard to derive lower bound that are statistically meaningful and easy to interpret. In particular, we have deliberately avoided arguments that result in lower bounds that depend on the eigenvalues of the covariance matrix (an exception are expressions that depend on the trace).

We begin with a review of lower bounds on the variance of $\ell_p$-norms for standard Gaussian random vectors with isotropic covariance matrix. While this case is of limited use in statistical applications, it is nonetheless an important reference point that guides our (partial) analysis of the general an-isotropic case below.

\begin{proposition}[Theorem A,~\citeauthor{lytova2017variance},~\citeyear{lytova2017variance}]\label{prop:Lyotova2017variance}
	Let $Z \in \mathbb{R}^d$ be a standard normal random vector with mean zero and identity covariance matrix. Denote by $\xi$ the $1- 1/d$ quantile of $|Z_1|$, i.e. $\mathrm{P}(|Z_1| \leq \xi ) = 1- 1/d$. Then, there exists an absolute constant $d_0 > 0$ such that for all $d > d_0$,
	\begin{itemize}
		\item[(i)] for all $p \in \left[1, \frac{2 \log d}{\log(2e)}\right]$, 
		\begin{align*}
			\mathrm{Var}(\|Z\|_p) \asymp \frac{2^p}{p} d^{2/p-1};
		\end{align*}
		\item[(ii)] for all $p \in \left(\frac{2 \log d}{\log(2e)}, \xi^2\right)$,
		\begin{align*}
			\mathrm{Var}(\|Z\|_p) \asymp \frac{\exp\left(-\frac{p}{2e}d^{2/p} + \log d\right)}{\sqrt{\log d} \left(\sqrt{\log d} + p - \frac{2 \log d}{\log(2e)}\right)};
		\end{align*}
		\item[(iii)] for all $p \in [\xi^2, \infty]$,
		\begin{align*}
			\mathrm{Var}(\|Z\|_p) \asymp \frac{1}{\log d}\left(1 - \frac{\xi^2 - \xi}{p} \right).
		\end{align*}
	\end{itemize}
\end{proposition}
Notice that case (i) really contains three qualitatively different cases: for $p \in [1,2)$ the variance $\mathrm{Var}(\|Z\|_p)$ increases in the dimension, for $p=2$ the variance is independent of the dimension, and for $p \in (2, 2 \log d / \log(2e)]$ the variance decreases polynomially fast in the dimension. This implies that if the data has identity covariance matrix, $\ell_p$-statistics with $p \in (2, 2 \log d / \log(2e)]$ are Gaussian approximable (in the sense of Theorem~\ref{theorem:CLT-Lp-Norm}) only under very stringent growth conditions on $d$ and $n$ which are often incompatible with high-dimensional settings $d \gg n$. This collapse of the variance is an instance of the so-called ``super-concentration phenomenon''~\citep[e.g.][]{chatterjee2014Superconcentration}. Combined with Theorem~\ref{theorem:CLT-Lp-Norm} we conjecture that super-concentration precludes universality laws (such as high-dimensional CLTs). A rigorous proof of this claim would require a matching lower bound on the Kologorov-Smirnov distance in Theorem~\ref{theorem:CLT-Lp-Norm}. Alternatively, since by Lemma~\ref{lemma:AntiConcentration-SeparableProcess} Gaussian processes are either super-concentrated or anti-concentrated, the claim would follow if we made rigorous~\possessivecite{lecam1986asymptotic} heuristic connection between CLTs and anti-concentration inequalities. We leave these open problems to future research. For the present paper these observations suggest that the Gaussian process bootstrap procedure for $\ell_p$-norms with small exponents $p\in (2, 2\log d/\log(2e)]$ can only work if the data are not too independent and/ or too noisy, i.e. the data have to concentrate in some sense on a low-dimensional subspace (be it low-rank covariance matrix, variance decay, or bounded effective rank). This is the guiding principle behind the conditions and bounds in Lemmas~\ref{lemma:example:SubGaussianData}--\ref{lemma:example:HeavyTailedData}.

We are now ready to state the main result of this section: lower bounds on $\mathrm{Var}(\|Z\|_p)$ when $Z \in \mathbb{R}^d$ is a Gaussian random vector with an-isotropic covariance matrices $\Sigma \neq I_d$. We emphasize that this result is not a simple generalization of Proposition~\ref{prop:Lyotova2017variance}. The proof of Proposition~\ref{prop:Lyotova2017variance} relies heavily on the isotropy of the Gaussian random vector and cannot be extended to the an-isotropic case. Moreover, our lower bounds are tighter than those that would follow from direct applications of the general results in~\cite{cacoullus1982upper} and~\cite{houdre1995covariance}.

\begin{theorem}\label{theorem:LowerBoundVariance-LpNorm}
	Let $Z \in \mathbb{R}^d$ be a centered Gaussian random vector with marginals $Z_k \sim N(0, \sigma_k^2)$ and $\sigma_{(1)}^2 \leq \ldots \leq \sigma_{(d)}^2$ be the ordered values of $\sigma_1^2, \ldots, \sigma_d^2$. Further, let $|\mathrm{Corr}(Z_j, Z_k)| \leq \rho$ for some $\rho \in [0,1]$ and all $1 \leq j, k \leq d$. Then, for all $d \geq 1$,
	\begin{itemize}
		\item[(i)] for all $p \in [1,2]$,
		\begin{align*}
			\mathrm{Var}(\|Z\|_p) \geq \frac{\pi}{3^2} \left( \frac{1}{d} \sum_{i=1}^d \sigma_i^p\right)^{2/p}  d^{2/p-1};
		\end{align*}
		\item[(ii)] for all $p \in (2, 2\log d)$,
		\begin{align*}
			\mathrm{Var}(\|Z\|_p) \geq \frac{\pi}{6} \frac{p^2}{2^{3p}}  \left( \frac{1}{d} \sum_{i=1}^d \sigma_i^{2p}\right)^{1/p} d^{2/p-1};
		\end{align*}
		\item[(iii)] for all $p \in [2 \log d, \infty)$,
		\begin{align*}
			\mathrm{Var}\left(\|Z\|_p\right)  \geq \frac{1}{15^2}\left(\frac{\sigma_{(1)}^2}{\sigma_{(1)} + \sigma_{(d)} \sqrt{\log d}}\right)^2 \frac{1}{11^6}\left(\frac{1}{ 1  + \rho \sqrt{\log d} }\right)^2;
		\end{align*}
		\item[(iv)] for $p=\infty$,
		\begin{align*}
			\mathrm{Var}(\|Z\|_\infty) \geq  \frac{1}{15^2}\left(\frac{\sigma_{(1)}^2}{\sigma_{(1)} + \sigma_{(d)} \sqrt{\log d}}\right)^2.
		\end{align*}
	\end{itemize}
\end{theorem}
Clearly, above result is not as complete as the one in Proposition~\ref{prop:Lyotova2017variance}. However, it is sufficient for a wide range of statistical applications. Moreover, one easily verifies that if $Z$ is isotropic the lower bounds are of the correct order in most cases when compared to Proposition~\ref{prop:Lyotova2017variance} and results in~\cite{biau2015HighDimPNorms}. The notable exception is the regime of $p \in (2 \log d / \log(2e),\infty)$.~\cite{lytova2017variance} show that $\xi^2 = 2\log d - o(\sqrt{\log d})$. Therefore, cases (ii) and (iii) of Proposition~\ref{prop:Lyotova2017variance} imply that for any $\delta \in (0,1)$ and for all $d \geq d_0$, there exists a constant $c_\delta > 0$ (depending on $\delta$ only) such that $\mathrm{Var}(\|Z\|_{(2 -\delta)\log d}) \leq  d^{-c_\delta}$ but also $\mathrm{Var}(\|Z\|_{(2 + \delta)\log d}) \geq  \frac{c_\delta}{\log d}$. Thus, the variances of $\ell_p$-norms of isotropic Gaussian random vectors can exhibit qualitatively dramatically different behavior across equivalent $\ell_p$-norms for $p \in (2 \log d / \log(2e),\infty)$. Our lower bounds for an-isotropic Gaussian random vectors do not reflect this. Though, interestingly, our lower bounds depend on the (maximum) correlation between the entries of the Gaussian random vector. 

The following refinement of above result for $p \in \{1, 2, \infty\}$ is useful in applications. The lower bound for $p =2$ is (up to an absolute constant) optimal~\citep[e.g.][]{valettas2019tightness}.

\begin{theorem}[Refinement of Theorem~\ref{theorem:LowerBoundVariance-LpNorm}]\label{theorem:LowerBoundVariance-LpNorm-Refinements}
	Recall assumptions and notation from Theorem~\ref{theorem:LowerBoundVariance-LpNorm}. Further, let $\Sigma$ be the covariance matrix of the Gaussian random vector $Z \in \mathbb{R}^d$. We have the following refinements: For all $d \geq 1$,
	\begin{itemize}
		\item[(i)] for $p=1$,
		\begin{align*}
			\mathrm{Var}(\|Z\|_1) \geq \frac{\pi}{2} \mathrm{tr}(\Sigma);
		\end{align*}
		\item[(ii)] for $p=2$,
		\begin{align*}
			\mathrm{Var}(\|Z\|_2) \geq \frac{\mathrm{tr}(\Sigma^2)}{\mathrm{tr}(\Sigma)} \vee \frac{\sum_{i=1}^d \sigma_i^4}{\sum_{i=1}^d \sigma_i^2};
		\end{align*}
		\item[(iii)] for $p = \infty$,
		\begin{align*}
			\mathrm{Var}(\|Z\|_\infty) \geq  \frac{1}{12}\left(\frac{\bar{\sigma}_d}{1 + \sqrt{\log d}}\right)^2,
		\end{align*}
		where $\bar{\sigma}_d = \left(1 + \sqrt{\log d}\right) / \left(1/ \sigma_{(1)} + \max_{1 \leq k \leq d} \left(1 + \sqrt{\log k}\right)/\sigma_{(k)}\right)$. 
	\end{itemize}
\end{theorem}

We conclude this section with three auxiliary results which, when combined, imply Theorem~\ref{theorem:LowerBoundVariance-LpNorm}.

\begin{lemma}\label{lemma:LowerBoundVariance-LpNorm-2}
	Let $Z \in \mathbb{R}^d$ be a centered Gaussian random vector. 
	\begin{align*}
		\mathrm{Var}(\|Z\|_p) \geq \frac{\sum_{k=1}^d\mathrm{Var}(Z_k)^p}{\mathrm{E}[\|Z\|_p^{2p-2}]} \times 
		\begin{cases}
			\frac{2^p}{p^2}\Gamma^2\left(\frac{p}{2} + 1\right), \quad{} & \mathrm{if} \quad{} p \in [1,2]\\
			\frac{2^{-p}}{3} \Gamma^2\left(\frac{p}{2} + 1\right), \quad{} & \mathrm{if} \quad{} p \in (2, \infty).
		\end{cases}
	\end{align*}
\end{lemma}

\begin{lemma}\label{lemma:LowerBoundVariance-LpNorm-6}
	Let $Z \in \mathbb{R}^d$ be a centered Gaussian random vector with $0 < \sigma_{(1)}^2 \leq \mathrm{E}[Z_k^2] \leq \sigma_{(d)}^2 < \infty$ and $|\mathrm{Corr}(Z_j, Z_k)| \leq \rho$ for all $1 \leq j, k \leq d$. For all $p \geq 2 \log d$,
	\begin{align*}
		\mathrm{Var}\left(\|Z\|_p\right)  \geq \left(\frac{\sigma_{(1)}^2}{\sigma_{(1)} + \mathrm{E}[\|Z\|_p]}\right)^2\frac{e^{-3} 48^{-4}}{1 + \log(2d) \rho}.
	\end{align*}
\end{lemma}

\begin{lemma}\label{lemma:LowerBoundVariance-LpNorm-3}
	Let $Z \in \mathbb{R}^d$ be a centered Gaussian random vector with marginals $Z_k \sim N(0, \sigma_k^2)$ and $\sigma_{(1)}^2 \leq \ldots \leq \sigma_{(d)}^2$ be the ordered values of $\sigma_1^2, \ldots, \sigma_d^2$. 
	\begin{align*}
		\mathrm{Var}(\|X\|_\infty) \geq \frac{\bar{\sigma}_d^2}{12(2 + 2 \sqrt{\log d})^2} \geq \frac{\sigma_{(1)}^2}{15^2 \log d},
	\end{align*}
	where  $\bar{\sigma}_d = \left(2 + 2 \sqrt{\log d}\right) / \left(1/ \sigma_{(1)} + \max_{1 \leq k \leq d} \left(1 + \sqrt{2\log k}\right)/\sigma_{(k)}\right) \geq \sigma_{(1)}$.
\end{lemma}

\subsection{Gaussian multiplier and Gaussian process bootstrap}\label{subsec:ResultsEPT}
In this section we summarize the key theoretical results on the \emph{Gaussian multiplier bootstrap} and the \emph{Gaussian process bootstrap} from our companion papers~\cite{giessing2022anticoncentration, giessing2023Gaussian} for easy reference. The theoretical results on bootstrapping high-dimensional $\ell_p$-statistics in Section~\ref{subsec:GeneralResults} are simple consequences of these more general results.

In the following, we denote by $X, X_1, X_2, \ldots$ a sequence of i.i.d. random variables taking values in a measurable space $(S, \mathcal{S})$ with common distribution $P$, i.e. $X_i : S^\infty \rightarrow S$, $i \geq 1$, are the coordinate projections of the infinite product probability space $(\Omega, \mathcal{A}, \mathbb{P}) = \left(S^\infty, \mathcal{S}^\infty, P^\infty\right)$ with law $\mathbb{P}_{X_i} = P$. If auxiliary variables independent of the $X_i$'s are involved, the underlying probability space is assumed to be of the form $(\Omega, \mathcal{A}, \mathbb{P}) = \left(S^\infty, \mathcal{S}^\infty, P^\infty\right) \times (Z, \mathcal{Z}, Q)$. We define the empirical measures $P_n$ associated with observations $X_1, \ldots, X_n$ as random measures on $\left(S^\infty,\mathcal{S}^\infty\right)$ given by $P_n(\omega) := n^{-1}\sum_{i=1}^n \delta_{X_i(\omega_i)}$ for all $\omega \in S^\infty$, where $\delta_{x}$ is the Dirac measure at $x$. 

For a class $\mathcal{F}$ of measurable functions from $S$ onto the real line $\mathbb{R}$ we define the \emph{empirical process indexed by $\mathcal{F}$} as
\begin{align*}
	\mathbb{G}_n(f) : = \frac{1}{\sqrt{n}}\sum_{i=1}^n (f(X_i) - P f\big), \quad{} f \in \mathcal{F}.
\end{align*}
Further, we denote by $\{G_P(f) : f \in \mathcal{F}\}$ the \emph{Gaussian $P$-bridge process} with  mean zero and the same covariance function $C_P : \mathcal{F} \times \mathcal{F} \rightarrow \mathbb{R}$ as the process $\{f(X): f \in \mathcal{F}\}$,
\begin{align*}
	(f, g) \mapsto C_P(f, g) := \mathrm{E}[G_P(f) G_P(g)] = Pfg - (Pf)(Pg).
\end{align*}
and by $\{Z_Q(f) : f \in \mathcal{F}\}$ the \emph{Gaussian $Q$-motion} with mean zero and covariance function
\begin{align*}
	(f, g) \mapsto \mathrm{E}[Z_Q(f) Z_Q(g)] = Qfg.
\end{align*}
In passing, we note that, above representations readily imply that Gaussian $P$-bridge processes and Gaussian $Q$-motions are linear in $f \in \mathcal{F}$. (This can be proved in the same way as Theorem 2.1 in~\cite{dudley2014uniform}.)

For probability measures $Q$ on $(S, \mathcal{S})$ we define the $L_q(Q)$-norm, $q \geq 1$, for function $f \in \mathcal{F}$ by $\|f\|_{Q,q} = (Q|f|^q)^{1/q}$ and the $L_2(Q)$-semimetric for functions $f, g \in \mathcal{F}$ by $e_Q(f, g) := \|f- g\|_{Q,2}$. We denote by $L_q(S, \mathcal{S}, Q)$, $q\geq 1$, the space of all real-valued measurable functions $f$ on $(S, \mathcal{S})$ with finite $L_q(Q)$-norm. We write $L_2(\mathcal{F})$ for the space of all linear, square integrable functionals $f^*: \mathcal{F} \rightarrow \mathbb{R}$.

The following theorem provides a non-asymptotic bound on the Kolmogorov distance of the laws of $\|\mathbb{G}_n\|_{\mathcal{F}_n}$ and $\|G_P\|_{\mathcal{F}_n}$.

\begin{lemma}[\citeauthor{giessing2023Gaussian},~\citeyear{giessing2023Gaussian}]\label{lemma:CLT-DimFree}
	For each $n \geq 1$ let $\mathcal{F}_n \subset L_2(S, \mathcal{S}, P)$ be totally bounded w.r.t. a metric $\rho$ and have envelope $F_n \in  L_3(S, \mathcal{S}, P)$. Suppose that there exist functions $\psi_n, \phi_n$ such that 
	\begin{align}\label{eq:corollary:CLT-Max-Norm-Simultaneous-0}
		\mathrm{E} \|G_P\|_{\mathcal{F}_{n,\delta}'} \lesssim 	\psi_n(\delta)\sqrt{\mathrm{Var}(\|G_P\|_{\mathcal{F}_n})} \quad{} \quad{} \text{and} \quad{} \quad{} \big\|\|\mathbb{G}_n\|_{\mathcal{F}_{n,\delta}'} \big\|_{P,1} \lesssim 	\phi_n(\delta)\sqrt{\mathrm{Var}(\|G_P\|_{\mathcal{F}_n})},
	\end{align}
	where $\mathcal{F}_{n, \delta}' = \{ f - g: f, g \in \mathcal{F}_n, \: \rho(f, g) < \delta\|F_n\|_{P,2}\}$. Let $r_n^2 = \inf\big\{  \psi_n(\delta) \vee \phi_n(\delta): \delta > 0\big\}$. Then, for each $n \geq 1$ and $M_n \geq 0$,
	\begin{align*}
		&\sup_{s \geq 0} \Big|\mathbb{P}\left\{\|\mathbb{G}_n\|_{\mathcal{F}_n} \leq s \right\} -\mathbb{P}\left\{\|G_P\|_{\mathcal{F}_n} \leq s \right\} \Big|\\
		&\quad{}\quad{}\quad{} \lesssim \frac{\|F_n\|_{P, 3}}{n^{1/6}\sqrt{\mathrm{Var}(\|G_P\|_{\mathcal{F}_n})}} + \frac{\|F_n\mathbf{1}\{F_n > M_n\}\|_{P,3}^3}{\|F_n\|_{P,3}^3} + \frac{ \mathrm{E} \|G_P\|_{\mathcal{F}_n} + M_n}{\sqrt{n\mathrm{Var}(\|G_P\|_{\mathcal{F}_n})}} + r_n,
	\end{align*}
	where $\lesssim$ hides an absolute constant independent of $n, r_n, \mathcal{F}_n, F_n, M_n, \psi_n, \phi_n$, and $P$.
\end{lemma}
For a thorough discussion of this lemma we refer to~\cite{giessing2023Gaussian}. Here, we focus on the following important consequence: The result implies that in order to approximate the sampling distribution of $\|\mathbb{G}_n\|_{\mathcal{F}_n}$ it suffices to approximate the distribution of the Gaussian $P$-bridge process $G_P$.

Conceptually, there are two ways in which we can approach this problem: Either we take a nonparametric perspective, i.e. we estimate the $P$-measure via the empirical measure $P_n$ and then sample from the Gaussian $P_n$-bridge process. This leads to the well-known \emph{Gaussian multiplier bootstrap procedure}, a close cousin of the wild bootstrap. Or we take a parametric perspective, i.e. we estimate the covariance function of the Gaussian process $G_P$ and then sample from a Gaussian $Q_n$-motion constructed from the (truncated) Karhunen-Lo{\`e}ve expansion based on the estimated covariance function. This leads to the \emph{Gaussian process bootstrap procedure}, proposed by~\cite{giessing2023Gaussian}.

Let $\xi, \xi_1, \ldots, \xi_n \in \mathbb{R}$ be i.i.d. $N(0,1)$ random variables that are independent of $X, X_1, \ldots, X_n$. We define the \emph{Gaussian multiplier bootstrap process} as
\begin{align*}
	G_{P_n}(f) : = \frac{1}{\sqrt{n}}\sum_{i=1}^n \xi_i \big(f(X_i) - P_n f\big), \quad{} \text{where} \quad{} P_nf = \frac{1}{n} \sum_{i=1}^n f(X_i), \quad{} f \in \mathcal{F}_n.
\end{align*}
Notice that $G_{P_n}$ is really just the explicit representation of the Gaussian $P_n$-bridge process based on the empirical measure $P_n$. In particular, $G_{P_n}$ has mean zero and empirical covariance function
\begin{align*}
	(f, g ) \mapsto \mathrm{E}[G_{P_n}(f) G_{P_n}(g)] &= P_n fg - (P_nf)(P_ng) \\
	&= \frac{1}{n}\sum_{i=1}^n f(X_i)g(X_i) - \left(\frac{1}{n}\sum_{i=1}^n f(X_i)\right) \left(\frac{1}{n}\sum_{i=1}^n g(X_i)\right).
\end{align*}
Thus, provided that the empirical covariance function is uniformly close to the population covariance function $C_P$, it is reasonable to approximate $G_P$ via $G_{P_n}$. However, if the covariance function $C_P$ has a known structure, we might be better off by using an estimator that is more adaptive to this additional information. The Gaussian process bootstrap is an attempt in this direction.

In the following we consider the special case in which $\mathcal{F}_n$ can be identified with a compact subset of $\mathbb{R}^d$. (Essentially, we interpret the function class $\mathcal{F}_n$ as the dual of a compact subset of $\mathbb{R}^d$. Since $(\mathbb{R}^d)' = \mathbb{R}^d$, we use the same notation for $\mathcal{F}_n$ and the compact subset of $\mathbb{R}^d$.) This setup suffices for developing the Gaussian multiplier bootstrap in the context of $\ell_p$-norms and allows us to avoid questions about measurability of the Gaussian process. For a more general theory and a more careful treatment of measurability issues we refer to~\cite{giessing2023Gaussian}.

By the Moore-Aronszajn theorem the associated reproducing kernel Hilbert space (RKHS) $\mathcal{H}(C_P)$ of a Gaussian process with covariance function $C_P: \mathcal{F}_n \times \mathcal{F}_n \rightarrow \mathbb{R}$ is the completion of the linear span
\begin{align*}
	\mathcal{H} = \left\{ h: \mathcal{F}_n \rightarrow \mathbb{R} : h(\cdot) = \sum_{j=1}^m a_i C_P(s_j, \cdot ), \: a_1, \ldots, a_m \in \mathbb{R}, \: s_1, \ldots, s_m \in \mathcal{F}_n, \: m \geq 1 \right\}
\end{align*}
relative to the norm $\|\cdot\|_\mathcal{H}$ induced by the inner product
\begin{align*}
	\langle f, g \rangle_{\mathcal{H}} = \left\langle \sum_{j=1}^l a_j C_P(s_j, \cdot), \: \sum_{k=1}^m b_k C_P(t_k, \cdot) \right\rangle_{\mathcal{H}} := \sum_{j=1}^l \sum_{k=1}^m a_j b_k C_P(s_j, t_k).
\end{align*}
Let $(\phi_j)_{j\geq 1}$ be an orthonormal basis of the RKHS $\mathcal{H}(C_P)$ endowed with norm $\|\cdot\|_{\mathcal{H}}$. Then, the Karhunen-Lo{\`e}ve expansion of $G_P$ (in $L_2(\Omega, \mathcal{A}, \mathbb{P})$) is given by
\begin{align}\label{eq:subsec:ResultsEPT-1}
	G_p(f) = \sum_{j=1}^\infty \xi_j \phi_j(f), \quad{} f \in \mathcal{F}_n,
\end{align}
where  $(\xi_j)_{j\geq 1}$ is an orthonormal sequence of mean zero Gaussian random variables. Let $(\widehat{\phi}_j)_{j\geq 1}$ be estimates of $(\phi_j)_{j\geq 1}$ based on $X_1, \ldots, X_n$ and define the \emph{Gaussian process bootstrap} as the $Q_n$-motion with $L_2(\Omega, \mathcal{A}, \mathbb{P})$-representation
\begin{align}\label{eq:subsec:ResultsEPT-2}
	Z_{Q_n}(f) = \sum_{j=1}^\infty \xi_j \widehat{\phi}_j(f), \quad{} f \in \mathcal{F}_n.
\end{align}
While the representation in eq.~\eqref{eq:subsec:ResultsEPT-2} is neat, the crux is how to compute the estimates $(\widehat{\phi}_j)_{j \geq 1}$ and, in doing so, exploit structure and information about the function class $\mathcal{F}_n$ and the covariance function $C_P$. For a detailed discussion and examples we refer to~\cite{giessing2023Gaussian}. In the following we only present the case which is relevant for deriving Theorems~\ref{theorem:Bootstrap-Lp-Norm} and~\ref{theorem:Bootstrap-Lp-Norm-Quantiles}.

Since $\mathcal{F}_n$ is a compact subset of $\mathbb{R}^d$, we have, for all $f, g \in \mathcal{F}_n$,
\begin{align*}
	C_P(f,g) = P(fg) - (Pf)(Pg) = f'\left(\mathrm{E}[XX'] - \mathrm{E}[X]\mathrm{E}[X]'\right)g \equiv f' \Sigma g,
\end{align*}
where $\Sigma \in \mathbb{R}^{d \times d}$ is symmetric and positive semidefinite. Let $(\lambda_j)_{j \geq 1}$ and $(\psi_j)_{j \geq 1}$ be, respectively, the eigenvalues and normalized eigenfunctions of the operator $\mathcal{K} : L_2(\mathcal{F}) \rightarrow L_2(\mathcal{F})$ defined by $\mathcal{K} \psi(\cdot) = \int_\mathcal{F} C_P(s,\cdot) \psi(s) ds$. Then, by Mercer's theorem, the Karhunen-Lo{\`e}ve expansion of $G_P$ in eq.~\eqref{eq:subsec:ResultsEPT-1} can be written as
\begin{align}\label{eq:subsec:ResultsEPT-3}
	G_p(f) = \sum_{j=1}^\infty \sqrt{\lambda_j} \xi_j \psi_j(f) 
	\equiv (\Sigma^{1/2}Z)'f, \quad{} f \in \mathcal{F}_n,
\end{align}
where $Z \sim N(0, I_d)$. Therefore, the Gaussian process bootstrap procedure in eq.~\eqref{eq:subsec:ResultsEPT-2} reduces to sampling from the $Q_n$-motion
\begin{align}\label{eq:subsec:ResultsEPT-4}
	Z_{Q_n}(f) = (\widehat{\Sigma}_n^{1/2}Z)'f, \quad{} f \in \mathcal{F}_n,
\end{align}
where $\widehat{\Sigma}_n$ is any positive semidefinite estimate of $\Sigma$, possibly constructed by exploiting known low-rank structure, bandedness, sparsity, or independence relations.
 
To unify the notation, let $\{Y_n(f) : f \in \mathcal{F}_n\}$ be a Gaussian process with mean zero, covariance function depending (in some way) on the empirical measure $P_n$, and bounded and 1-Lipschitz continuous sample paths w.r.t. the intrinsic standard deviation metric $\rho_{Y_n}(f, g) = (\mathrm{E}[(Y_n(f) - Y_n(g))^2])^{1/2}$. This covers the Gaussian multiplier bootstrap process $G_{P_n}$ as well as the Gaussian process bootstrap $Z_{Q_n}$. Therefore, the following result can be used to obtain non-asymptotic bounds on the Kolmogorov distance between the laws of $\|\mathbb{G}_n\|_{\mathcal{F}_n}$ and $\|G_{P_n}\|_{\mathcal{F}_n}$ and the laws of $\|\mathbb{G}_n\|_{\mathcal{F}_n}$ and $\|Z_{Q_n}\|_{\mathcal{F}_n}$, respectively.

\begin{lemma}[\citeauthor{giessing2023Gaussian},~\citeyear{giessing2023Gaussian}]\label{lemma:Bootstrap-DimFree}
	Let $\mathcal{F}_n \subset L_2(S, \mathcal{S}, P)$ be totally bounded w.r.t. a metric $\rho$ and have envelope $F_n \in L_2(S, \mathcal{S}, P)$. Let $\mathcal{F}_{n, \delta}' = \{ f - g: f, g \in \mathcal{F}_n, \: \rho(f, g) \vee \rho_{Y_n}(f,g) < \delta\left(\|F_n\|_{P,2} \wedge \mathrm{E}\|Y_n\|_{\mathcal{F}_n}\right)\}$ and suppose that there exists a function  $\psi_n$ such that
	\begin{align*}
		\mathrm{E} \|G_P\|_{\mathcal{F}_{n,\delta}'} \:  \vee \: \mathrm{E}\big[\|Y_n\|_{\mathcal{F}_{n,\delta}'} \mid X_1, \ldots, X_n\big] \: \vee \: \big\|\|\mathbb{G}_n\|_{\mathcal{F}_{n,\delta}'} \big\|_{P,1} \lesssim \psi_n(\delta) \sqrt{\mathrm{Var}(\|G_P\|_{\mathcal{F}_n})} .
	\end{align*}
	Let $r_n^2 = \inf\big\{  \psi_n(\delta): \delta > 0\big\}$. Then, for each $n \geq 1$ and $M_n \geq 0$,
	\begin{align*}
		&\sup_{s \geq 0} \Big|\mathbb{P}\left\{\|\mathbb{G}_n\|_{\mathcal{F}_n} \leq s \right\} -\mathbb{P}\left\{\|Y_n\|_{\mathcal{F}_n} \leq s \mid X_1, \ldots, X_n \right\} \Big|\\
		&\quad{}\quad{}\quad{} \lesssim \frac{\|F_n\|_{P, 3}}{n^{1/6}\sqrt{\mathrm{Var}(\|G_P\|_{\mathcal{F}_n})}} + \frac{\|F_n\mathbf{1}\{F_n > M_n\}\|_{P,3}^3}{\|F_n\|_{P,3}^3}\\
		&\quad{}\quad{}\quad{} \quad{} \quad{} + \frac{ \mathrm{E} \|G_P\|_{\mathcal{F}_n} + M_n}{\sqrt{n\mathrm{Var}(\|G_P\|_{\mathcal{F}_n})}} + \left( \frac{\Delta_{\mathcal{F}_n}(P, P_n)}{ \mathrm{Var}(\|G_P\|_{\mathcal{F}_n} )} \right)^{1/3} + r_n,
	\end{align*}
	where $\Delta_{\mathcal{F}_n}(P, P_n) := \sup_{f, g \in \mathcal{F}_n} \big| \mathrm{E}[G_P(f)G_P(g)] - \mathrm{E}[Y_n(f)Y_n(g)] \big|$, and $\lesssim$ hides an absolute constant independent of $n, r_n, \mathcal{F}_n, F_n, M_n, \psi_n, P_n$, and $P$.
\end{lemma}
Next, for $\alpha \in (0,1)$ arbitrary, denote the $\alpha$-quantile of the conditional law of $\|Y_n\|_{\mathcal{F}_n}$ given $X_1, \ldots, X_n$ by
\begin{align*}
	c_{n, P_n}(\alpha) &:= \inf \left\{s \geq 0: \mathbb{P}\left\{  \|Y_n\|_{\mathcal{F}_n} \leq s \mid X_1, \ldots, X_n\right\} \geq \alpha \right\}.
\end{align*}
Combing Lemma~\ref{lemma:Bootstrap-DimFree} with standard arguments~\citep[e.g.][]{chernozhukov2013GaussianApproxVec} we conclude that $c_{n, P_n}(\alpha)$ is a valid approximation of the $\alpha$-quantile of the sampling distribution of $\|\mathbb{G}_n\|_{\mathcal{F}_n}$:

\begin{lemma}[\citeauthor{giessing2023Gaussian},~\citeyear{giessing2023Gaussian}]\label{lemma:Bootstrap-DimFree-AsympSize}
	Consider the setup of Lemma~\ref{lemma:Bootstrap-DimFree}. Let $(\Theta_n)_{n \geq 1} \in \mathbb{R}$ be a sequence of arbitrary random variables, not necessarily independent of $X_1, \ldots, X_n$. Then, for $M_n \geq 0$, $n \geq 1$,
	\begin{align*}
		&\sup_{\alpha \in (0,1)} \Big|\mathbb{P}\left\{ \|\mathbb{G}_n\|_{\mathcal{F}_n} + \Theta_n \leq c_{n, P_n}(\alpha) \right\}  - \alpha \Big| \\
		&\quad{}\quad{}\quad{} \lesssim \frac{\|F_n\|_{P, 3}}{n^{1/6}\sqrt{\mathrm{Var}(\|G_P\|_{\mathcal{F}_n})}} + \frac{\|F_n\mathbf{1}\{F_n > M_n\}\|_{P,3}^3}{\|F_n\|_{P,3}^3}+ \frac{ \mathrm{E} \|G_P\|_{\mathcal{F}_n} + M_n}{\sqrt{n\mathrm{Var}(\|G_P\|_{\mathcal{F}_n})}} + r_n\\
		& \quad{}\quad{}\quad{}\quad{}\quad{} + \inf_{\delta > 0}\left\{ \left(\frac{\delta }{\mathrm{Var}(\|G_P\|_{\mathcal{F}_n})} \right)^{1/3}  + \mathrm{P}\left(\Delta_{\mathcal{F}_n}(P, P_n) > \delta\right)\right\}\\
		&\quad{}\quad{}\quad{} \quad{} \quad{} \quad{} + \inf_{\eta > 0} \left\{ \frac{\eta }{ \sqrt{\mathrm{Var}(\|G_P\|_{\mathcal{F}_n})} } + \mathrm{P}\left(|\Theta_n| > \eta \right)\right\},
	\end{align*}
	where $\lesssim$ hides an absolute constant independent of $n, r_n, \mathcal{F}_n, F_n, M_n, \psi_n, P_n, P$, and the law of $(\Theta_n)_{n\geq 1}$.
\end{lemma}

We conclude this section with several technical results involving empirical process (notation) that we use throughout the proofs of the results in the main text.

\begin{lemma}[\citeauthor{giessing2023Gaussian},~\citeyear{giessing2023Gaussian}]\label{lemma:Bootstrap-Sup-Empirical-Process-Quantil-Comparison}
	For $\alpha \in (0,1)$ arbitrary let
	 \begin{align*}
		c_{n, P}(\alpha) &:= \inf \left\{s \geq 0: \mathbb{P}\left\{  \|G_P\|_{\mathcal{F}_n} \leq s \right\} \geq \alpha \right\}, \quad{} \text{and}\\
		c_{n, P_n}(\alpha) &:= \inf \left\{s \geq 0: \mathbb{P}\left\{  \|Y_n\|_{\mathcal{F}_n} \leq s \mid X_1, \ldots, X_n\right\} \geq \alpha \right\}.
	\end{align*}
	Then, for all $\delta > 0$,
	\begin{align*}
		&\inf_{\alpha \in (0,1)} \mathbb{P}\left\{ c_{n, P_n}(\alpha) \leq c_{n, P}(\pi_{n,P}(\delta) + \alpha) \right\} \geq 1 - \mathrm{P}\left( \Delta_{\mathcal{F}_n}(P, P_n) > \delta \right), \quad{} \text{and}\\
		&\inf_{\alpha \in (0,1)} \mathbb{P}\left\{ c_{n, P}(\alpha) \leq c_{n, P_n}(\pi_{n,P}(\delta) + \alpha) \right\} \geq 1 - \mathrm{P}\left( \Delta_{\mathcal{F}_n}(P, P_n) > \delta \right),
	\end{align*}
	where $\Delta_{\mathcal{F}_n}(P, P_n) := \sup_{f, g \in \mathcal{F}_n} \big| \mathrm{E}[G_P(f)G_P(g)] - \mathrm{E}[Y_n(f)Y_n(g)] \big|$, $\pi_{n,P}(\delta) := K\delta^{1/3} \left(\mathrm{Var}(\|G_P\|_{\mathcal{F}_n})\right)^{-1/3}$ and $K > 0$ is an absolute constant.
\end{lemma}

\begin{lemma}[\citeauthor{giessing2023Gaussian},~\citeyear{giessing2023Gaussian}]\label{lemma:Bootstrap-Sup-Empirical-Process-Quantil-Comparison-Non-Gaussian}
	Let $\{Z_n(f) = Z_n^1(f) + Z_n^2(f) : f \in \mathcal{F}_n\}$ be an arbitrary stochastic process. For $\alpha \in (0,1)$ arbitrary define $c_{n, Z_n}(\alpha) := \inf \{s \geq 0: \mathbb{P}\{\|Z_n\|_{\mathcal{F}_n} \leq s \mid X_1, \ldots, X_n\} \geq \alpha \}$. Then, for all $\delta > 0$,
	\begin{align*}
		&\inf_{\alpha \in (0,1)} \mathbb{P}\left\{ c_{n, Z_n}(\alpha) \leq c_{n, P}(\kappa_{n,P}(\delta) +  \delta + \alpha) \right\} \geq 1 - \mathbb{P} \left\{ \gamma_{n, Z_n^1} + \rho_{n,Z_n^2}(\delta) > \delta \right\},  \quad{} \text{and}\\
		&\inf_{\alpha \in (0,1)} \mathbb{P}\left\{ c_{n, P}(\alpha) \leq c_{n, Z_n}(\kappa_{n,P}(\delta) +  \delta + \alpha) \right\} \geq 1 - \mathbb{P} \left\{ \gamma_{n, Z_n^1} + \rho_{n,Z_n^2}(\delta) > \delta \right\},
	\end{align*}
	where $\gamma_{n,Z_n^1} := \sup_{s \geq 0} | \mathbb{P}\{\|Z_n^1\|_{\mathcal{F}_n} \leq s \mid X_1, \ldots, X_n\} - \mathbb{P}\{\|G_P\|_{\mathcal{F}_n} \leq s\}  |$, $\kappa_{n,P}(\delta) :=  K \delta /\sqrt{\mathrm{Var}(\|G_P\|_{\mathcal{F}_n})}$, $\rho_{n,Z_n^2}(\delta) :=  \mathbb{P}\left\{ \|Z_n^2\|_{\mathcal{F}_n} > \delta \mid X_1, \ldots, X_n \right\}$, and $K > 0$ is an absolute constant.
\end{lemma}

\begin{lemma}[\citeauthor{giessing2022anticoncentration},~\citeyear{giessing2022anticoncentration}]\label{lemma:AntiConcentration-SeparableProcess}
	Let $X = (X_u)_{u \in U}$ be a centered separable Gaussian process indexed by a semi-metric space $U$. Set $Z = \sup_{u \in U}X_u$ and assume that $0 \leq Z < \infty$ a.s. For all $\varepsilon \geq 0$,
	\begin{align*}
		\frac{\varepsilon/\sqrt{12}}{ \sqrt{\mathrm{Var}(Z) + \varepsilon^2/ 12}} \leq \sup_{t\geq 0} \mathrm{P}\left( t \leq  Z \leq t + \varepsilon  \right) \leq \frac{\varepsilon\sqrt{12}}{ \sqrt{\mathrm{Var}(Z) + \varepsilon^2/ 12}}.
	\end{align*}
	The result remains true if $Z$ is replaced by $\widetilde{Z} = \sup_{u \in U}|X_u|$.
\end{lemma}

\begin{lemma}[\citeauthor{giessing2022anticoncentration},~\citeyear{giessing2022anticoncentration}]\label{lemma:BoundsVariance-SeparableProcess}
	Let $X = (X_u)_{u \in U}$ be a separable Gaussian process indexed by a semi-metric space $U$ such that $\mathrm{E}[X_u] = 0$, $0 < \underline{\sigma}^2 \leq \mathrm{E}[X_u^2] \leq \bar{\sigma}^2 < \infty$, and $|\mathrm{Corr}(X_u, X_v)| \leq \rho$ for all $u \neq v \in U$. Set $Z = \sup_{u \in U}X_u$ and assume that $Z < \infty$ a.s. Then, $ 0 \leq \mathrm{E}[Z] < \infty$ and there exist absolute constants $c, C > 0$ such that
	\begin{align*}
		\frac{1}{C} \left(\frac{\underline{\sigma}}{1 + \mathrm{E}[Z/\underline{\sigma}]}\right)^2 \leq \mathrm{Var}(Z) \leq C \left[ \bar{\sigma}^2 \wedge \left( \bar{\sigma}^2 \rho +  \left(\frac{\bar{\sigma}}{ (\mathrm{E}[Z/ \bar{\sigma}] - c)_+}\right)^2 \right) \right],
	\end{align*}
	with the convention that ``$1/0 = \infty$''. The result remains true if $Z$ is replaced by $\widetilde{Z} = \sup_{u \in U}|X_u|$.
\end{lemma}

\begin{lemma}[\citeauthor{ding2015multiple},~\citeyear{ding2015multiple}]\label{lemma:ding2015multiple}
	Let $X = (X_u)_{u \in U}$ be a centered Gaussian process indexed by a semi-metric space $U$ such that $\sup_{u \in U} \mathrm{Var}(X_u) < \infty$ and $Z: = \sup_{u \in U} X_u < \infty$ a.s. For $t \in (0,1)$ set $U_t = \{u \in U: X_u \geq t \mathrm{E}[Z]\}$. Then, for $\lambda > 0$ arbitrary,
	\begin{align*}
		\mathrm{P} \left( \sup_{u \in U_t} X_u' \geq \sqrt{1- t^2} \mathrm{E}[Z] + \frac{\lambda}{\sqrt{1-t^2}} \right) \leq \frac{\mathrm{Var}(Z)}{\lambda^2},
	\end{align*}
	where $X' = (X_u')_{u \in U}$ is an independent copy of $X$.
\end{lemma}

\begin{lemma}\label{lemma:ModulusContinuityRn}
	Let $\mathcal{F}_\delta' = \{f - g : f, g \in \mathcal{F},\: \rho_2(f, g) \leq \delta\}$, where $\rho_2^2(f, g) = P(f-g)^2$ and $\mathcal{F} = \{ x \mapsto u'x :  \|u\|_q = 1, \: u \in \mathbb{R}^d\}$, $q \geq 1$. There exist functions $\upsilon_n, \psi_n, \phi_n, $ such that $\mathrm{E} \big[\|Z_{Q_n}\|_{\mathcal{F}_{\delta}'} \mid X_1, \ldots, X_n \big] \lesssim 	\upsilon_n(\delta) \sqrt{\mathrm{Var}(\|G_P\|_{\mathcal{F}})}$, $\mathrm{E} \|G_P\|_{\mathcal{F}_{\delta}'} \lesssim 	\psi_n(\delta) \sqrt{\mathrm{Var}(\|G_P\|_{\mathcal{F}})}$, $\big\|\|\mathbb{G}_n\|_{\mathcal{F}_{\delta}'} \big\|_{P,1} \lesssim 	\phi_n(\delta)\sqrt{\mathrm{Var}(\|G_P\|_{\mathcal{F}})}$, and
	$\inf\big\{ \upsilon_n(\delta) \vee  \psi_n(\delta) \vee \phi_n(\delta): \delta > 0\big\} = 0$.
\end{lemma}

\subsection{Bounds on moments and quantiles and naive large deviation principles}\label{subsec:Miscellaneous}
In this section we collect miscellaneous technical lemmas used throughout the proofs of the main results of this paper.

\begin{lemma}[Lower and upper bounds on moments of $\ell_p$-norms of Gaussian random vectors]\label{lemma:EVNormGaussian}
	Let $Z\in \mathbb{R}^d$ be a centered Gaussian random vector with marginals $Z_k \sim N(0, \sigma_k^2)$. Then, for $1 \leq p \leq \infty$,
	\begin{align*}
		\left( \sum_{k=1}^d \sigma_k^p\right)^{1/p} \bigvee \sqrt{\log d}\left(\min_{1 \leq k \leq d} \sigma_k^2\right)^{1/2}\lesssim \mathrm{E}[\|Z\|_p] \lesssim \sqrt{p}\left( \sum_{k=1}^d \sigma_k^p\right)^{1/p} \bigwedge  d^{1/p}\sqrt{\log d}\left(\max_{1 \leq k \leq d} \sigma_k^2\right)^{1/2},
	\end{align*}
	where $\lesssim$ hides an absolute constant independent of $d, p$, and the covariance matrix of $Z$.
\end{lemma}

\begin{lemma}[Lower and upper bounds on the quantiles of $\ell_p$-norms of Gaussian random vectors]\label{lemma:UpperBoundQuantiles}
	Let $Z \sim N(0, \Sigma)$. For all $\alpha \in (0,1/2]$,
	\begin{align*}
		&- \left(\|\Sigma^{1/2}\|_{2 \rightarrow p} \wedge \sqrt{\mathrm{Var}(\|Z\|_p)} \right)\\
		&\quad{}\quad{}\quad{}\quad{}\quad{}\quad{}\quad{}\leq c_{n,p}(1-\alpha) - \mathrm{E}[\|Z\|_p] \\
		& \quad{}\quad{}\quad{} \quad{}\quad{}\quad{}\quad{}\quad{}\quad{}\quad{}\quad{}\quad{}\quad{}\leq  \sqrt{2\log(1/\alpha)} \|\Sigma^{1/2}\|_{2 \rightarrow p} \wedge \sqrt{(1/\alpha) \mathrm{Var}(\|Z\|_p)}.
	\end{align*}
	In fact, the upper bound holds for all $\alpha \in (0, 1)$.
\end{lemma}

\begin{lemma}[Reverse Lyapunov inequality for norms of log-concave random vectors]\label{lemma:ReverseLyapunovLogConcave}
	Let $X \in \mathbb{R}^d$ be a log-concave random vector and $\|\cdot\|$ any norm on $\mathbb{R}^d$. For any $1 \leq s \leq t < \infty$, 
	\begin{align*}
		\frac{1}{t}\left(\mathrm{E}[\|X\|^t]\right)^{1/t} \leq \frac{6}{s} \left(\mathrm{E}[\|X\|^s]\right)^{1/s}.
	\end{align*}
\end{lemma}
\begin{remark}
	The constant 6 is not optimal, as can be seen by the smaller constant of the reverse Lyapunov inequality for Gaussian measures~\citep[e.g.][Proposition  A.2.4]{vandervaart1996weak}. The single most important consequence of this result is that norms of log-concave random vectors have sub-exponential norms; in particular, $\left\| \|X\| \right\|_{\psi_1} \leq 6 \mathrm{E}[\|X\|]$.
\end{remark}

\begin{lemma}[LDP for Gaussian random vectors in changing dimensions]\label{lemma:LDPGaussian}
	Let $Z \in \mathbb{R}^d$ be a centered Gaussian random vector with positive semi-definite covariance matrix $\Sigma \neq \mathbf{0} \in \mathbb{R}^{d \times d}$. Let $1 \leq p \leq \infty$ be arbitrary and $(t_n)_{n \geq 1}$ be such that $t_n \rightarrow  + \infty$ and $\limsup_{n \rightarrow \infty} \mathrm{E}[\|Z\|_p]/(t_n \|\Sigma^{1/2}\|_{2 \rightarrow p}) < 1$. Then,
	\begin{align*}
		- \log  \mathrm{P} \left(\|Z\|_p > t_n \|\Sigma^{1/2}\|_{2\rightarrow p} \right) \asymp t_n^2 \quad{}\quad{} \mathrm{as} \quad{} \quad{} n \rightarrow \infty.
	\end{align*}
\end{lemma}
\begin{remark}
	The key point is that the dimension of the vector $Z \in \mathbb{R}^d$ is allowed to depend in an arbitrary fashion on $n$ provided that $t_n$ diverges and $\limsup_{n \rightarrow \infty} \mathrm{E}[\|Z\|_p]/(t_n \|\Sigma^{1/2}\|_{2 \rightarrow p}) < 1$. In particular, $d$ is allowed to grow (arbitrarily fast) as $n \rightarrow \infty$. It is possible that this fact is already known to some specialists; however, we could not locate it in the literature.
\end{remark}
\begin{remark}
	Since the Gaussian isoperimetric inequality and an elementary lower bound on the density of the Gaussian tail probability are key ingredients of the proof, any extension to non-Gaussian random vectors needs to be carried out very differently (e.g. along the lines of the Ellis-G{\"a}rtner theorem).
\end{remark}

\begin{lemma}[Lower LDP for strictly log-concave random vectors in changing dimensions]\label{lemma:LDPStrictly-Log-Concave}
	Let $X \in \mathbb{R}^d$ be a centered random vector with positive semi-definite covariance matrix $\Sigma \neq \mathbf{0} \in \mathbb{R}^{d \times d}$ and strictly log-concave density $f = e^{-\varphi}$ where $\varphi'' \geq \lambda_d I$ with $\lambda_d > 0$. Let $1 \leq p \leq \infty$ be arbitrary and $(t_n)_{n \geq 1}$ be a sequence of numbers in $\mathbb{R}_+$ such that $\limsup_{n \rightarrow \infty}(\sqrt{\lambda_d} \wedge 1)\mathrm{E}[\|X\|_p]/(t_n \|\Sigma^{1/2}\|_{2 \rightarrow p}) < 1$. Then,
	\begin{align*}
		- \log  \mathrm{P} \left(\|X\|_p > t_n \|\Sigma^{1/2}\|_{2\rightarrow p}/\sqrt{\lambda_d \wedge 1} \right) \gtrsim t_n^2 \quad{}\quad{} \mathrm{as} \quad{} \quad{} n \rightarrow \infty.
	\end{align*}
\end{lemma}

\begin{lemma}[Lower LDP for log-concave random vectors in changing dimensions]\label{lemma:LDP-Log-Concave}
	Let $X \in \mathbb{R}^d$ be a centered random vector with positive semi-definite covariance matrix $\Sigma \neq \mathbf{0} \in \mathbb{R}^{d \times d}$ and log-concave density. Let $1 \leq p \leq \infty$ be arbitrary and $(t_n)_{n \geq 1}$ be a sequence of numbers in $\mathbb{R}_+$ such that $\limsup_{n \rightarrow \infty}\mathrm{E}[\|X\|_p]/(t_n \|\Sigma^{1/2}\|_{2 \rightarrow p}) < 1$. Then,
	\begin{align*}
		- \log  \mathrm{P} \left(\|X\|_p > t_n \|\Sigma^{1/2}\|_{2\rightarrow p} \right) \gtrsim t_n \quad{}\quad{} \mathrm{as} \quad{} \quad{} n \rightarrow \infty.
	\end{align*}
\end{lemma}

\begin{lemma}[Upper bound on the variance of $\ell_p$-norms of uniform random vectors]\label{lemma:BoundsVarianceUniformRV} Let $U \sim \mathrm{Unif}(S^{n-1})$ and $\Gamma \in \mathbb{R}^{d \times n}$. For all $p \in [1, \infty]$,
	\begin{align*}
		\mathrm{Var}(\|\Gamma U\|_p) \lesssim \frac{\|\Gamma\|_{2 \rightarrow p}^2}{n^2},
	\end{align*}
	where $\lesssim$ hides an absolute constant independent of $p, n, d$, and $\Gamma$.
\end{lemma}
\begin{remark}
	This lemma establishes the superconcentration of $\ell_p$-norms of uniform random vectors.
\end{remark}
\newpage
\section{Proofs of results in the main text}\label{sec:Proofs}

\subsection{Proofs of Section~\ref{subsec:Assumptions}}\label{subsec:proofs-Assumptions}

\begin{proof}[\textbf{Proof of Lemma~\ref{lemma:example:SubGaussianData}}]
	Observe that
	\begin{align*}
		\widehat{\Omega}_n - \Omega = \frac{1}{n}\sum_{i=1}^n R(X_i - \mu_n)(X_i - \mu_n)'R' - \Omega - \left(\frac{1}{n}\sum_{i=1}^n  R(X_i - \mu_n)\right)\left(\frac{1}{n}\sum_{i=1}^n R(X_i - \mu_n)\right)' \equiv \mathbf{I} + \mathbf{II}.
	\end{align*}
	Consider case (i). Recall that for any $V \in \mathbb{R}^d$ sub-Gaussian with mean zero and covariance matrix $\Omega$, and all $s \geq 2$, 
	\begin{align*}
		\sqrt{\mathrm{tr}(\Omega)} = \left(\mathrm{E}[\|V\|_2^2]\right)^{1/2} \overset{(a)}{\lesssim}\left(\mathrm{E}[\|V\|_2^s]\right)^{1/s} \overset{(b)}{\lesssim} \sqrt{\frac{s}{2}} \left(\mathrm{E}[\|V\|_2^2]\right)^{1/2} =   \sqrt{\frac{s}{2}} \sqrt{\mathrm{tr}(\Omega)},
	\end{align*}
	where $\lesssim$ hides an absolute constant independent of $d,t, \Omega$ and (a) follows from H{\"o}lder's inequality and (b) follows from Fubini's theorem and integrating over the exponential tail probability (see also proof of Lemma~\ref{lemma:ReverseLyapunovLogConcave}). Thus, for all $s \geq 3$,
	\begin{align*}
		\sqrt{\|\Omega\|_{op}}\sqrt{\mathrm{r}(\Omega)} \lesssim \left(\mathrm{E}[\|R(X- \mu_n)\|_2^s]\right)^{1/s} \vee \left(\mathrm{E}[\|Z\|_2^s]\right)^{1/s} \lesssim \sqrt{s\|\Omega\|_{op}}\sqrt{\mathrm{r}(\Omega)}.
	\end{align*}
	Moreover, since the linear function $u \mapsto X'u$ from $\mathbb{R}^d$ to $\mathbb{R}$ is pre-Gaussian, Theorem 4 in~\cite{koltchinskiiConcentration2017} applies and we conclude via Markov's inequality that
	\begin{align*}
		\|\mathbf{I}\|_{2 \rightarrow 2} = \|\mathbf{I}\|_{op}  = O_p\left( \|\Omega\|_{op} \left( \sqrt{\frac{r(\Omega)}{n}} \vee \frac{r(\Omega)}{n} \right) \right).
	\end{align*}
	Also, by Hoeffding's inequality, the union bound, and Markov's inequality,
	\begin{align*}
		\|\mathbf{II}\|_{2 \rightarrow 2} = \|\mathbf{II}\|_{op}  =  \left\|\frac{1}{n}\sum_{i=1}^n R(X_i - \mu_n) \right\|_2^2 = O_p \left( \|\Omega\|_{op} \frac{r(\Omega)}{n}\right),
	\end{align*}
	because, by straightforward adaptation of the arguments of the proof of Theorem 4 in~\cite{koltchinskiiConcentration2017} (i.e. Talagrand's majorizing measure bounds for Gaussian processes, e.g. Theorem 2.5 in~\cite{shahar2010empirical}),
	\begin{align*}
		\mathrm{E}\left[\left\|\frac{1}{n}\sum_{i=1}^n R(X_i - \mu_n) \right\|_2 \right] \lesssim \sqrt{\frac{\mathrm{tr}(\Omega)}{n}}.
	\end{align*}
	Lastly, by Theorem~\ref{theorem:LowerBoundVariance-LpNorm-Refinements},
	\begin{align*}
		\mathrm{Var}(\|Z\|_2) \geq \frac{\mathrm{tr}(\Omega^2)}{\mathrm{tr}(\Omega)} = \|\Omega\|_{op}  \frac{\mathrm{r}(\Omega^2)}{\mathrm{r}(\Omega)}.
	\end{align*}
	Combine these inequalities to verify that Assumptions~\ref{assumption:ControlThirdMoments}--\ref{assumption:ConsistentCovariance} hold under the stated conditions.
	
	Consider case (ii). Since $X \in \mathbb{R}^d$ is sub-Gaussian with mean $\mu_n$, for all $s \geq 1$,
	\begin{align*}
		\omega_{(t)} \lesssim \left(\mathrm{E}[\|R(X- \mu_n)\|_\infty^s]\right)^{1/s}  \vee \left(\mathrm{E}[\|Z\|_\infty^s]\right)^{1/s} \lesssim s  \sqrt{ \omega_{(t)}^2 \log t},
	\end{align*}
	where $\lesssim$ hides an absolute constant independent of $d,t,\Omega$. Again, by Theorem 4 in~\cite{koltchinskiiConcentration2017} and Markov's inequality,
	\begin{align*}
		\|\mathbf{I}\|_{1 \rightarrow \infty} = \max_{1 \leq j,k \leq t} |\mathbf{I}_{jk}|  = O_p\left( \omega_{(t)}^2 \left( \sqrt{\frac{ \log t}{n}} \vee \frac{ \log t}{n} \right) \right),
	\end{align*}
	by Hoeffding's inequality, the union bound, and Markov's inequality,
	\begin{align*}
		\|\mathbf{II}\|_{1 \rightarrow \infty } = \left\|\frac{1}{n}\sum_{i=1}^n R(X_i - \mu_n) \right\|_\infty^2 = O_p \left(\omega_{(t)}^2 \frac{\log t}{n}\right),
	\end{align*}
	because
	\begin{align*}
			\mathrm{E}\left[\left\|\frac{1}{n}\sum_{i=1}^n R(X_i - \mu_n) \right\|_\infty \right] \lesssim \omega_{(t)}\sqrt{ \frac{\log t}{n}}.
	\end{align*}
	By Theorem~\ref{theorem:LowerBoundVariance-LpNorm-Refinements},
	\begin{align*}	
		\mathrm{Var}(\|Z\|_\infty) \gtrsim \left(\frac{\omega_{(1)}^2}{\omega_{(1)} + \omega_{(t)} \sqrt{\log t}}\right)^2.
	\end{align*}
	Combine these bounds with $s= 1/6$ to conclude that Assumptions~\ref{assumption:ControlThirdMoments}--\ref{assumption:ConsistentCovariance} hold under the stated conditions.
\end{proof}

\begin{proof}[\textbf{Proof of Lemma~\ref{lemma:example:LogConcaveData}}]
	Recall that
	\begin{align*}
		\widehat{\Omega}_n - \Omega &= \frac{1}{n}\sum_{i=1}^n R(X_i - \mu_n)(X_i - \mu_n)'R' - \Omega - \left(\frac{1}{n}\sum_{i=1}^n  R(X_i - \mu_n)\right)\left(\frac{1}{n}\sum_{i=1}^n R(X_i - \mu_n)\right)'\\
		& \equiv \breve{\Omega}_n - \Omega - \breve{X}_n \breve{X}_n'.
	\end{align*}	
	Consider case (i). For any $V \in \mathbb{R}^d$ with log-concave density with mean zero and covariance matrix $\Omega$, we have, for all $s \geq 2$, 
	\begin{align*}
		\sqrt{\mathrm{tr}(\Omega)} = \left(\mathrm{E}[\|V\|_2^2]\right)^{1/2} \overset{(a)}{\lesssim}\left(\mathrm{E}[\|V\|_2^s]\right)^{1/s} \overset{(b)}{\lesssim} \sqrt{\frac{s}{2}} \left(\mathrm{E}[\|V\|_2^2]\right)^{1/2} =   \sqrt{\frac{s}{2}} \sqrt{\mathrm{tr}(\Omega)},
	\end{align*}
	where $\lesssim$ hides an absolute constant independent of $d,t, \Omega$ and (a) follows from H{\"o}lder's inequality and (b) follows from Lemma~\ref{lemma:ReverseLyapunovLogConcave}. Thus, for all $s \geq 3$,
	\begin{align*}
		\sqrt{\|\Omega\|_{op}}\sqrt{\mathrm{r}(\Omega)} \lesssim \left(\mathrm{E}[\|R(X- \mu_n)\|_2^s]\right)^{1/s} \vee \left(\mathrm{E}[\|Z\|_2^s]\right)^{1/s} \lesssim \sqrt{s\|\Omega\|_{op}}\sqrt{\mathrm{r}(\Omega)}.
	\end{align*}
	Moreover, since the linear function $u \mapsto X'u$ from $\mathbb{R}^d$ to $\mathbb{R}$ is pre-Gaussian, Gaussian symmetrization, followed by Theorem 4 in~\cite{koltchinskiiConcentration2017} and an iteration step imply via Markov's inequality that
	\begin{align}\label{eq:lemma:example:LogConcaveData-1}
		\|\breve{\Omega}_n - \Omega \|_{2 \rightarrow 2} = \|\breve{\Omega}_n - \Omega \|_{op}   = O_p\left(\sqrt{ \|\Omega\|_{op}} \left(\left(\frac{r(\Omega)}{n}\right)^{1/4} \vee \sqrt{ \frac{r(\Omega)}{n}} \right) \right).
	\end{align}
	Indeed, let $g_1, \ldots, g_n$ be i.i.d. standard normal random variables and compute
	\begin{align*}
		\mathrm{E}[\|\breve{\Omega}_n - \Omega \|_{op}] & \lesssim \mathrm{E}\left[\sup_{\|u\|_2=1} \left| \frac{1}{n}\sum_{i=1}^n g_i \big(R(X_i-\mu_n)'u\big)^2 \right| \right]\\
		&\overset{(a)}{\lesssim} \mathrm{E}\left[\sqrt{\frac{\|\breve{\Omega}_n \|_{op}\mathrm{tr}(\breve{\Omega}_n)}{n}} \vee \frac{\mathrm{tr}(\breve{\Omega}_n )}{n} \right]\\
		&\lesssim  \left(\mathrm{E}\left[\|\breve{\Omega}_n\|_{op}\right]\right)^{1/2}\sqrt{\frac{\mathrm{tr}(\Omega)}{n}} \vee \frac{\mathrm{tr}(\Omega)}{n}\\
		&\lesssim  \left(\mathrm{E}\left[\|\breve{\Omega}_n - \Omega\|_{op}\right]\right)^{1/2}\sqrt{\frac{\mathrm{tr}(\Omega)}{n}} \vee \| \Omega\|_{op}^{1/2}\sqrt{\frac{\mathrm{tr}(\Omega)}{n}} \vee \frac{\mathrm{tr}(\Omega)}{n},
	\end{align*}
	where (a) holds by Theorem 4 in~\cite{koltchinskiiConcentration2017} and the other inequalities follow from Gaussian symmetrization, Cauchy-Schwarz, and the triangle inequality. Solving above inequality for $\mathrm{E}[\|\breve{\Omega}_n- \Omega\|_{op}]$ yields the bound in eq.~\eqref{eq:lemma:example:LogConcaveData-1} via Markov's inequality.
	
	Next, since the linear function $u \mapsto X'u$ from $\mathbb{R}^d$ to $\mathbb{R}$ is pre-Gaussian, Gaussian symmetrization, followed by Theorem 2.5 in~\cite{shahar2010empirical} and similar steps as in the proof of Theorem 4 in~\cite{koltchinskiiConcentration2017} yield,
	\begin{align}\label{eq:lemma:example:LogConcaveData-2}
		\|\breve{X}_n \breve{X}_n'\|_{2 \rightarrow 2} = \|\breve{X}_n \breve{X}_n'\|_{op}  =  \left\|\frac{1}{n}\sum_{i=1}^n R(X_i - \mu_n) \right\|_2^2 = O_p \left( \|\Omega\|_{op} \frac{r(\Omega)}{n}\right).
	\end{align}
		Indeed, let $g_1, \ldots, g_n$ be i.i.d. standard normal random variables and compute
	\begin{align*}
		\mathrm{E}\left[ \left\|\frac{1}{n}\sum_{i=1}^n R(X_i - \mu_n) \right\|_2 \right] \lesssim \mathrm{E}\left[\sup_{\|u\|_2=1} \left| \frac{1}{n}\sum_{i=1}^n g_i R(X_i-\mu_n)'u \right| \right] \lesssim \mathrm{E}\left[\sqrt{\frac{\mathrm{tr}(\breve{\Omega}_n)}{n}}\: \right] \lesssim \sqrt{\frac{\mathrm{tr}(\Omega)}{n}}.
	\end{align*}
	Moreover, by Theorem~\ref{theorem:LowerBoundVariance-LpNorm-Refinements},
	\begin{align*}
		\mathrm{Var}(\|Z\|_2) \geq \frac{\mathrm{tr}(\Omega^2)}{\mathrm{tr}(\Omega)} = \|\Omega\|_{op}  \frac{\mathrm{r}(\Omega^2)}{\mathrm{r}(\Omega)}.
	\end{align*}
	Combine these inequalities to verify that Assumptions~\ref{assumption:ControlThirdMoments}--\ref{assumption:ConsistentCovariance} hold under the stated conditions.
	
	Consider case (ii). Since $X \in \mathbb{R}^d$ has log-concave density, by Lemma~\ref{lemma:ReverseLyapunovLogConcave}, for all $s \geq 2$,
	\begin{align*}
		\omega_{(t)} \lesssim \left(\mathrm{E}[\|R(X- \mu_n)\|_\infty^s]\right)^{1/s}  \vee \left(\mathrm{E}[\|Z\|_\infty^s]\right)^{1/s} \lesssim s \left(\mathrm{E}[\|R(X- \mu_n)\|_\infty] \vee \mathrm{E}[\|Z\|_\infty]\right)
	\end{align*}
	where $\lesssim$ hides an absolute constant independent of $d,t,\Omega,\varphi$. Next, since $X \in \mathbb{R}^d$ has log-concave distribution so has $(R(X- \mu_n))_k$ for $1 \leq k \leq d$. Hence, by Lemma 1 in~\cite{cule2010theoretical} there exist constants $a > 0, b \in \mathbb{R}$ (depending on the marginals of density $f= e^{-\varphi}$) such that for any $1/K < a$, 
	\begin{align*}
		\mathrm{E}[\|R(X- \mu_n)\|_\infty] \leq \omega_{(t)} K  \mathrm{E}[\|\widetilde{X}\|_\infty] &\leq \omega_{(t)} K \left( (\log t) + \max_{1 \leq k \leq t} \mathrm{E}\big[ e^{|\widetilde{X}_k|} \big]  \right)\\
		& \leq   \omega_{(t)} K \left( (\log t) + \frac{e^b}{a} \frac{K}{K-1}\right)\\
		& \equiv \omega_{(t)} \log (t K_\varphi),
	\end{align*}
	where $\widetilde{X} = \left( \big(R(X - \mu_n) \big)_k/(K\omega_k) \right)_{k=1}^d$ and $K_\varphi > 1$ is a constant depending on (the marginals of) the density $f= e^{-\varphi}$. Also, as usual,
	\begin{align*}
		\mathrm{E}[\|Z\|_\infty] \lesssim \sqrt{ \omega_{(t)}^2 \log t},
	\end{align*}
	and, repeating the iterative argument that led to eq.~\eqref{eq:lemma:example:LogConcaveData-1}, we obtain
	\begin{align*}
		\|\breve{\Omega}_n- \Omega\|_{1, \infty} = \max_{1 \leq j,k \leq t} |\breve{\omega}_{jk}- \omega_{jk}|  = O_p\left( \omega_{(t)} \left( \left( \frac{ \log t}{n} \right)^{1/4} \vee \sqrt{\frac{ \log t}{n}} \right) \right),
	\end{align*}
	and, by the arguments that gave eq.~\eqref{eq:lemma:example:LogConcaveData-2},
	\begin{align*}
		\|\breve{X}_n \breve{X}_n'\|_{1 \rightarrow \infty } = \left\|\frac{1}{n}\sum_{i=1}^n R(X_i - \mu_n) \right\|_\infty^2 = O_p \left(\omega_{(t)}^2 \frac{\log t}{n}\right),
	\end{align*}
	and, by Theorem~\ref{theorem:LowerBoundVariance-LpNorm-Refinements},
	\begin{align*}	
		\mathrm{Var}(\|Z\|_\infty) \gtrsim \left(\frac{\omega_{(1)}^2}{\omega_{(1)} + \omega_{(t)} \sqrt{\log t}}\right)^2.
	\end{align*}
	We combine these bounds with $s= 1/6$ and conclude that Assumptions~\ref{assumption:ControlThirdMoments}--\ref{assumption:ConsistentCovariance} hold under the stated conditions.
\end{proof}

\begin{proof}[\textbf{Proof of Lemma~\ref{lemma:example:HeavyTailedData}}]
	We have
	\begin{align*}
		\widehat{\Omega}_n - \Omega &= \frac{1}{n}\sum_{i=1}^n R(X_i - \mu_n)(X_i - \mu_n)'R' - \Omega - \left(\frac{1}{n}\sum_{i=1}^n  R(X_i - \mu_n)\right)\left(\frac{1}{n}\sum_{i=1}^n R(X_i - \mu_n)\right)'\\
		& \equiv \breve{\Omega}_n - \Omega - \breve{X}_n \breve{X}_n'.
	\end{align*}	
	Consider case (i). By Theorem~\ref{theorem:LowerBoundVariance-LpNorm-Refinements},
	\begin{align*}
		\mathrm{Var}(\|Z\|_2) \geq \frac{\mathrm{tr}(\Omega^2)}{\mathrm{tr}(\Omega)} = \|\Omega\|_{op}  \frac{\mathrm{r}(\Omega^2)}{\mathrm{r}(\Omega)},
	\end{align*}
	and by Theorem 5.48~\cite{vershynin2011introduction},
	\begin{align*}
		\|\breve{\Omega}_n- \Omega\|_{2 \rightarrow 2} = \|\breve{\Omega}_n- \Omega\|_{op}  = O_p\left( \|\Omega\|_{op} \left(\sqrt{ \frac{M_n(\log n)}{n}}  \vee  \frac{M_n(\log n)}{n} \right) \right).
	\end{align*}
	Moreover, as in the proof of Lemma~\ref{lemma:example:LogConcaveData} (i), Gaussian symmetrization yields,
	\begin{align*}
		\|\breve{X}_n \breve{X}_n'\|_{2 \rightarrow 2} = \|\breve{X}_n \breve{X}_n'\|_{op}  =  \left\|\frac{1}{n}\sum_{i=1}^n R(X_i - \mu_n) \right\|_2^2 = O_p \left( \|\Omega\|_{op} \frac{r(\Omega)}{n}\right).
	\end{align*}
	Notice that $r(\Omega) \leq M_n$. Hence, above three inequalities imply that Assumptions~\ref{assumption:ControlThirdMoments}--\ref{assumption:ConsistentCovariance} hold under the stated conditions.
	
	Case (ii) follows from similar arguments as used in the proof of Lemma~\ref{lemma:example:LogConcaveData} (ii). We omit the repetitive details.
\end{proof}

\subsection{Proofs of Section~\ref{subsec:AsympSizeLevel}}

\begin{proof}[\textbf{Proof of Theorem~\ref{theorem:SizeAlphaTest}}]
	Trivial. Apply Theorem~\ref{theorem:Bootstrap-Lp-Norm-Quantiles} with $\Theta_n \equiv 0$ for all $n \geq 1$ and simplify using Assumptions~\ref{assumption:ControlThirdMoments}--\ref{assumption:ConsistentCovariance}.
\end{proof}

\begin{proof}[\textbf{Proof of Corollary~\ref{corollary:theorem:SizeAlphaTest}}]
	Apply Theorem~\ref{theorem:Bootstrap-Lp-Norm-Quantiles} with 
	\begin{align*}
		\Theta_n : = \| \widehat{R}_n(S_n - \sqrt{n} \mu_0)\|_p  -   \| R S_n - \sqrt{n} r)\|_p.
	\end{align*}
	Note that by the reverse triangle inequality
	\begin{align*}
		|\Theta_n| \leq \| (\widehat{R}_n- R)(S_n - \sqrt{n} \mu_0)\|_p \leq \|\widehat{R}_n - R\|_{q \rightarrow p} \|S_n -\sqrt{n}  \mu_0\|_p.
	\end{align*}
	Then, proceed as in the proof of Theorem~\ref{theorem:SizeAlphaTest}.
\end{proof}

\subsection{Proofs of Section~\ref{subsec:AsympCorrectness}}
\begin{proof}[\textbf{Proof of Theorem~\ref{theorem:UnbiasedTest}}]
	We begin with the following three facts: First, since the $X_i$'s are i.i.d. and have a log-concave density and since log-concavity is preserved under affine transformations and convolutions, the density of $n^{-1/2} \sum_{i=1}^n R(X_i-\mu)$ is also log-concave~\citep[e.g.][Proposition 3.1 and 3.5]{saumard2014log-concavity}. Second, by monotonicity of the logarithm, a log-concave function has convex super-level sets and is therefore quasi-concave. Third, since the $X_i$'s are symmetric around their mean $\mu$,  $n^{-1/2} \sum_{i=1}^n R(X_i-\mu)$ is symmetric around $0$.
	
	Define (to shorten the notation of the main paper)
	\begin{align*}
		c_p^*(\alpha) &:= \inf \left\{s \geq 0: \mathrm{P}\left(\|Z\|_p \leq s \right)  \geq \alpha \right\}, \quad{}\text{and}\\
		c_{n,p}^*(\alpha) &:= \inf \left\{s \geq 0: \mathrm{P}\left(\|Z_n\|_p\leq s \mid X_1, \ldots, X_n \right)  \geq \alpha \right\},
	\end{align*}
	where $Z \sim N(0, \Omega)$ and $Z_n \sim N(0, \widehat{\Omega}_n)$ where $\widehat{\Omega}_n = \widehat{\Omega}_n(X_1, \ldots, X_n; R)$ is an estimate of $\Omega = R \Sigma R'$.
	
	Now, let $\alpha \in (0, 1)$ and $\delta > 0$ be arbitrary. Set $\delta_n := \delta^3/K^3 \mathrm{Var}(\|Z\|_p)$. Fix a sequence $\mu \in \mathcal{H}_1$. As in the proof of Theorem~\ref{theorem:ConsistencyLocalAlternatives} (iii) we have, by inequality~\ref{eq:theorem:ConsistencyLocalAlternatives-0-1},
	\begin{align}\label{eq:theorem:UnbiasedTest-1}
		\mathrm{P}_\mu\left(T_{n,p} > c_{n,p}^*(1-\alpha)\right)  &\geq \mathrm{P}_\mu\left(T_{n,p} > c_p^*(\pi_{n,p}(\delta_n) + 1- \alpha) \right) - \mathrm{P}\left(\|\widehat{\Omega}_n- \Omega\|_{q \rightarrow p} > \delta_n \right)\nonumber\\
		&= \mathrm{P}_\mu\left(T_{n,p} > c_p^*(\pi_{n,p}(\delta_n) + 1- \alpha) \right) + o(1).
	\end{align}
	Whence, we compute
	\begin{align}\label{eq:theorem:UnbiasedTest-2}
		&\mathrm{P}_{\mu_n}\left(T_{n,p} > c_p^*(\pi_{n,p}(\delta_n) + 1- \alpha) \right)\nonumber\\
		&\quad{} = \mathrm{P}_{\mu_n}\left(\left\|\frac{1}{\sqrt{n}}\sum_{i=1}^n R(X_i - \mu) + \sqrt{n} (R\mu - r) \right\|_p > c_p^*(\pi_{n,p}(\delta_n) + 1- \alpha) \right) \nonumber\\
		&\quad{} \overset{(a)}{\geq} \mathrm{P}_{\mu_n}\left(\left\|\frac{1}{\sqrt{n}}\sum_{i=1}^n R(X_i - \mu) \right\|_p > c_p^*(\pi_{n,p}(\delta_n) + 1- \alpha) \right)\nonumber\\
		&\quad{} \geq \mathrm{P}\big(\|Z\|_p > c_p^*(\pi_{n,p}(\delta_n) + 1- \alpha) \big) - \sup_{s \geq 0} \left|\mathrm{P}_\mu\left(\left\|R (S_n - \sqrt{n} \mu)\right\|_p\leq s \right) - \mathrm{P}\left(\|Z\|_p \leq s\right)  \right|\nonumber\\
		&\quad{} = 1- (\pi_p(\delta_n) + 1- \alpha) + o(1)\nonumber\\
		&\quad{} = \alpha - \delta + o(1),
	\end{align}
	where (a) by Anderson's lemma~\cite[][Corollary 2]{anderson1955integral} and the three facts stated at the beginning. We combine eq.~\eqref{eq:theorem:UnbiasedTest-1} and~\eqref{eq:theorem:UnbiasedTest-2} to conclude that
	\begin{align*}
		\inf_{ \mu \in \mathcal{H}_1} \mathrm{P}_{\mu_n}\left(T_{n,p} > c_{n,p}^*(1-\alpha)\right) \geq \alpha - \delta + o(1).
	\end{align*}
	To complete the proof, take $n \rightarrow \infty$ followed by $\delta \rightarrow 0$.
\end{proof}

\begin{proof}[\textbf{Proof of Theorem~\ref{theorem:ConsistencyLocalAlternatives}}]
	Define (to shorten the notation of the main paper)
	\begin{align*}
		c_p^*(\alpha) &:= \inf \left\{s \geq 0: \mathrm{P}\left(\|Z\|_p \leq s \right)  \geq \alpha \right\}, \quad{}\text{and}\\
		c_{n,p}^*(\alpha) &:= \inf \left\{s \geq 0: \mathrm{P}\left(\|Z_n\|_p\leq s \mid X_1, \ldots, X_n \right)  \geq \alpha \right\},
	\end{align*}
	where $Z \sim N(0, \Omega)$ and $Z_n \sim N(0, \widehat{\Omega}_n)$ where $\widehat{\Omega}_n = \widehat{\Omega}_n(X_1, \ldots, X_n; R)$ is an estimate of $\Omega = R \Sigma R'$.
	
	Consider statement (i). 
	Fix $(\mu_n)_{n \in \mathbb{N}} \in \mathcal{A}_p$ and write
	\begin{align*}
		\mathrm{P}_{\mu_n}\left(T_{n,p} > c_{n,p}^*(1- \alpha) \right) &= \mathrm{P}_{\mu_n}\left(\left\|\frac{1}{\sqrt{n}}\sum_{i=1}^n (R X_i - r) \right\|_p  > c_{n,p}^*(1- \alpha) \right)\\
		&= \mathrm{P}_{\mu_n}\left( \left\| R(S_n - \sqrt{n}\mu_n)\right\|_p + \Theta_n  > c_{n,p}^*(1- \alpha)\right),
	\end{align*}
	where $\Theta_n : = \left\|RS_n - \sqrt{n}r\right\|_p - \left\| R(S_n - \sqrt{n}\mu_n)\right\|_p$.
	Notice that
	\begin{align*}
		\Theta_n \leq \left| \left\|RS_n - \sqrt{n}r\right\|_p - \left\| R(S_n - \sqrt{n}\mu_n)\right\|_p \right| \leq \sqrt{n}\left\|R\mu_n - r\right\|_p \overset{(a)}{=} o\Big(\sqrt{\mathrm{Var}(\|Z\|_p)}\Big),
	\end{align*}
	where (a) holds by definition of $\mathcal{C}_p$. Thus, by Theorem~\ref{theorem:Bootstrap-Lp-Norm-Quantiles},
	\begin{align*}
		\lim_{n \rightarrow \infty}\sup_{(\mu_n)_{n \in \mathbb{N}} \in \mathcal{A}_p }\mathrm{P}_{\mu_n}\left(T_{n,p} > c_{n,p}^*(1- \alpha) \right)  = \alpha.
	\end{align*}
	This proves statement (i).
	
	Next, consider statement (iii). For $1 \leq p, q \leq \infty$ conjugate exponents such that $ 1/p + 1/q = 1$, Lemma~\ref{lemma:Bootstrap-Sup-Empirical-Process-Quantil-Comparison} implies that
	\begin{align}
		&\sup_{\alpha \in (0,1)} \mathrm{P}\left( c_{n, p}^*(\alpha) \leq c_p^*(\pi_{n,p}(\delta) + \alpha) \right) \geq 1 - \mathrm{P}\left(\|\widehat{\Omega}_n- \Omega\|_{q \rightarrow p} > \delta \right),\label{eq:theorem:ConsistencyLocalAlternatives-0-1}\\
		&\sup_{\alpha \in (0,1)} \mathrm{P}\left(c_p^*(\alpha) \leq c_{n, p}^*(\pi_{n,p}(\delta) + \alpha) \right) \geq 1 - \mathrm{P}\left(\|\widehat{\Omega}_n- \Omega\|_{q \rightarrow p} > \delta \right),\label{eq:theorem:ConsistencyLocalAlternatives-0-2}
	\end{align}
	where $\pi_{n,p}(\delta) := K\delta^{1/3} \left(\mathrm{Var}(\|Z\|_p)\right)^{-1/3}$ and $K > 0$ is an absolute constant. 
	Fix $\alpha \in (0,1)$ and define $\delta_\alpha := \alpha^3/(2K)^3\mathrm{Var}(\|Z\|_p)$ such that $1/ \big(\alpha - \pi_p(\delta_\alpha)\big) \leq 2/\alpha \Leftrightarrow \pi_{n,p}(\delta_\alpha) \leq \alpha/2$. Fix a sequence $(\mu_n)_{n \in \mathbb{N}} \in \mathcal{C}_p$. Since $1 \geq \mathrm{P}(A) + \mathrm{P}(B) - \mathrm{P}(A \cap B)$ for any events $A$ and $B$, we have
	\begin{align}\label{eq:theorem:ConsistencyLocalAlternatives-1}
		&\mathrm{P}_{\mu_n}\left(T_{n,p} > c_{n,p}^*(1-\alpha)\right) \nonumber\\
		&\quad{} \geq \mathrm{P}_{\mu_n}\left(T_{n,p} > c_{n,p}^*(1-\alpha), \:  c_{n,p}^*(1-\alpha) \leq c_p^*(\pi_{n,p}(\delta_\alpha) + 1- \alpha)\right)  \nonumber\\
		&\quad{} \geq \mathrm{P}_{\mu_n}\left(T_{n,p} > c_p^*(\pi_{n,p}(\delta_\alpha) + 1- \alpha) \right)  + \mathrm{P} \left(c_{n,p}^*(1-\alpha) \leq c_p^*(\pi_{n,p}(\delta_\alpha) + 1- \alpha)\right) -1\nonumber\\
		&\quad{} \overset{(a)}{\geq} \mathrm{P}_{\mu_n}\left(T_{n,p} > c_p^*(\pi_{n,p}(\delta_\alpha) + 1- \alpha) \right) - \mathrm{P}\left(\|\widehat{\Omega}_n- \Omega\|_{q \rightarrow p} > \delta_\alpha \right)\nonumber\\
		&\quad{} = \mathrm{P}_{\mu_n}\left(T_{n,p} > c_p^*(\pi_{n,p}(\delta_\alpha) + 1- \alpha) \right) + o(1),
	\end{align}
	where (a) follows from inequality~\eqref{eq:theorem:ConsistencyLocalAlternatives-0-1}. We now lower bound the first term on the right hand side in above display. By the reverse triangle inequality and Lemma~\ref{lemma:UpperBoundQuantiles},
	\begin{align}\label{eq:theorem:ConsistencyLocalAlternatives-2}
		&\mathrm{P}_{\mu_n}\left(T_{n,p} > c_p^*(\pi_{n,p}(\delta_\alpha) + 1- \alpha) \right) \nonumber\\
		&\quad{} = \mathrm{P}_{\mu_n}\left(\left\|\frac{1}{\sqrt{n}}\sum_{i=1}^n R(X_i - \mu_n) + \sqrt{n}(R\mu_n - r)  \right\|_p > c_p^*(\pi_{n,p}(\delta_\alpha) + 1- \alpha)\right) \nonumber\\
		&\quad{} \geq \mathrm{P}_{\mu_n}\left(\left\|\frac{1}{\sqrt{n}}\sum_{i=1}^n R(X_i - \mu_n) \right\|_p < \sqrt{n}\|R \mu_n- r\|_p - c_p^*(\pi_{n,p}(\delta_\alpha) + 1- \alpha)\right)\nonumber\\
		&\quad{} \geq \mathrm{P}\left(\|Z\|_p  < \sqrt{n}\|R \mu_n- r\|_p - c_p^*(\pi_{n,p}(\delta_\alpha) + 1- \alpha) \right)\nonumber\\
		&\quad{}\quad{} - \sup_{s \geq 0} \left|\mathrm{P}_{\mu_n}\left(\left\|R (S_n - \sqrt{n} \mu_n)\right\|_p\leq s \right) - \mathrm{P}\left(\|Z\|_p \leq s\right)  \right|\nonumber\\
		&\quad{} \geq \mathrm{P}\left(\|Z\|_p < \sqrt{n}\|R\mu_n- r\|_p - \mathrm{E}[\|Z\|_p]  - \sqrt{1/\big(\alpha- \pi_{n,p}(\delta_\alpha)\big)\mathrm{Var}(\|Z\|_p)}\right) + o(1).
	\end{align}
	By Markov's inequality and the construction of $\delta_\alpha > 0$ the first term on the right hand side in above display can be bounded by
	\begin{align}\label{eq:theorem:ConsistencyLocalAlternatives-3}
		&\inf_{\mu_n\in \mathcal{C}_p}\mathrm{P}_{\mu_n}\left(\|Z\|_p < \sqrt{n}\|R \mu_n- r\|_p - \mathrm{E}[\|Z\|_p]  - \sqrt{2 \mathrm{Var}(\|Z\|_p) / \alpha} \right) \nonumber\\
		&\quad{} \overset{(a)}{\geq}\inf_{\mu_n\in \mathcal{C}_p}\mathrm{P}_{\mu_n}\left(\|Z\|_p < \sqrt{n}\|R \mu_n- r\|_p - C_\alpha \mathrm{E}[\|Z\|_p] \right) \nonumber\\
		&\quad{}\geq 1 -  \frac{(1 +C_\alpha)\mathrm{E}[\|Z\|_p]}{\sqrt{n}\|R \mu_n- r\|_p} \nonumber\\
		&\quad{} \overset{(b)}{\rightarrow} 1 \quad{} \mathrm{as} \quad{} n \rightarrow \infty,
	\end{align}
	where (a) holds by Proposition A.2.4 in~\cite{vandervaart1996weak} 
	and $C_\alpha > 0$ is an absolute constant depending only on $\alpha > 0$, and (b) holds since $(\mu_n)_{n \in \mathbb{N}} \in \mathcal{C}_p$. To conclude the proof of statement (iii) combine eq.~\eqref{eq:theorem:ConsistencyLocalAlternatives-1}--\eqref{eq:theorem:ConsistencyLocalAlternatives-3}.
	
	Lastly, consider statement (ii). Without loss of generality we can restrict our attention to alternatives $(\mu_n)_{n \geq 1} \in \mathcal{B}_p \cap \mathcal{A}_p^c$; the case of $(\mu_n)_{n \geq 1} \in \mathcal{B}_p \cap \mathcal{A}_p$ follows from statement (i). Let $\alpha \in (0, 1/2)$ be arbitrary, choose $\delta_\alpha := (1-2\alpha)^3/K^3 \mathrm{Var}(\|Z\|_p)$ such that $\pi_{n,p}(\delta_\alpha) = 1- 2\alpha$, and set $\beta := \alpha + \pi_{n,p}(\delta_\alpha)$. Fix a sequence $(\mu_n)_{n \in \mathbb{N}} \in \mathcal{B}_p\cap \mathcal{A}_p^c$. Analogous to eq.~\eqref{eq:theorem:ConsistencyLocalAlternatives-1} we have, by inequality~\eqref{eq:theorem:ConsistencyLocalAlternatives-0-2},
	\begin{align*}
		\mathrm{P}_{\mu_n}\left(T_{n,p} > c_p^*(1-\beta)\right) \geq \mathrm{P}_{\mu_n}\left(T_{n,p} > c_{n,p}^*(\pi_{n,p}(\delta_\alpha) + 1- \beta) \right) - \mathrm{P}\left(\|\widehat{\Omega}_n- \Omega\|_{q \rightarrow p} > \delta_\alpha\right),
	\end{align*}
	which is equivalent to
	\begin{align}\label{eq:theorem:ConsistencyLocalAlternatives-4}
		\mathrm{P}_{\mu_n}\left(T_{n,p} > c_{n,p}^*(1 - \alpha) \right) \leq \mathrm{P}_{\mu_n}\left(T_{n,p} > c_{n,p}(1- \pi_{n,p}^*(\delta_\alpha) - \alpha \right) + \mathrm{P}\left(\|\widehat{\Omega}_n- \Omega\|_{q \rightarrow p} > \delta_\alpha \right).
	\end{align}
	Now, we upper bound the first term on the right hand side in above display. By the triangle inequality and Lemma~\ref{lemma:UpperBoundQuantiles},
	\begin{align}\label{eq:theorem:ConsistencyLocalAlternatives-5}
		&\mathrm{P}_{\mu_n}\left(T_{n,p} > c_p^*(1- \pi_{n,p}(\delta_\alpha) - \alpha) \right) \nonumber\\
		&\quad{} = \mathrm{P}_{\mu_n}\left(\left\|\frac{1}{\sqrt{n}}\sum_{i=1}^n R(X_i - \mu_n) + \sqrt{n}(R\mu_n - r)  \right\|_p > c_p^*(1- \pi_{n,p}(\delta_\alpha) - \alpha)\right) \nonumber\\
		&\quad{} \leq \mathrm{P}_{\mu_n}\left(\left\|\frac{1}{\sqrt{n}}\sum_{i=1}^n R(X_i - \mu_n) \right\|_p > c_p^*(1- \pi_{n,p}(\delta_\alpha) - \alpha) - \sqrt{n}\|R \mu_n- r\|_p \right)\nonumber\\
		&\quad{} \leq \mathrm{P}\left(\|Z\|_p  > c_p^*(1- \pi_{n,p}(\delta_\alpha) - \alpha) - \sqrt{n}\|R \mu_n- r\|_p  \right)\nonumber\\
		&\quad{}\quad{} + \sup_{s \geq 0} \left|\mathrm{P}_{\mu_n}\left(\left\|R (S_n - \sqrt{n} \mu_n)\right\|_p\leq s \right) - \mathrm{P}\left(\|Z\|_p \leq s\right)  \right|\nonumber\\
		&\quad{} \leq \mathrm{P}\left(\|Z\|_p  > c_p^*(\alpha) - \sqrt{n}\|R \mu_n- r\|_p  \right)+ o(1)\nonumber\\
		&\quad{} \leq \mathrm{P}\left(\|Z\|_p > \mathrm{E}[\|Z\|_p]  - \sqrt{1/(1-\alpha)\mathrm{Var}(\|Z\|_p)} - \sqrt{n}\|R\mu_n- r\|_p \right) + o(1).
	\end{align}
	By definition $\mathcal{B}_p$ and since $\|Z\|_p$ has no point mass there exists $N > 0$ such that for all $n \geq N$ the first term on the far right hand side of eq.~\eqref{eq:theorem:ConsistencyLocalAlternatives-5} can be upper bounded by
	\begin{align}\label{eq:theorem:ConsistencyLocalAlternatives-6}
		&\sup_{\mu_n \in \mathcal{B}_p\cap \mathcal{A}_p^c} \mathrm{P}_{\mu_n}\left(\|Z\|_p >  \mathrm{E}[\|Z\|_p] - C_\alpha\sqrt{n}\|R\mu_n- r\|_p \right) \nonumber\\
		&\quad{}\leq \sup_{\mu_n \in \mathcal{B}_p\cap \mathcal{A}_p^c} \mathrm{P}_{\mu_n}\left(\|Z\|_p >  \widetilde{C}_\alpha \mathrm{E}[\|Z\|_p] \right) \nonumber\\
		&\quad{} < 1,
	\end{align}
	where $C_\alpha, \widetilde{C}_\alpha > 0$ are absolute constants depending on $\alpha > 0$ only. Combine eq.~\eqref{eq:theorem:ConsistencyLocalAlternatives-4}--\eqref{eq:theorem:ConsistencyLocalAlternatives-6} to complete the proof of statement (ii).
\end{proof}

\begin{proof}[\textbf{Proof of Proposition~\ref{lemma:example:SparseDenseAlternatives}}]
	Combine Lemma~\ref{lemma:EVNormGaussian} with Theorem~\ref{theorem:ConsistencyLocalAlternatives} (ii) and (iii).
\end{proof}

\subsection{Proofs of Section~\ref{subsec:BahadurSlope}}

\begin{proof}[\textbf{Proof of Theorem~\ref{theorem:BahadurARE}}]
	We begin with the following claim: Under Assumptions~\ref{assumption:ControlThirdMoments}--\ref{assumption:ConsistentCovariance} we can substitute the expected value of the Gaussian proxy statistic $\mathrm{E}[\|Z\|_p]$ for $\mathrm{E}[\|n^{-1/2}\sum_{i=1}^n R(X_i - \mu_n)\|_p ]$.
	
	Indeed, recall from~\cite{giessing2023Gaussian} that the proof of Theorem~\ref{theorem:CLT-Lp-Norm} (Gaussian approximation) relies on bounding
	\begin{align*}
		\left| \mathrm{E}\left[ h\left(\left\|\frac{1}{\sqrt{n} }\sum_{i=1}^n R(X_i - \mu_n)\right\|_p \right) - h\left(\|Z\|_p\right)\right] \right| = \left| \mathrm{E}\left[ h\left(T_{n,p}\right) - h\left(\|Z\|_p\right)\right] \right| \lesssim B_n
	\end{align*}
	for some smooth function $h$ with bounded derivatives and an upper bound $B_n > 0$. It is easy to verify that the proof (and hence the upper bound) of Theorem~\ref{theorem:CLT-Lp-Norm} also applies to the identity function $h(x) = x$. The claim now follows since $B_n \rightarrow 0$ as $n \rightarrow \infty$ whenever Assumptions~\ref{assumption:ControlThirdMoments}--\ref{assumption:ConsistentCovariance} hold. In particular, it follows that $\limsup_{n \rightarrow \infty} \mathrm{E}[\|n^{-1/2} \sum_{i=1}^n R(X_i - \mu_n) \|_p]/(\sqrt{n} \|R\mu_n - r\|_p) < 1$.
	
	Now, consider statement (i). By Lemma~\ref{lemma:LDP-Log-Concave} for $1 \leq p < \log d$,
	\begin{align*}
		- \frac{1}{\sqrt{n}}\log\Big( \sup_{ \mu \in \mathcal{H}_0}\big(1 - F_{n, \mu}(T_{n,p}) \big)\Big) & \gtrsim \frac{1}{\sqrt{n}}\frac{T_{n,p}}{\|\Omega^{1/2}\|_{2 \rightarrow p}} \\
		&\gtrsim \frac{ \|R\mu_n - r\|_p}{\|\Omega^{1/2}\|_{2 \rightarrow p}} - \frac{ \left|\frac{1}{\sqrt{n}} T_{n,p} - \|R\mu_n - r\|_p\right|}{\|\Omega^{1/2}\|_{2 \rightarrow p}} \\
		&\overset{(a)}{\gtrsim} \frac{ \|R\mu_n - r\|_p}{\|\Omega^{1/2}\|_{2 \rightarrow p}} + O_p\left( \frac{ \mathrm{E}\left[ \|Z\|_p \right] }{ n^{1/2}\|\Omega^{1/2}\|_{2 \rightarrow p}} \right) \\ 
		&\overset{(b)}{\gtrsim} \frac{ \|R\mu_n - r\|_p}{\|\Omega^{1/2}\|_{2 \rightarrow p}} + o_p(1), 
	\end{align*}
	where (a) holds since by the reverse triangle inequality and above claim we have
	\begin{align*}
		\mathrm{E}\left[\left|\frac{1}{\sqrt{n}} T_{n,p} - \|R\mu_n - r\|_p\right| \right] \leq \mathrm{E}\left[\left\|\frac{1}{n} \sum_{i=1}^n R(X_i - \mu_n)\right\|_p \right] \leq \frac{1}{\sqrt{n}}\mathrm{E}[\|Z\|_p] + o(n^{-1/2}),
	\end{align*}	
	and (b) follows from Assumptions~\ref{assumption:ControlThirdMoments}--\ref{assumption:ConsistentCovariance} and since by the Gaussian Poincar{\'e} inequality $\mathrm{Var}(\|Z\|_p) \leq \|\Omega^{1/2}\|_{2 \rightarrow p}$.
	
	The proofs of statements (ii) and (iii) follow in the same way by substituting Lemma~\ref{lemma:LDP-Log-Concave} with Lemma~\ref{lemma:LDPStrictly-Log-Concave} and Lemma~\ref{lemma:LDPGaussian}, respectively.
\end{proof}

\begin{proof}[\textbf{Proof of Proposition~\ref{lemma:example:BahadurSlope}}]
	Since the Bahadur slope is an asymptotic concept for $n,d,t \rightarrow \infty$, we can assume that  $\log t \geq 2$. Also, denote by $\omega_{(t)}^2$ the largest diagonal entry of $\Omega$. 
	
	Proof of statement (i). Let $p \in [\log t, \infty]$ be arbitrary. By Theorem~\ref{theorem:ConsistencyLocalAlternatives}, $\mathcal{D}_{\delta, s} \subseteq \mathcal{C}_{\log t} \cap \mathcal{C}_p $. Moreover, $\|\Omega^{1/2}\|_{2 \rightarrow p} = \sup_{\|u\|_2 = 1}\|\Omega^{1/2}u\|_p = \omega_{(t)} \sup_{\|u\|_2=1} \|u\|_p$ and hence for all $p \in [2, \infty]$, 
	\begin{align}\label{eq:lemma:example:BahadurSlope-1}
		1 = \sup_{\|u\|_2=1} \|u\|_\infty = \|\Omega^{1/2}\|_{2 \rightarrow \infty}/ \omega_{(t)} \leq \|\Omega^{1/2}\|_{2 \rightarrow p}/ \omega_{(t)} \leq \|\Omega^{1/2}\|_{2 \rightarrow 2}/ \omega_{(t)} = \sup_{\|u\|_2=1} \|u\|_2 = 1.
	\end{align}
	Now, the statement (i) follows from the fact that $\|R\mu_n-r\|_{\log t} \geq \|R\mu_n-r\|_p$ for all $p \in [\log t, \infty]$ and all $(\mu_n)_{n \geq 1} \in \mathcal{D}_{\delta,s}$.
	
	Proof of statement (ii). First, consider $p \in [2, q)$. By Theorem~\ref{theorem:ConsistencyLocalAlternatives}, $\mathcal{D}_{\delta, s} \subseteq \mathcal{C}_2 \cap \mathcal{C}_p $. Moreover, eq.~\eqref{eq:lemma:example:BahadurSlope-1} continues to hold and $\|R\mu_n-r\|_2 \geq \|R\mu_n-r\|_p$ for all $p \in [2, q)$ and $(\mu_n)_{n \geq 1} \in \mathcal{D}_{\delta,s}$. This proves statement (ii) for $p \in [2, q)$. Next, consider $p \in [1, 2)$. Notice that $\|u\|_2 \leq \|u\|_p \leq d^{1/p -1/2}\|u\|_2$ and $\|\Omega^{1/2}\|_{2 \rightarrow p}=\omega_{(t)} \sup_{\|u\|_2=1}\|u\|_p =\omega_{(t)} d^{1/p -1/2}$. Thus, for all $(\mu_n)_{n \geq 1} \in \mathcal{D}_{\delta, s}$,
	\begin{align*}
		\frac{\|R\mu_n-r\|_p}{\|\Omega^{1/2}\|_{2 \rightarrow p}} = \frac{\|R\mu_n-r\|_p}{\omega_{(t)} d^{1/p - 1/2}}\leq \frac{\|R\mu_n-r\|_2}{\omega_{(t)}} = \frac{\|R\mu_n-r\|_2}{\|\Omega^{1/2}\|_{2 \rightarrow 2}}.
	\end{align*}
	This proves statement (ii) for all $p \in [1, 2)$. 
\end{proof}

\subsection{Proofs of Section~\ref{subsec:PowerEnhancement}}

\begin{proof}[\textbf{Proof of Proposition~\ref{theorem:SizeAlphaTestL1LInfty}}]
	Note that 
	\begin{align*}
		W_n := T_{n,2} + T_{n,\log t} &= \sup_{\|u\|_2 =1} (RS_n - \sqrt{n} r)'u + \sup_{\|v\|_{\log t/(\log t - 1)} =1} (RS_n - \sqrt{n} r)'v\\
		& = \sup_{ \substack{\|u\|_2 = 1,\\ \|v\|_{\log t/(\log t - 1)} = 1}} (RS_n - \sqrt{n} r)'(u + v).
	\end{align*}
	Under the null hypothesis above supremum is a centered empirical process indexed by functions
	\begin{align*}
		\mathcal{F}_n = \left\{f\equiv f_{u,v} : f_{u,v}(x) = x'(u +v), \: \|u\|_2 =1, \:  \|v\|_{\log t/(\log t - 1)} = 1, \: u,v \in \mathbb{R}^t \right\},
	\end{align*}
	with envelop $F(x) = \|x\|_2 + \|x\|_{\log t}$. By Harris' association inequality $\mathrm{Cov}(\|Z\|_2, \|Z\|_{\log t})\geq 0$ and hence $\mathrm{Var}(W_n^*) \geq \mathrm{Var}(\|Z\|_2) + \mathrm{Var}(\|Z\|_{\log t})$. Thus, under the stated assumptions the claim follows from Theorem~\ref{theorem:Bootstrap-Lp-Norm-Quantiles}.
\end{proof}

\begin{proof}[\textbf{Proof of Proposition~\ref{theorem:ConsistencySizeAlphaTestL1LInfty}}]
	The claim follows by similar arguments as those used to prove Theorem~\ref{theorem:ConsistencyLocalAlternatives} (iii). Note that we apply Lemma~\ref{lemma:UpperBoundQuantiles} with $\sqrt{\mathrm{Var}(\|Z\|_2 + \|Z\|_{\log t})}$ and then use Cauchy-Schwarz (or triangle inequality) $\sqrt{\mathrm{Var}(\|Z\|_2 + \|Z\|_{\log t})} \leq \sqrt{\mathrm{Var}(\|Z\|_2)} + \sqrt{\mathrm{Var}(\|Z\|_{\log t})}$ to simplify the expression.
\end{proof}

\subsection{Proofs of Section~\ref{subsec:EllipticallyDistributedData}}
\begin{proof}[\textbf{Proof of Proposition~\ref{theorem:SizeAlphaTestEllipticallyDistributed}}]
	Throughout the proof, we denote by $r(M) := \mathrm{tr}(M)/\|M\|_{op}$ the effective rank of a matrix $M$. Moreover, we write $Z \sim N(0, \Omega)$ and $\Omega = \mathrm{E}[\widetilde{X}\widetilde{X}']$.	Denote by $R_k$ the $k$th row in $R$ and compute
	\begin{align*}
		\mathrm{tr}( \Omega ) = \mathrm{tr}\left( \mathrm{E}\left[ \frac{(RX- r)(RX-r)'}{\|RX - r\|_2^2}\right]\right) = \sum_{k=1}^d  \mathrm{E}\left[ \frac{(R_k'X - r_k)^2}{\|RX-r\|_2^2}\right] = 1.
	\end{align*}
	Therefore, by Theorem~\ref{theorem:LowerBoundVariance-LpNorm-Refinements},
	\begin{align*}
		\mathrm{Var}(\|Z\|_2) \geq \frac{\mathrm{tr}(\Omega^2)}{\mathrm{tr}(\Omega)} =  \|\Omega\|_F^2,
	\end{align*}
	and, hence,
	\begin{align*}
		\frac{ \mathrm{E}[\|Z\|_2]}{\sqrt{\mathrm{Var}(\|Z\|_2)}} \leq \frac{\mathrm{tr}(\Omega)}{\sqrt{\mathrm{tr}(\Omega^2)}} = \frac{1}{\|\Omega\|_F}. 
	\end{align*}
	Also,
	\begin{align*}
		\|\Omega\|_{op} =  \sup_{\|u\|_2 = 1}\left\| \mathrm{E}\left[ \frac{(RX- r)(RX-r)'}{\|RX - r\|_2^2}\right] u\right\|_2 \leq  \sup_{\|u\|_2 = 1} \mathrm{E}\left[ \frac{|(RX-r)'u|}{\|RX - r\|_2} \right]  \leq 1.
	\end{align*}
	Hence, by Theorem 5.48~\cite{vershynin2011introduction},
	\begin{align*}
		\|\widetilde{\Omega}_n- \Omega\|_{2 \rightarrow 2} = \|\widetilde{\Omega}_n- \Omega\|_{op} & = O_p\left( \|\Omega\|_{op} \left(\sqrt{ \frac{\mathrm{r}(\Omega)(\log n)}{n}}  \vee  \frac{\mathrm{r}(\Omega)(\log n)}{n} \right) \right)=  O_p\left( \sqrt{ \frac{\log n}{n}} \right),
	\end{align*}
	where we have used that  $r(\Omega) = \mathrm{tr}(\Omega)/\|\Omega\|_{op} = 1/\|\Omega\|_{op}$.
	
	Since $\|\widetilde{X}\|_2 = 1$ and, under the null hypothesis, $\mathrm{E}[\widetilde{X}] = 0$ (because $X$ is elliptically distributed), it follows from Theorem~\ref{theorem:Bootstrap-Lp-Norm-Quantiles} and above bounds on trace, variance and operator norm that there exists an absolute constant $C > 0$ such that, for all $n > 1$,
	\begin{align*}
		&\sup_{\alpha \in (0,1)} \Big|\mathrm{P}\left(V_n \leq c_V^*(\alpha; \widetilde{\Omega}_n) \right)  - \alpha \Big| \\
		& \lesssim \frac{1}{n^{1/6}\sqrt{\mathrm{Var}(\| Z\|_2)}} +  \frac{ \mathrm{E}[\|Z\|_2]}{\sqrt{n\mathrm{Var}(\|Z\|_2)}} + \left( \frac{\sqrt{\log n}}{\sqrt{n}\mathrm{Var}( \|Z\|_2)}\right)^{1/3}  + \mathrm{P}\left(\frac{ \|\widetilde{\Omega}_n- \Omega\|_{2 \rightarrow 2}}{\mathrm{Var}( \|Z\|_2)} >  \frac{C\sqrt{\log n}}{\sqrt{n}\mathrm{Var}( \|Z\|_2)} \right)\\
		&\lesssim \frac{1}{n^{1/6}\|\Omega\|_F} +  \frac{1}{\sqrt{n}\|\Omega\|_F} + \left(\frac{\log n}{n \|\Omega\|_F^2}  \right)^{1/6}  + o(1).
	\end{align*}
	To conclude the proof, we lower bound $\|\Omega\|_F$. We have
	\begin{align*}
		\|RX - r\|_2^2 = \theta^2 {U^{(s)}}'\Gamma'R'R \Gamma U^{(s)} \leq \theta^2 \lambda_{\max}(\Gamma'R'R\Gamma) = \theta^2 \lambda_{\max}(R\Gamma\Gamma'R') = \theta^2 \|R\Gamma\Gamma'R'\|_{op},
	\end{align*}
	where $\lambda_{\max}(M)$ denotes the largest eigenvalue of $M$ and the last equality holds because $R\Gamma\Gamma'R'$ is normal. Now, compute
	\begin{align*}
		\|\Omega\|_F^2 = \sum_{k=1}^d \sum_{j =1}^d \left(\mathrm{E} \left[\frac{X_kX_j}{\|X\|_2^2}\right] \right)^2 &\geq \sum_{k=1}^d \sum_{j =1}^d \left( \frac{\mathrm{E}[R_k'\Gamma  U^{(s)} R_j'\Gamma U^{(s)}]}{\|R\Gamma\Gamma'R'\|_{op}} \right)^2\\
		& =  \sum_{k=1}^d \sum_{j =1}^d \left( \frac{s^{-1}R_k'\Gamma\Gamma' R_j}{\|R\Gamma\Gamma'R'\|_{op}} \right)^2\\
		& = \left(\frac{\|R\Gamma\Gamma' R'\|_F}{s\|R\Gamma\Gamma'R'\|_{op}}\right)^2.
	\end{align*}
	Combine this lower bound with the assumptions in the theorem to conclude.
\end{proof}

\begin{proof}[\textbf{Proof of Proposition~\ref{theorem:ConsistencySizeAlphaTestEllipticallyDistributed}}]
	The proof follows closely the one of Theorem~\ref{theorem:ConsistencyLocalAlternatives} (iii). We only point out the parts that need to be modified. 
	
	For $\alpha \in (0,1)$ arbitrary define (to shorten the notation of the main paper)
	\begin{align*}
		c^*(\alpha) &:= \inf \left\{s \geq 0: \mathrm{P}\left(\|Z\|_2 \leq s \right)  \geq \alpha \right\}, \quad{}\text{and}\\
		c_n^*(\alpha) &:= \inf \left\{s \geq 0: \mathrm{P}\left(\|Z_n\|_2\leq s \mid X_1, \ldots, X_n \right)  \geq \alpha \right\},
	\end{align*}
	where $Z \sim N(0, \Omega)$ and $Z_n \sim N(0, \widetilde{\Omega}_{\mu_n})$.
	By eq.~\ref{eq:theorem:ConsistencyLocalAlternatives-1} in the proof of Theorem~\ref{theorem:ConsistencyLocalAlternatives} (iii) we have
	\begin{align}\label{eq:theorem:ConsistencySizeAlphaTestEllipticallyDistributed-1}
		&\mathrm{P}_{\mu_n}\left(V_n > c_n^*(1-\alpha)\right) \leq \mathrm{P}_{\mu_n}\left(V_n > c^*(\pi_{n,2}(\delta_\alpha) + 1- \alpha) \right) + o(1),
	\end{align}
	where $\pi_{n,2}(\delta) = K\delta^{1/3} \left(\mathrm{Var}(\|Z\|_2)\right)^{-1/3}$ and $K > 0$ is an absolute constant.
	
	Now, by an applications of the reverse triangle inequality and Lemma~\ref{lemma:UpperBoundQuantiles},
	\begin{align}\label{eq:theorem:ConsistencySizeAlphaTestEllipticallyDistributed-3}
		&\mathrm{P}_{\mu_n}\left(V_n > c^*(\pi_{n,2}(\delta_\alpha) + 1- \alpha) \right) \nonumber\\
		&\quad{} = \mathrm{P}_{\mu_n}\left(\left\|\frac{1}{\sqrt{n}}\sum_{i=1}^n \frac{R(X_i - \mu_n)}{\|R(X_i - \mu_n)\|_2} + \frac{1}{\sqrt{n}}\sum_{i=1}^n \frac{RX_i -r }{\|RX_i - r\|_2} \right. \right. \nonumber \\
		&\quad{}\quad{}\quad{}\quad{}\quad{}\left. \left. - \frac{1}{\sqrt{n}}\sum_{i=1}^n \frac{R(X_i - \mu_n)}{\|R(X_i - \mu_n)\|_2} \right\|_2 > c^*(\pi_{n,2}(\delta_\alpha) + 1- \alpha)\right) \nonumber\\
		\begin{split}
			&\quad{} \geq \mathrm{P}_{\mu_n}\left(\left\|\frac{1}{\sqrt{n}}\sum_{i=1}^n \frac{R(X_i - \mu_n)}{\|R(X_i - \mu_n)\|_2} \right\|_2 \right.\\
			&\quad{}\quad{}\quad{} \quad{} \quad{} \left.  < \left\| \frac{1}{\sqrt{n}}\sum_{i=1}^n \frac{RX_i - r}{\|RX_i - r\|_2} - \frac{1}{\sqrt{n}}\sum_{i=1}^n \frac{R(X_i - \mu_n)}{\|R(X_i - \mu_n)\|_2}\right\|_2 - c^*(\pi_{n,2}(\delta_\alpha) + 1- \alpha)\right).
		\end{split}
	\end{align}
	For $t > 0$ arbitrary, denote by $\Omega_t$ the event
	\begin{align*}
		&\left\| \frac{1}{\sqrt{n}}\sum_{i=1}^n \frac{RX_i - r}{\|RX_i - r\|_2} - \frac{1}{\sqrt{n}}\sum_{i=1}^n \frac{R(X_i - \mu_n)}{\|R(X_i - \mu_n)\|_2}\right\|_2\\
		&\quad{}\leq \mathrm{E}_{\mu_n}\left[ \left\| \frac{1}{\sqrt{n}}\sum_{i=1}^n \frac{RX_i - r}{\|RX_i - r\|_2} - \frac{1}{\sqrt{n}}\sum_{i=1}^n \frac{R(X_i - \mu_n)}{\|R(X_i - \mu_n)\|_2}\right\|_2\right] -t.
	\end{align*}
	Since $\mathrm{P}(A) \geq \mathrm{P}(A \cap B) \geq \mathrm{P}(C \cap B) \geq P(C) - P(B^c)$ for arbitrary events $A, B$ and $C \subseteq A$, we can lower bound the probability in eq.~\eqref{eq:theorem:ConsistencySizeAlphaTestEllipticallyDistributed-3} by
	\begin{align}\label{eq:theorem:ConsistencySizeAlphaTestEllipticallyDistributed-4}
		&\mathrm{P}_{\mu_n}\left( \Omega_t^c, \: \left\|\frac{1}{\sqrt{n}}\sum_{i=1}^n \frac{R(X_i - \mu_n)}{\|R(X_i - \mu_n)\|_2} \right\|_2 \right.\nonumber \\
		&\quad{}\quad{}\quad{}\left.  < \left\| \frac{1}{\sqrt{n}}\sum_{i=1}^n \frac{RX_i - r}{\|RX_i - r\|_2} - \frac{1}{\sqrt{n}}\sum_{i=1}^n \frac{R(X_i - \mu_n)}{\|R(X_i - \mu_n)\|_2}\right\|_2 - c^*(\pi_{n,2}(\delta_\alpha) + 1- \alpha)\right)\nonumber \\
		&\geq\mathrm{P}_{\mu_n}\left( \Omega_t^c, \: \left\|\frac{1}{\sqrt{n}}\sum_{i=1}^n \frac{R(X_i - \mu_n)}{\|R(X_i - \mu_n)\|_2} \right\|_2 \right.\nonumber \\
		&\quad{}\quad{}\quad{} \left.  < \mathrm{E}_{\mu_n}\left[ \left\| \frac{1}{\sqrt{n}}\sum_{i=1}^n \frac{RX_i - r}{\|RX_i - r\|_2} - \frac{1}{\sqrt{n}}\sum_{i=1}^n \frac{R(X_i - \mu_n)}{\|R(X_i - \mu_n)\|_2}\right\|_2\right] - t - c^*(\pi_{n,2}(\delta_\alpha) + 1- \alpha)\right)\nonumber \\
		\begin{split}
			&\geq \mathrm{P}_{\mu_n}\left(\left\|\frac{1}{\sqrt{n}}\sum_{i=1}^n \frac{R(X_i - \mu_n)}{\|R(X_i - \mu_n)\|_2} \right\|_2 \right.\\
			&\quad{}\quad{}\quad{} \left.  < \sqrt{n} \left\| \mathrm{E}_{\mu_n}\left[\frac{RX - r}{\|RX - r\|_2} -  \frac{R(X - \mu_n)}{\|R(X - \mu_n)\|_2} \right] \right\|_2 - t - c^*(\pi_{n,2}(\delta_\alpha) + 1- \alpha)\right)\\
			&\quad{} - \mathrm{P}_{\mu_n} (\Omega_t).
		\end{split}
	\end{align}
	Using the classical bounded differences inequality~\citep[e.g.][Theorem 6.2]{boucheron2013concentration} we easily find that $\mathrm{P}_{\mu_n} (\Omega_t) \leq e^{-t^2/8}$ for all $t > 0$. Moreover, $\mathrm{E}_{\mu_n}\left[ R(X - \mu_n)/\|R(X - \mu_n)\|_2\right] = 0$ because the data are elliptically distributed with mean $\mu_n \in \mathbb{R}^d$. Also, by construction of $\delta_\alpha$, Lemma~\ref{lemma:UpperBoundQuantiles}, and Proposition A.2.4 in~\cite{vandervaart1996weak}
	\begin{align*}
		c^*(\pi_{n,2}(\delta_\alpha) + 1- \alpha) \leq \mathrm{E}[\|Z\|_2] + \sqrt{1/\big(\alpha- \pi_{n,2}(\delta_\alpha)\big)\mathrm{Var}(\|Z\|_p)} \leq C_\alpha  \mathrm{E}[\|Z\|_2] \leq C_\alpha,
	\end{align*}
	where $C_\alpha > 0$ is an absolute constant depending only on $\alpha > 0$. Thus, setting $t = \frac{\sqrt{n}}{2}\| \mathrm{E}_{\mu_n}[(RX - r)/ \|RX - r\|_2]\|_2]$ we can lower bound eq.~\eqref{eq:theorem:ConsistencySizeAlphaTestEllipticallyDistributed-4} by
	\begin{align*}
		&\mathrm{P}_{\mu_n} \left(\|Z\|_2 <  \frac{\sqrt{n}}{2} \left\| \mathrm{E}_{\mu_n}\left[\frac{RX - r}{\|RX - r\|_2}\right] \right\|_2 - C_\alpha \right)\\
		&\quad{}  - \sup_{s \geq 0} \left|\mathrm{P}_{\mu_n}\left(\left\|R (S_n - \sqrt{n} \mu_n)\right\|_2\leq s \right) - \mathrm{P}\left(\|Z\|_2 \leq s\right)  \right| -\exp\left( -\frac{n}{32}\left\| \mathrm{E}_{\mu_n}\left[\frac{RX - r}{\|RX - r\|_2}\right] \right\|_2^2 \right)\\
		& \geq 1 - \frac{2 + 2C_\alpha}{ \delta_n } - \sup_{s \geq 0} \left|\mathrm{P}_{\mu_n}\left(\left\|R (S_n - \sqrt{n} \mu_n)\right\|_2\leq s \right) - \mathrm{P}\left(\|Z\|_2 \leq s\right)  \right| - e^{-\delta_n^2/32} \\
		&\rightarrow 1 \quad{} \mathrm{as} \quad{} n \rightarrow \infty,
	\end{align*}
	where $\delta_n :=  \sqrt{n} \left\| \mathrm{E}_{\mu_n}\left[\frac{RX - r}{\|RX - r\|_2}\right] \right\|_2 \rightarrow \infty$ since $(\mu_n)_{n \geq 1} \in \mathcal{E}_2$. This completes the proof.
\end{proof}

\subsection{Proofs of Section~\ref{subsec:ApproximateSampleAverages}}
\begin{proof}[\textbf{Proof of Proposition~\ref{theorem:ApproximateSampleAverages}}]
	Trivial. Redo the proofs of Theorems~\ref{theorem:SizeAlphaTest} and~\ref{theorem:ConsistencyLocalAlternatives} and Propositions~\ref{theorem:ConfidenceSets},~\ref{theorem:SizeAlphaTestL1LInfty}, and~\ref{theorem:ConsistencySizeAlphaTestL1LInfty} using Theorem~\ref{theorem:Bootstrap-Lp-Norm-Quantiles} with $\Theta_n = \|W_n\|_p$.
\end{proof}

\subsection{Proofs of Section~\ref{subsec:HeuristicsBias}}

\begin{proof}[\textbf{Proof of Theorem~\ref{theorem:BoundMonteCarloError}}]
	Let $T^*_{n,p}$ and $T_{n,p}^{**}$ be the random variables with quantile functions $u \mapsto c^*_p(u; \widehat{\Omega}_n)$ and $u \mapsto c^*_{p,B}(u; \widehat{\Omega}_n)$, respectively. Denote by $M_{n,p} = \esssup_{z, \mu} f_{n,p, \mu}(z)$, where $f_{n,p, \mu}$ is the density of the test statistics $T_{n,p}$ when the $X_i$'s have mean $\mu$. Write $\mathcal{X} = \{X_1, \ldots, X_n\}$ and $\mathcal{T} = \{T_{n,p,1}^*, \ldots, T_{n,p,B}\}$. Now, compute, for $\alpha \in (0,1)$ arbitrary,
	\begin{align*}
		&\left| \mathrm{E}_\mu\left[ \varphi_{\alpha,B}(T_{n,p},\widehat{\Omega}_n) \right] -\mathrm{E}_\mu\left[ \varphi_{\alpha}(T_{n,p},\widehat{\Omega}_n) \right] \right| \nonumber\\ 
		&\quad{}\quad{} = \left|  \mathrm{E}_\mu\left[ \int_{c^*_{p,B}(1- \alpha;\widehat{\Omega}_n)}^{c^*_p(1- \alpha;\widehat{\Omega}_n)} f_{n,p, \mu}(u) du\right]  \right|  \nonumber\\
		&\quad{}\quad{}\leq M_{n,p}\mathrm{E} \left[ \left|c^*_{p,B}(1- \alpha;\widehat{\Omega}_n) - c^*_p(1- \alpha;\widehat{\Omega}_n) \right|  \right] \nonumber\\
		&\quad{}\quad{} \leq M_{n,p} \mathrm{E}\big[\inf \mathrm{E}\left[ \left|U - V \right|\right] \big],
	\end{align*}
	where the last line holds by Theorem 8.1 in~\cite{major1978on} and the infimum is taken over all joint (product) probability distributions of $(U,V)$ such that $U \overset{d}{=} T_{n,p}^{**} - \mathrm{E}[T_{n,p}^*\mid \mathcal{X}]\mid \{\mathcal{X}, \mathcal{T} \}$ and $V \overset{d}{=} T_{n,p}^*  -  \mathrm{E}[T_{n,p}^*\mid \mathcal{X}]\mid \mathcal{X}$. By Theorem 1 in~\cite{fournier2015on}
	\begin{align*}
		\inf \mathrm{E}\left[ \left|U - V \right|\right] \lesssim B^{-1/2} (\log B)  \sqrt{ \mathrm{Var}( T_{n,p}^* \mid \mathcal{X})}.
	\end{align*}	
	This completes the proof. (As an aside, technically, Theorem 1 in~\cite{fournier2015on} does not yield above bound because, using their notation, it does not apply to the case $q = 2p$. However, as the authors point out on top of p. 709 it is easy to modify their proof for the case $q = 2p$ and, since our statistic is one-dimensional, i.e. $d = 1$, we only need to rework Step 2 on p. 717. This modification is trivial and we leave the details to the reader.)
\end{proof}

\subsection{Proofs of Section~\ref{subsec:CovarianceEstimation}}

\begin{proof}[\textbf{Proof of Lemma~\ref{lemma:example:ApproxSparseCovariance}}]
	We begin with the following general observation: For $\lambda > 0$ arbitrary such that $\|\widehat{\Omega}_n^{\mathrm{naive}} - \Omega\|_{1 \rightarrow \infty} \leq \lambda/2$, with probability one,
	\begin{align}\label{eq:lemma:example:ApproxSparseCovariance-1}
		\left\|  \mathcal{T}_{\lambda}^+(\widehat{\Omega}_n^{\mathrm{naive}}) -  \Omega\right\|_{op} \overset{(a)}{\leq} 2 \left\|  \mathcal{T}_{\lambda}(\widehat{\Omega}_n^{\mathrm{naive}}) -  \Omega\right\|_{op} \overset{(b)}{\leq} 8 R_\gamma \lambda^{1- \gamma},
	\end{align}
	where (a) follows from the triangle inequality~\citep[][p. 275]{avella-medina2018robust} and (b) from the first part of Theorem 6.27 in~\cite{wainwright2019high}. Let $\lambda > 0$ be arbitrary and compute
	\begin{align}\label{eq:lemma:example:ApproxSparseCovariance-2}
		&\mathrm{P}\left( \left\|  \mathcal{T}_{\lambda}^+(\widehat{\Omega}_n^{\mathrm{naive}}) -  \Omega\right\|_{op} > 8 R_\gamma \lambda^{1- \gamma} \right) \nonumber \\
		&\quad \overset{(a)}{\leq} \mathrm{P}\left( \left\|  \mathcal{T}_{\lambda}^+(\widehat{\Omega}_n^{\mathrm{naive}}) -  \Omega\right\|_{op} > 8 R_\gamma \lambda^{1- \gamma}, \: \|\widehat{\Omega}_n^{\mathrm{naive}} - \Omega\|_{1 \rightarrow \infty} > \lambda/2\right)\nonumber \\
		&\quad \leq \mathrm{P}\left( \|\widehat{\Omega}_n^{\mathrm{naive}} - \Omega\|_{1 \rightarrow \infty} > \lambda/2\right),
	\end{align}
	where (a) follows from inequality~\eqref{eq:lemma:example:ApproxSparseCovariance-1}. 
	
	Now, consider case (i). By Markov's inequality applied to eq.~\eqref{eq:lemma:example:ApproxSparseCovariance-2} and the arguments of the proof of case (ii) in Lemma~\ref{lemma:example:SubGaussianData}, we have
	\begin{align*}
		\|\widehat{\Omega}_n^{\mathrm{naive}} - \Omega\|_{1 \rightarrow \infty} = O_p\left(  \omega_{(t)}^2 \left( \sqrt{\frac{ \log t}{n}} \vee \frac{ \log t}{n} \right) \right).
	\end{align*}
	Above two inequalities combined imply that
	\begin{align*}
		\left\|  \mathcal{T}_{\lambda}^+(\widehat{\Omega}_n^{\mathrm{naive}}) -  \Omega\right\|_{op} = O_p \left(  R_\gamma \omega_{(t)}^{2(1-\gamma)} \left( \sqrt{\frac{ \log t}{n}} \vee \frac{ \log t}{n} \right)^{1- \gamma} \right).
	\end{align*}
	Combine this with the lower bound on the variance in the proof of Lemma~\ref{lemma:example:SubGaussianData} (i) and note that $\omega_{(t)}^2 \leq \|\Omega\|_{op}$. This completes the proof of the first statement.
	
	Next, consider case (ii). By Markov's inequality applied to eq.~\eqref{eq:lemma:example:ApproxSparseCovariance-2} and the arguments of the proof of case (ii) in Lemma~\ref{lemma:example:LogConcaveData}, we have
	\begin{align*}
		\|\widehat{\Omega}_n^{\mathrm{naive}} - \Omega\|_{1 \rightarrow \infty} = O_p\left( \omega_{(t)} \left( \left( \frac{ \log t}{n} \right)^{1/4} \vee \sqrt{\frac{ \log t}{n}} \right) \right).
	\end{align*}
	This inequality and eq.~\eqref{eq:lemma:example:ApproxSparseCovariance-1} imply that
	\begin{align*}
		\left\|  \mathcal{T}_{\lambda}^+(\widehat{\Omega}_n^{\mathrm{naive}}) -  \Omega\right\|_{op} = O_p \left(  R_\gamma \omega_{(t)}^{1-\gamma} \left(\left( \frac{ \log t}{n} \right)^{1/4} \vee \sqrt{\frac{ \log t}{n}} \right)^{1- \gamma} \right).
	\end{align*}
	Combine this bound with the lower bound on the variance in the proof of Lemma~\ref{lemma:example:LogConcaveData} (i). This proves of the second statement.
	
	Lastly, consider case (iii). Follows ass case (ii) form Markov's inequality applied to eq.~\eqref{eq:lemma:example:ApproxSparseCovariance-2} and the arguments of the proof of case (ii) in Lemma~\ref{lemma:example:HeavyTailedData}. We skip the repetitive details.	
\end{proof}

\subsection{Proofs of Section~\ref{subsec:SphericalBootstrap}}

\begin{proof}[\textbf{Proof of Theorem~\ref{theorem:SphericalBootstrapTest}}]
	We only show the proof for $\mu \in \mathcal{H}_0$. The cases $\mu \in \mathcal{H}_0^c$ can be proved in the same way by adding the additional steps from the respective proofs of Theorem~\ref{theorem:ConsistencyLocalAlternatives} and Propositions~\ref{theorem:ConsistencySizeAlphaTestL1LInfty} and~\ref{theorem:ConsistencySizeAlphaTestEllipticallyDistributed}. Let $q \geq 1$ be the conjugate exponent to $p$, i.e. $1/p+ 1/q = 1$ and $Z \sim N(0, \Omega)$. By Theorem~\ref{theorem:Bootstrap-Lp-Norm-Quantiles},
	\begin{align*}
		&\sup_{\alpha \in (0,1)}\left| \mathrm{E}_\mu\left[ 	\varphi_{\alpha}^s(T_{n,p},\widehat{\Gamma}_n)\right] - \alpha \right| \\
		&\quad =\sup_{\alpha \in (0,1)} \Big|\mathrm{P}_\mu\left(T_{n,p} + c_p^*(1-\alpha; \widehat{\Omega}_n) - c_p^s(1-\alpha; \widehat{\Gamma}_n) \leq c_p^*(1-\alpha; \widehat{\Omega}_n) \right)  - \alpha \Big| \\
		&\quad{}\lesssim \frac{(\mathrm{E}_\mu[\|RX- \mu\|_p^3])^{1/3}}{n^{1/6}\sqrt{\mathrm{Var}(\|Z\|_p)}} + \frac{\mathrm{E}_\mu \left[ \|RX-\mu\|_p^3 \mathbf{1}\{\|RX-\mu\|_p^3 > n\: \mathrm{E}_\mu[\|RX-\mu\|_p^3]\}\right]}{\mathrm{E}_\mu\left[ \|RX-\mu\|_p^3\right]}+  \frac{ \mathrm{E}[\|Z\|_p]}{\sqrt{n\mathrm{Var}(\|Z\|_p)}}\\
		&\quad{}\quad{}\quad{}+ \inf_{\delta > 0}\left\{ \left(\frac{\delta }{ \mathrm{Var}( \|Z\|_p)} \right)^{1/3}  + \mathrm{P}_\mu\left( \|\widehat{\Omega}_n- \Omega\|_{q \rightarrow p} > \delta\right)\right\} \\
		&\quad{}\quad{} \quad{} \quad{}+ 
		\sup_{\alpha \in (0,1)}\inf_{\eta > 0} \left\{ \frac{\eta }{ \sqrt{\mathrm{Var}( \|Z\|_p)} } + \mathrm{P}_\mu\left(\left|c_p^*(\alpha; \widehat{\Omega}_n) - c_p^s(\alpha; \widehat{\Gamma}_n)\right| > \eta \right)\right\}\\
		&\quad \equiv \mathbf{I} + \mathbf{II} + \mathbf{III}.
	\end{align*}
	Under the assumptions of Theorems~\ref{theorem:SizeAlphaTest} and~\ref{theorem:ConsistencyLocalAlternatives} and Propositions~\ref{theorem:ConfidenceSets},~\ref{theorem:SizeAlphaTestL1LInfty}--\ref{theorem:ApproximateSampleAverages} the terms $\mathbf{I}$ and $\mathbf{II}$ are negligible. We therefore only need to analyze the third term $\mathbf{III}$. By Theorem 8.1 in~\cite{major1978on} we have
	\begin{align}\label{eq:theorem:SphericalBootstrapTest-1}
		\sup_{\alpha \in (0,1)}\left|c_p^*(\alpha; \widehat{\Omega}_n) - c_p^s(\alpha; \widehat{\Gamma}_n)\right| &\leq \int_0^1 \left|c_p^*(u; \widehat{\Omega}_n) - c_p^s(u; \widehat{\Gamma}_n)\right| du \nonumber\\
		&= \inf \mathrm{E}\left[ \left| \sqrt{s} \|\widehat{\Gamma}_n U^{(s)}\|_p - \|\widehat{\Gamma}_n Z^{(s)}\|_p \right| \mid X_1, \ldots, X_n  \right],
	\end{align}
	where the infimum is taken over all couplings $(U^{(s)}, Z^{(s)})$ with marginals $U^{(s)} \sim \mathrm{Unif}(\mathbb{S}^{s-1})$ and $Z^{(s)} \sim N(0, I_s)$. Recall the stochastic representation $U^{(s)} \overset{d}{=} Z^{(s)}/ \|Z^{(s)}\|_2$. Hence, we can upper bound eq.~\eqref{eq:theorem:SphericalBootstrapTest-1} by
	\begin{align}\label{eq:theorem:SphericalBootstrapTest-2}
		\mathrm{E}\left[ \left| \sqrt{s} \|\widehat{\Gamma}_n U^{(s)}\|_p - \|\widehat{\Gamma}_n Z^{(s)}\|_p \right| \mid X_1, \ldots, X_n  \right] &\leq \mathrm{E}\left[ \left\| \sqrt{s} \widehat{\Gamma}_n U^{(s)} - \widehat{\Gamma}_n Z^{(s)} \right\|_p \mid X_1, \ldots, X_n  \right]  \nonumber\\
		&\leq \|\widehat{\Gamma}_n \|_{2 \rightarrow p}  \mathrm{E}\left[ \left\| \left( \frac{\sqrt{s}}{\|Z^{(s)}\|_2}- 1\right)Z^{(s)}\right\|_2 \right] \nonumber\\
		&\leq \|\widehat{\Gamma}_n \|_{2 \rightarrow p}  \left(\mathrm{E}\left[ 2s - \sqrt{s} \|Z^{(s)}\|_2 \right] \right)^{1/2} \nonumber\\
		&\leq \frac{ \|\widehat{\Gamma}_n \|_{2 \rightarrow p}}{2s^{1/4}},
	\end{align}
	where the last inequality follows by Sterling's approximation to the Gamma function (recall that $\mathrm{E}[\|Z^{(s)}\|_2] = \Gamma\left(\frac{s+1}{2}\right)/\Gamma\left(\frac{s}{2}\right)$). To conclude the proof, combine eq.~\eqref{eq:theorem:SphericalBootstrapTest-1} and~\eqref{eq:theorem:SphericalBootstrapTest-2} with $\mathbf{III}$.
\end{proof}

\begin{proof}[\textbf{Proof of Theorem~\ref{theorem:BoundSphericalMonteCarloError}}]
	Repeating verbatim the proof of Theorem~\ref{theorem:BoundMonteCarloError} we obtain	
	\begin{align*}
		&\sup_{\alpha \in (0,1)} \sup_{\mu} \left| \mathrm{E}_\mu\left[ \varphi_{\alpha,B}^s(T_{n,p},\widehat{\Gamma}_n) \right] -\mathrm{E}_\mu\left[ \varphi_{\alpha}^s(T_{n,p},\widehat{\Gamma}_n) \right] \right|\\
		&\quad{}\quad{} \lesssim  B^{-1/2} (\log B) \sqrt{  M_{n,p}^2 \mathrm{E}\big[ \mathrm{Var}(\|\widehat{\Gamma}_nU^{(s)}\|_p \mid X_1, \ldots, X_n)\big]},
	\end{align*}
	where $\lesssim$ hides an absolute constant independent of $p, n, d, t, s, B$, and the distribution of the $X_i$'s. The theorem now follows from Lemma~\ref{lemma:BoundsVarianceUniformRV}. 
\end{proof}
\newpage
\section{Proofs of results in the supplementary materials}\label{sec:ProofsSupplement}

\subsection{Proofs of Section~\ref{subsec:GeneralResults}}
\begin{proof}[\textbf{Proof of Theorem~\ref{theorem:CLT-Lp-Norm}}]
	Trivial. Let $1 \leq p,q \leq \infty$ be conjugate exponents such that $1/p + 1/q = 1$. Apply Lemma~\ref{lemma:CLT-DimFree} with  $\mathcal{F} = \{ x \mapsto u'x :  \|u\|_q = 1, \: u \in \mathbb{R}^d\}$ and $ \{G_P(f) : f \in \mathcal{F}_n\} \equiv \{Z'u : \|u\|_q = 1, \: u \in \mathbb{R}^d\}$, where $Z \sim N(0, \Sigma)$. By Lemma~\ref{lemma:ModulusContinuityRn}, $r_n = 0$. Set $M = n^{1/3} (\mathrm{E}[\|X \|_p^3])^{1/3}$. This completes the proof.
\end{proof}
\begin{proof}[\textbf{Proof of Theorem~\ref{theorem:Bootstrap-Lp-Norm}}]
	Trivial. Apply Lemma~\ref{lemma:Bootstrap-DimFree} with $M = n^{1/3} (\mathrm{E}[\|X \|_p^3])^{1/3}$ and $Y_n = Z_{Q_n}$ (as defined in eq.~\eqref{eq:subsec:ResultsEPT-4}) and Lemma~\ref{lemma:ModulusContinuityRn}. See also the proof of Theorem~\ref{theorem:CLT-Lp-Norm} for a few more details.
\end{proof}

\begin{proof}[\textbf{Proof of Theorem~\ref{theorem:Bootstrap-Lp-Norm-Quantiles}}]
	Trivial. Apply Lemma~\ref{lemma:Bootstrap-DimFree-AsympSize} with $Y_n = Z_{Q_n}$ (as defined in eq.~\eqref{eq:subsec:ResultsEPT-4}) and Lemma~\ref{lemma:ModulusContinuityRn}.
\end{proof}

\subsection{Proofs of Section~\ref{subsec:LowerBoundVariance}}

\begin{proof}[\textbf{Proof of Theorem~\ref{theorem:LowerBoundVariance-LpNorm}}]
	First, if $p \in [1,2]$, then $2-2/p \leq 1$. Therefore, by Jensen's inequality,
	\begin{align*}
		\mathrm{E}[\|Z\|_p^{2(p-1)}] \leq \mathrm{E}\left[\|Z\|_p^p\right]^{2-2/p} = \left(\sum_{k=1}^d\sigma_k^p\right)^{2-2/p} \left(\frac{2^{p/2}}{\sqrt{\pi}} \Gamma\left(\frac{p+1}{2}\right)  \right)^{2-2/p}.
	\end{align*}
	Thus, Lemma~\ref{lemma:LowerBoundVariance-LpNorm-2}, Cauchy-Schwarz, and Sterling's approximation of the Gamma function yield
	\begin{align*}
		\mathrm{Var}(\|Z\|_p) &\geq \frac{\pi}{p^2} \left(\frac{\Gamma\left(p/2 + 1\right)}{\Gamma\left(p/2+1/2\right)}\right)^2 \left(\frac{2^{p/2}}{\sqrt{\pi}} \Gamma\left(\frac{p+1}{2}\right)  \right)^{2/p} \frac{\sum_{k=1}^d\sigma_k^{2p}}{\left(\sum_{k=1}^d\sigma_k^p\right)^{2-2/p}} \\
		&\geq \frac{\pi^{1-1/p}}{p} \left(\Gamma\left(\frac{p+1}{2}\right)  \right)^{2/p} \left( \frac{1}{d} \sum_{k=1}^d\sigma_k^p\right)^{2/p} d^{2/p-1}\\
		&\geq \frac{\pi}{6}\left( \frac{1}{d} \sum_{k=1}^d \sigma_k^p\right)^{2/p} p^{-1/2} d^{2/p-1}.
	\end{align*}
	Second, if $p \in (2, 2 \log d)$, then $2p-2 \geq p$. Therefore, 
	\begin{align*}
		\mathrm{E}[\|Z\|_p^{2(p-1)}] \leq d^{1-2/p}\mathrm{E}\left[\|Z\|_{2p-2}^{2p-2}\right] =  d^{1-2/p} \left(\sum_{k=1}^d\sigma_k^{2p-2} \right) \left(\frac{2^{p-1}}{\sqrt{\pi}} \Gamma\left(\frac{2p-1}{2}\right)  \right).
	\end{align*}
	Thus, by Lemma~\ref{lemma:LowerBoundVariance-LpNorm-2},
	\begin{align*}
		\mathrm{Var}(\|Z\|_p) &\geq \frac{\pi^{1/2}/3}{2^{2p-1}} \left(\frac{\Gamma\left(p/2 + 1\right)}{\Gamma\left(p-1/2\right)}\right)^2 \Gamma\left(\frac{2p-1}{2}\right) \left(\frac{\sum_{k=1}^d\sigma_k^{2p}}{\sum_{k=1}^d\sigma_i^{2p-2}}\right) d^{2/p-1} \\
		& \geq \frac{\pi}{6}   \left( \frac{1}{d} \sum_{k=1}^d \sigma_k^{2p}\right)^{1/p} 2^{-3p} p^2 d^{2/p-1}.
	\end{align*}
	Third, if $p \in [2 \log d, \infty)$, then $\mathrm{E}[\|Z\|_p] \leq e \mathrm{E}[\|Z\|_\infty] \leq e \sqrt{2 \log d}$. We combine this with Lemma~\ref{lemma:LowerBoundVariance-LpNorm-6} and simplify the expression.
	\noindent\\
	Lastly, if $p = \infty$, we directly apply Lemma~\ref{lemma:BoundsVariance-SeparableProcess} with the explicit constants given in the proof and $\mathrm{E}[\|Z\|_\infty] \leq e \sqrt{2 \log d}$.
\end{proof}

\begin{proof}[\textbf{Proof of Theorem~\ref{theorem:LowerBoundVariance-LpNorm-Refinements}}]
	Case $p=1$ follows from Lemma~\ref{lemma:LowerBoundVariance-LpNorm-2} and case $p = \infty$ is a restatement of Lemma~\ref{lemma:LowerBoundVariance-LpNorm-3}. We only need to proof case $p=2$. A naive application of Lemma~\ref{lemma:LowerBoundVariance-LpNorm-2} with for $p=2$ gives
	\begin{align*}
		\mathrm{Var}(\|Z\|_2) \geq \frac{\sum_{k=1}^d \sigma_k^4}{\sum_{k=1}^d \sigma_i^2}.
	\end{align*}
	Further, recall that $\|Z\|_2 \overset{d}{=}\|\Lambda^{1/2}G\|_2$, where $G$ is a standard Gaussian random vector and $\Lambda$ is a diagonal matrix with the eigenvalues of $\Sigma$ on its diagonal. Thus, by Lemma~\ref{lemma:LowerBoundVariance-LpNorm-2} we also have
	\begin{align*}
		\mathrm{Var}(\|Z\|_2) = \mathrm{Var}(\|\Lambda^{1/2}G\|_2) \geq \frac{\sum_{k=1}^d \lambda_k^2}{\sum_{k=1}^d \lambda_k} = \frac{\mathrm{tr}(\Sigma^2)}{\mathrm{tr}(\Sigma)}.
	\end{align*}
\end{proof}
\begin{proof}[\textbf{Proof of Lemma~\ref{lemma:LowerBoundVariance-LpNorm-2}}]
	
	Recall the following numerical result, which is an immediate consequence of the Hadamard-Hermite inequality~\citep[][Theorem 9]{bullen1988means} applied to the convex function $f(t) = t^{\theta -1} $, $t, \theta \geq 0$: For $a, b, \theta > 0$,
	\begin{align*}
		2|a^\theta - b^\theta| \leq \theta |a-b| \left(a^{\theta-1} + b^{\theta-1}\right).
	\end{align*}
	Let $Y$ be an independent copy of $X$. By Cauchy-Schwarz and above inequality with $a = \|X\|_p$, $b = \|Y\|_p$ and $\theta = p$ we have
	\begin{align*}
		\mathrm{E}\left[(\|X\|_p - \|Y\|_p)^2\right]^{1/2} \mathrm{E}\left[(\|X\|_p^{p-1} + \|Y\|_p^{p-1})^2\right]^{1/2} \geq \frac{2}{p} \mathrm{E}\left[|\|X\|_p^p - \|Y\|_p^p|\right],
	\end{align*}
	and, thus,
	\begin{align}\label{eq:lemma:LowerBoundVariance-LpNorm-2-1}
		\mathrm{Var}(\|X\|_p) \geq \frac{2}{p^2} \frac{\mathrm{E}\left[|\|X\|_p^p - \|Y\|_p^p|\right]^2}{\mathrm{E}\left[(\|X\|_p^{p-1} + \|Y\|_p^{p-1})^2\right]}.
	\end{align}
	(See also eq. (3.4) in~\cite{paouris2018Dvoretzky}.) Since $X$ and $Y$ are identically distributed, the denominator in above expression can be upper bounded by $4 \mathrm{E}[\|X\|_p^{2(p-1)}]$. The remainder of the proof is concerned with deriving a lower bound on the numerator.
	
	Let $\varepsilon_1, \ldots, \varepsilon_d$ be i.i.d. Rademacher random variables independent of $X$ and $Y$. Then, $|X_k|^p - |Y_k|^p \overset{d}{=} \varepsilon_k \big||X_k|^p - |Y_k|^p \big|$ all $1 \leq k \leq d$, and Khintchine's inequality for (conditional) Rademacher averages yields
	\begin{align}\label{eq:lemma:LowerBoundVariance-LpNorm-2-2}
		\mathrm{E}\left[|\|X\|_p^p - \|Y\|_p^p|\right] = \mathrm{E}\left[\left|\sum_{k=1}^d\varepsilon_k\big||X_k|^p - |Y_k|^p \big| \right|\right] \geq \frac{1}{\sqrt{2}} \mathrm{E}\left[\left(\sum_{k=1}^d 
		(|X_k|^p - |Y_k|^p)^2 \right)^{1/2}\right].
	\end{align}
	Next, recall Minkowski's integral inequality~\citep[][Theorem 202]{hardy1988inequalities}: Let $(U, \mu)$ and $(V, \nu)$ be measure spaces and $f : U \times V \rightarrow \mathbb{R}$ be a measurable map. Then, for $r \geq 1$ arbitrary,
	\begin{align*}
		\left(\int_U \left| \int_V f(u, v) d\nu(v) \right|^r d \mu(u)\right)^{1/r} \leq  \int_V \left(\int_U |f(u,v)|^rd\mu(u)\right)^{1/r} d \nu(v).
	\end{align*}
	For $\nu = \gamma \times \gamma$, where $\gamma$ is the Gaussian measure associated with random vector $X$, $V = \mathbb{R}^n \times \mathbb{R}^n$, $\mu$ the counting measure over the set $U = \{1, \ldots, n\}$, $f(u, v_1, \ldots, v_{2n}) = \big||v_u|^p - |v_{n + u}|^p\big|$, and $r = 2$, this inequality implies
	\begin{align}\label{eq:lemma:LowerBoundVariance-LpNorm-2-3}
		\mathrm{E}\left[\left(\sum_{k=1}^d (|X_k|^p - |Y_k|^p)^2 \right)^{1/2}\right] \geq \left(\sum_{k=1}^d\mathrm{E}\left[\big||X_k|^p - |Y_k|^p\big| \right]^2\right)^{1/2}.
	\end{align}
	Observe that the expression on the right hand side in above inequality depends only on the marginals $X_k, Y_k \sim N(0, \mathrm{Var}(X_k))$, $1 \leq k \leq n$. Therefore,
	\begin{align}\label{eq:lemma:LowerBoundVariance-LpNorm-2-4}
		\left(\sum_{k=1}^d \mathrm{E}\left[\big||X_k|^p - |Y_k|^p\big| \right]^2\right)^{1/2} = \mathrm{E}\left[\big||Z|^p - |Z'|^p\big| \right] \left(\sum_{k=1}^d\mathrm{Var}(X_k)^p \right)^{1/2},
	\end{align}
	where $Z, Z'$ are independent standard normal random variables. To lower bound the first factor on the right hand side in above inequality, we compute
	\begin{align*}
		\mathrm{E}\left[\big||Z|^p - |Z'|^p\big| \right] &= \int_0^\infty \int_{-\pi}^\pi \big||r \cos \varphi|^p - |r \sin\varphi|^p\big| \frac{re^{-r^2/2}}{2 \pi} dr d\varphi\\
		& = \frac{2^{p/2}}{2\pi}\Gamma\left(\frac{p}{2}+ 1\right) 8 \int_0^{\pi/4} (\cos \varphi)^p - (\sin\varphi)^p d \varphi.
	\end{align*}
	Observe that the integrand $I_p(\varphi) = (\cos \varphi)^p - (\sin\varphi)^p$ satisfies $I_p(0) = 1$, $I_p(\pi/4) = 0$, $I_p(\varphi) \geq 0$, and
	\begin{align*}
		I_p'(\varphi) &= - p (\sin \varphi) (\cos\varphi) \left[ (\sin \varphi)^{p-2} + (\cos \varphi)^{p-2}\right]\\
		I_p''(\varphi) &= - p(p-1)\left[(\cos \varphi)^2(\sin\varphi)^{p-2} - (\sin \varphi)^2(\cos \varphi)^{p-2}\right] - p \left[ (\cos \varphi)^p - (\sin \varphi)^p \right].
	\end{align*}
	Moreover, $I_p'(\varphi) \leq 0$ for all $p \in [1,\infty)$ and all $\varphi \in [0, \pi/4]$, and $I_p''(\varphi) \leq 0$ for all $p \in [1, 2]$ and all $\varphi \in [0, \pi/4]$. Thus, in the regime $p \in [1,2]$ the integrand is nonincreasing and concave on the interval $[0, \pi/4]$. Therefore, its graph over $[0,\pi/4]$ lies on or above the line connecting $(0,1)$ and $(\pi/4, 0)$, and
	\begin{align*}
		\frac{\pi}{4} \geq \int_0^{\pi/4} (\cos \varphi)^p - (\sin\varphi)^p d \varphi \geq \int_0^{\pi/4} \left(1 - \frac{4}{\pi}x\right) dx = \frac{\pi}{8}.
	\end{align*}
	In the regime $p \in (2, \infty)$, $I_p$ is nonincreasing on $[0,\pi/4]$. We can therefore fit a rectangle under the graph with lower left corner $(0,0)$ and upper right corner $(a, I_p(a))$ for any $a \in [0,\pi/4]$. For $a = \pi/6$ we obtain the simple expression
	\begin{align*}
		\int_0^{\pi/4} (\cos \varphi)^p - (\sin\varphi)^p d \varphi \geq \frac{\pi}{6} \left( \left(\frac{3}{4}\right)^{p/2} - \left(\frac{1}{4}\right)^{p/2}\right) \overset{(a)}{\geq} \frac{\pi}{6} \frac{p}{2} \left(\frac{1}{4}\right)^{p/2-1}\left(\frac{3}{4} - \frac{1}{2}\right) = \frac{\pi}{6} \frac{p}{2} 2^{-p} ,
	\end{align*}
	where (a) follows from a first-order Taylor approximation of the convex (because $p > 2$) map $x \mapsto x^{p/2}$.
	
	In summary, for $p \in [1,2]$,
	\begin{align}\label{eq:lemma:LowerBoundVariance-LpNorm-2-5}
		\mathrm{E}\left[\big||Z|^p - |Z'|^p\big| \right] \geq \frac{2^{p/2}}{2}\Gamma\left(\frac{p}{2}+ 1\right), 
	\end{align}
	while for $p \in (2, \infty)$,
	\begin{align*}
		\mathrm{E}\left[\big||Z|^p - |Z'|^p\big| \right] \geq \frac{1}{ 2^{p/2}} \frac{p}{3} \Gamma\left(\frac{p}{2}+ 1\right).
	\end{align*}
	
	Combine eq.~\eqref{eq:lemma:LowerBoundVariance-LpNorm-2-1}--\eqref{eq:lemma:LowerBoundVariance-LpNorm-2-5} and conclude that
	\begin{align*}
		\mathrm{Var}(\|Z\|_p) \geq \frac{\sum_{k=1}^d\mathrm{Var}(Z_k)^p}{\mathrm{E}[\|Z\|_p^{2(p-1)}]} \times 
		\begin{cases}
			\frac{2^p}{p^2}\Gamma^2\left(\frac{p}{2} + 1\right), \quad{} & \mathrm{if} \quad{} p \in [1,2]\\
			\frac{2^{-p}}{3} \Gamma^2\left(\frac{p}{2} + 1\right), \quad{} & \mathrm{if} \quad{} p \in (2, \infty).
		\end{cases}
	\end{align*}
\end{proof}
\begin{proof}[\textbf{Proof of Lemma~\ref{lemma:LowerBoundVariance-LpNorm-6}}]
	We model our proof after the proof of Theorem 1.8 in~\cite{ding2015multiple}. The necessary modifications are relatively straightforward; however, since the original proof of Theorem 1.8 is very condensed we fill in many details. First, consider the case $\mathrm{E}[\|Z\|_p]/\sigma_{(1)} \leq 1/2$. By Chebyshev's inequality and since $\|x\|_\infty \leq \|x\|_p$ for all $x \in \mathbb{R}^n$,
	\begin{align}\label{eq:LowerBoundVariance-LpNorm-6-1}
		\mathrm{Var}\left(\|Z\|_p/\sigma_{(1)} \right) &\geq \mathrm{P}(\|Z\|_p/\sigma_{(1)} \geq 1) \left(1 - \mathrm{E}[\|Z\|_p]/\sigma_{(1)} \right)^2 \nonumber\\
		&\geq  \frac{1}{4}\mathrm{P}(\|Z\|_p/\sigma_{(1)}\geq 1)  \nonumber \\
		&\geq  \frac{1}{4}\mathrm{P}(\|Z\|_\infty/\sigma_{(1)}\geq 1)  \nonumber \\
		&\geq  \frac{1}{4}\min_{1 \leq k \leq n}\mathrm{P}(Z_k/\sigma_{(1)} \geq 1).
	\end{align}
	By the lower bound on Mill's ratio for a standard normal random variable,
	\begin{align}\label{eq:LowerBoundVariance-LpNorm-6-2}
		\min_{1 \leq k \leq n}\mathrm{P}(Z_k/\sigma_{(1)} \geq 1) \geq \min_{1 \leq k \leq n}\mathrm{P}(Z_k/\sqrt{\mathrm{Var}(Z_k)} \geq 1) \geq \frac{1}{2} \frac{e^{-1/2}}{\sqrt{2\pi}} \geq \frac{1}{9}.
	\end{align}
	Combine eq.~\eqref{eq:LowerBoundVariance-LpNorm-6-1} and~\eqref{eq:LowerBoundVariance-LpNorm-6-2} and conclude that
	\begin{align}\label{eq:LowerBoundVariance-LpNorm-6-3}
		\mathrm{Var}\left(\|Z\|_p/\sigma_{(1)}\right) \left(1 + \mathrm{E}[\|Z\|_p/\sigma_{(1)}]\right)^2 \geq \mathrm{Var}\left(\|Z\|_p/\sigma_{(1)}\right) \geq \frac{1}{36}.
	\end{align}
	
	Next, consider the case $\mathrm{E}[\|Z\|_p]/\sigma_{(1)} > 1/2$. Let $t = \sqrt{1 - \left(2 \mathrm{E}[\|Z\|_p/\sigma_{(1)}]\right)^{-2}}$, $U^t_\infty = \{u \in \mathbb{R}^d: X'u \geq t \mathrm{E}[\|Z\|_\infty]\}$, and $B_q^d(1) = \{ u \in \mathbb{R}^d: \|u\|_q =1\}$ for $1/ p + 1/q = 1$. Note that 
	\begin{align*}
		\sqrt{1 - t^2}\mathrm{E}[\|Z\|_p/\sigma_{(1)}] + \frac{\left(4 \mathrm{E}[\|Z\|_p/\sigma_{(1)}]\right)^{-1}}{\sqrt{1 - t^2}} = 1.
	\end{align*}
	Thus, by Lemma~\ref{lemma:ding2015multiple},
	\begin{align}\label{eq:LowerBoundVariance-LpNorm-6-4}
		\left(4 \mathrm{E}[\|Z\|_p]/\sigma_{(1)}\right)^2 \mathrm{Var}\left(\|Z\|_p/\sigma_{(1)}\right) \geq \mathrm{P} \left( \sup_{u \in U^t_\infty \cap B^d_q(1)} |Y'u|/\sigma_{(1)} \geq 1 \right), 
	\end{align}
	where $Y$ is an independent copy of $Z$. Denote by $\mathcal{E}^d_1$ the set of extreme points of the cross-polytope $B_1^d(1)$. Note that $B^d_q(1) \cap \mathcal{E}^d_1 = \mathcal{E}^d_1$ for all $ q \geq 1$. To lower bound the probability on the right hand side in above display, we compute 
	\begin{align*}
		\mathrm{P} \left( \sup_{u \in U^t_\infty \cap B^d_q(1) } |Y'u|/\sigma_{(1)} \geq 1 \right) &\geq \mathrm{P} \left( \sup_{u \in U^t_\infty \cap \mathcal{E}^d_1} |Y'u|/\sigma_{(1)} \geq 1 \mid U^t_\infty \cap \mathcal{E}^d_1 \neq \varnothing \right) \mathrm{P}\left( U^t_\infty \cap \mathcal{E}^d_1 \neq \varnothing \right) \\
		& \geq \inf_{u \in \mathcal{E}^d_1} \mathrm{P} \left( Y'u/\sigma_{(1)} \geq 1 \mid U^t_\infty \cap \mathcal{E}^d_1 \neq \varnothing \right) \mathrm{P}\left( U^t_\infty \cap \mathcal{E}^d_1 \neq \varnothing \right)\\
		& = \min_{1 \leq k \leq d} \mathrm{P} \left( Y_k/\sigma_{(1)} \geq 1 \mid U^t_\infty \cap \mathcal{E}^d_1 \neq \varnothing \right) \mathrm{P}\left( U^t_\infty \cap \mathcal{E}^d_1 \neq \varnothing \right).
	\end{align*}
	Since $Y \ind Z$ and $\{ U^t_\infty \cap \mathcal{E}^d_1 \neq \varnothing \}$ depends on $Z$ only, we have, as in eq.~\eqref{eq:LowerBoundVariance-LpNorm-6-2},
	\begin{align}\label{eq:LowerBoundVariance-LpNorm-6-5} 
		\min_{1 \leq k \leq d} \mathrm{P}  \left( Y_k/\sigma_{(1)} \geq 1 \mid U^d_p \cap \mathcal{E}^d_1 \neq \varnothing \right) = \min_{1 \leq k \leq d} \mathrm{P}  \left( Y_k/\sigma_{(1)} \geq 1\right) \geq \frac{1}{9}.
	\end{align}
	Next, recall the following Paley-Zygmund-type lower bound on the tail probability of a non-negative random variable $V \geq 0$: For all $t \in (0,1)$,
	\begin{align*}
		\mathrm{P}(V > t \mathrm{E}[V]) \geq \frac{(1-t)^2(\mathrm{E}[V])^2}{\mathrm{Var}(V) + (1-t)^2(\mathrm{E}[V])^2}.
	\end{align*}
	Since $(1 - \sqrt{1-t^2})^2 \geq t^4/4$ for all $t\in [0,1]$, we have, by above inequality,
	\begin{align}\label{eq:LowerBoundVariance-LpNorm-6-6}
		\mathrm{P}\left(U^t_\infty \cap \mathcal{E}^d_1 \neq \varnothing \right) &= \mathrm{P} \left( \|Z\|_\infty \geq t \mathrm{E}[\|Z\|_\infty] \right) \nonumber\\
		&\geq \frac{(\mathrm{E}[\|Z\|_\infty/\sigma_{(d)}])^2}{4 (\mathrm{E}[\|Z\|_p/\sigma_{(d)}])^4\mathrm{Var}(\|Z\|_\infty/\sigma_{(d)}) + (\mathrm{E}[\|Z\|_\infty/\sigma_{(d)}])^2} \nonumber\\
		&\geq \frac{1}{4 d^{4/p} (\mathrm{E}[\|Z\|_\infty/\sigma_{(d)}])^2\mathrm{Var}(\|Z\|_\infty/\sigma_{(d)}) + 1}.
	\end{align}
	Since $a > 3b$ implies $a/(a-b) \leq 3/2$ for all $a, b > 0$, Lemma~\ref{lemma:BoundsVariance-SeparableProcess} yields
	\begin{align*}
		&(\mathrm{E}[\|Z\|_\infty/\sigma_{(d)}])^2\mathrm{Var}(\|Z\|_\infty/\sigma_{(d)}) \nonumber\\
		&\quad{}\quad{}\leq
		\begin{cases}
			108 (\log 2), & \mathrm{if}\:\:\mathrm{E}[\|Z\|_\infty/ \sigma_{(d)}] \leq 3^2 \sqrt{\log 2}\\
			(3/2)^2 60^2 +  11^2 (\mathrm{E}[\|Z\|_\infty/ \sigma_{(d)}])^2\rho , & \mathrm{if}\:\: \mathrm{E}[\|Z\|_\infty/ \sigma_{(d)}] \geq 3^2 \sqrt{\log 2}.
		\end{cases}
	\end{align*}
	Thus, for $p \geq 2 \log d$, we can lower bound the expression in eq.~\eqref{eq:LowerBoundVariance-LpNorm-6-6} by
	\begin{align}\label{eq:LowerBoundVariance-LpNorm-6-7}
		\frac{1}{180^2 e^2}\frac{1}{ 2+ \rho (\mathrm{E}[\|Z\|_\infty/ \sigma_{(d)}])^2}.
	\end{align}
	Combine eq.~\eqref{eq:LowerBoundVariance-LpNorm-6-4}--\eqref{eq:LowerBoundVariance-LpNorm-6-7} to conclude that
	\begin{align}\label{eq:LowerBoundVariance-LpNorm-6-8}
		e^2 48^4 (\mathrm{E}[\|Z\|_p/\sigma_{(1)}])^2 \mathrm{Var}\left(\|Z\|_p/\sigma_{(1)}\right) \geq \frac{1}{2+ \rho (\mathrm{E}[\|Z\|_\infty/ \sigma_{(d)}])^2}.
	\end{align}
	Now, eq.~\eqref{eq:LowerBoundVariance-LpNorm-6-3} and~\eqref{eq:LowerBoundVariance-LpNorm-6-8} imply for all $p \geq 2 \log d$,
	\begin{align}\label{eq:LowerBoundVariance-LpNorm-6-9}
		\mathrm{Var}\left(\|Z\|_p\right)  \geq \frac{1}{e^2 48^4} \left(\frac{\sigma_{(1)}^2}{\sigma_{(1)} + \mathrm{E}[\|Z\|_p]}\right)^2\frac{1}{ 2+ \rho (\mathrm{E}[\|Z\|_\infty/ \sigma_{(d)}])^2}.
	\end{align}
	The claim now follows from $\mathrm{E}[\|Z\|_\infty] \leq \sigma_{(d)} \sqrt{2 \log(2d)}$.
\end{proof}
\begin{proof}[\textbf{Proof of Lemma~\ref{lemma:LowerBoundVariance-LpNorm-3}}]
	By eq. (75) in Theorem 10 in~\cite{deng2020beyond} with $x_k = 1 + \sqrt{2 \log k}$ for $1 \leq k \leq d$, we have
	\begin{align*}
		M(Z) \leq 1/ \sigma_{(1)} + \max_{1 \leq k \leq d} \left(1 + \sqrt{2\log k}\right)/\sigma_{(k)}.
	\end{align*}
	Hence, the claim follows since $\mathrm{Var}(Z)M^2(Z) \lesssim 1$~\citep[e.g.][]{giessing2022anticoncentration} and above upper bound on $M(Z)$.
\end{proof}

\subsection{Proofs of Section~\ref{subsec:ResultsEPT}}

\begin{proof}[\textbf{Proof of Lemma~\ref{lemma:ModulusContinuityRn}}]	
	First, we consider the empirical process $\{\mathbb{G}_n(f): f \in \mathcal{F}\}$.  To fix notation, denote by $\Sigma \in \mathbb{R}^{d \times d}$ the covariance matrix of $X \in \mathbb{R}^d$. Since $\Sigma \in \mathbb{R}^{d \times d}$ is positive semi-definite, there exists $\Gamma \in \mathrm{R}^{d \times r}$ with $r = \mathrm{rank}(\Sigma)$ such that $\Gamma \Gamma' = \Sigma$. Next, let $ h \in \mathcal{F}_\delta'$ be arbitrary. By definition there exist $f, g \in \mathcal{F}$ such that for all $x \in \mathbb{R}^d$,
	\begin{align*}
		|h(x)| = |f(x) - g(x)| = |x'(u-v)| 
		& \leq  \|x\|_p \| (\Gamma'\Gamma)^{-1} \Gamma' \|_{2 \rightarrow q}\|\Gamma (u-v)\|_2\\
		& \leq  \|x\|_p \| (\Gamma'\Gamma)^{-1} \Gamma' \|_{2 \rightarrow q} \delta\\
		& = \|x\|_p \|\Gamma (\Gamma'\Gamma)^{-1} \|_{p \rightarrow 2} \delta,
	\end{align*}
	where $1 \leq p,q \leq \infty$ are conjugate exponents, i.e. $1/p + 1/q = 1$. Thus, 
	\begin{align*}
		F_\delta(x) : = \delta  \|  \Gamma (\Gamma'\Gamma)^{-1}\|_{p \rightarrow 2}\|x\|_p
	\end{align*}
	is an envelope for $\mathcal{F}_\delta'$. Whence, by Theorem 2.14.1 in~\cite{vandervaart1996weak},
	\begin{align*}
		\big\| \| \mathbb{G}_n\|_{\mathcal{F}_\delta'} \big\|_{P,1} \lesssim J(1, \mathcal{F}_\delta') (P F_\delta^2)^{1/2} \lesssim  \delta \: C_d  \:  \| \Gamma (\Gamma'\Gamma)^{-1}\|_{p \rightarrow 2}  \big(\mathrm{E}\|X\|_p^2\big)^{1/2},
	\end{align*}
	where the second inequality holds because $\mathcal{F}_\delta'$ is a VC-subgraph class and hence the uniform entropy integral $J(1,\mathcal{F}_\delta')$ can be upper bounded by a constant $C_d > 0$ independent of $\delta > 0$. Hence, we define
	\begin{align}\label{eq:lemma:ModulusContinuityRn-1}
		\phi_n(\delta) :=  \delta \: C_d  \:  \| \Gamma (\Gamma'\Gamma)^{-1}\|_{p \rightarrow 2} \left(\frac{\mathrm{E}[\|X\|_p^2]}{\mathrm{Var}(\|G_P\|_{\mathcal{F}})}\right)^{1/2}.
	\end{align}
	
	Second, observe that $\{Z'u : \|u\|_q = 1, \: u \in \mathbb{R}^d\}$ with $Z \sim N(0, \Sigma)$ is a version of the Gaussian $P$-bridge process $\{G_P(f) : f \in \mathcal{F}_n\}$. Thus, the same arguments as above yield
	\begin{align*}
		\mathrm{E} \| G_P\|_{\mathcal{F}_\delta'} \lesssim J(1, \mathcal{F}_\delta') \big(\mathrm{E}[ F_\delta^2(Z)] \big)^{1/2} \lesssim  \delta \:  C_d  \: \|\Gamma (\Gamma'\Gamma)^{-1} \|_{p \rightarrow 2}  \big(\mathrm{E}[\|Z\|_p^2]\big)^{1/2},
	\end{align*}
	and, therefore,
	\begin{align}\label{eq:lemma:ModulusContinuityRn-2}
		\psi_n(\delta) := \delta \: C_d  \:  \| \Gamma (\Gamma'\Gamma)^{-1}\|_{p \rightarrow 2} \left(\frac{\mathrm{E}[\|Z\|_p^2]}{\mathrm{Var}(\|G_P\|_{\mathcal{F}})}\right)^{1/2}.
	\end{align}
	
	Lastly, consider $\{Z_n'u : \|u\|_q = 1, \: u \in \mathbb{R}^d\}$ with $Z_n \sim N(0, \widehat{\Sigma}_n)$, where $\widehat{\Sigma}_n$ is a positive semidefinite estimate of $\Sigma$ based on the $X_i$'s only. As discussed in Section~\ref{subsec:ResultsEPT} this is a representation of the Gaussian $Q_n$-motion $\{Z_{Q_n}(f) : f \in \mathcal{F}_n \}$ with covariance function $(f, g) \mapsto f'\widehat{\Sigma}_n g$. Moreover, there exists $\widehat{\Gamma}_n \in \mathrm{R}^{d \times r}$ with $r = \mathrm{rank}(\widehat{\Sigma}_n)$ such that $\widehat{\Gamma}_n \widehat{\Gamma}_n' = \widehat{\Sigma}_n$. 
	Hence, repeating the arguments from above, we conclude that
	\begin{align*}
		\mathrm{E} \big[ \| Z_{Q_n}\|_{\mathcal{F}_\delta'} \mid X_1, \ldots, X_n \big] &\lesssim J(1, \mathcal{F}_\delta') \big(\mathrm{E}[ F_\delta^2(Z_n) \mid X_1, \ldots, X_n] \big)^{1/2} \\
		&\lesssim  \delta \:  C_d  \: \|\widehat{\Gamma}_n (\widehat{\Gamma}_n'\widehat{\Gamma}_n)^{-1}\|_{p \rightarrow 2}  \big(\mathrm{E}[\|Z_n\|_p^2 \mid X_1, \ldots, X_n]\big)^{1/2},
	\end{align*}
	and, thus,
	\begin{align}\label{eq:lemma:ModulusContinuityRn-3}
		\upsilon_n(\delta) := \delta \: C_d  \:  \|\widehat{\Gamma}_n (\widehat{\Gamma}_n'\widehat{\Gamma}_n)^{-1}\|_{p \rightarrow 2} \left(\frac{\mathrm{E}[\|Z_n\|_p^2]}{\mathrm{Var}(\|G_P\|_{\mathcal{F}})}\right)^{1/2}.
	\end{align}
	
	To conclude the proof, combine eq.~\eqref{eq:lemma:ModulusContinuityRn-1}--\eqref{eq:lemma:ModulusContinuityRn-3} and deduce that $\inf\big\{ \upsilon_n(\delta) \vee \psi_n(\delta) \vee \phi_n(\delta): \delta > 0\big\} = 0$.
\end{proof}

\subsection{Proofs of Section~\ref{subsec:Miscellaneous}}
\begin{proof}[\textbf{Proof of Lemma~\ref{lemma:EVNormGaussian}}]
	We have
	\begin{align}\label{eq:lemma:EVNormGaussian-1}
		\mathrm{E}[\|Z\|_p^p]^{1/p} / \sqrt{8\pi p} \overset{(a)}{\leq} \mathrm{E}[\|Z\|_p] \overset{(b)}{\leq} \mathrm{E}[\|Z\|_p^p]^{1/p},
	\end{align}
	where (a) follows from the reverse Liapunov inequality for moments of surpema of Gaussian processes~\citep[e.g.][Proposition A.2.4]{vandervaart1996weak} and (b) follows from Jensen's inequality. The explicit constant $1/\sqrt{8\pi p}$ is most likely not optimal; we have obtained it from the proof of Corollary 3.2 in~\cite{ledoux1991probability}. Also, for all $\ell_p$-norms with $1 \leq p \leq \infty$, we have
	\begin{align}\label{eq:lemma:EVNormGaussian-2}
		\mathrm{E}[\|Z\|_\infty] \leq \mathrm{E}[\|Z\|_p]  \leq d^{1/p} \mathrm{E}[\|Z\|_\infty] 
	\end{align}
	Straightforward computations and Sterling's formula now yield, for $1 \leq p < \infty$,
	\begin{align}\label{eq:lemma:EVNormGaussian-3}
		\mathrm{E}[\|Z\|_p] =  \left( \sum_{k=1}^d \sigma_k^p\right)^{1/p} \frac{2^{1/2}}{\pi^{1/(2p)}}\Gamma\left(\frac{p+1}{2}\right)^{1/p} \asymp \sqrt{p} \left( \sum_{k=1}^d \sigma_k^p\right)^{1/p},
	\end{align}
	while for $p = \infty$, Dudley's entropy bound and Sudakov's inequality yield,
	\begin{align}\label{eq:lemma:EVNormGaussian-4}
		\min_{1 \leq k \leq d} \sigma_k \sqrt{ \log d} \lesssim \mathrm{E}[\|Z\|_\infty] \lesssim \max_{1 \leq k \leq d} \sigma_k\sqrt{ \log d}.
	\end{align}
	To complete the proof combine eq.~\eqref{eq:lemma:EVNormGaussian-1}--\eqref{eq:lemma:EVNormGaussian-4}.
\end{proof}

\begin{proof}[\textbf{Proof of Lemma~\ref{lemma:UpperBoundQuantiles}}]
	First, we establish the upper bound for all $\alpha \in (0,1)$. Notice that the map $u \mapsto \|\Sigma^{1/2} u\|_p$ is Lipschitz continuous with respect to the Euclidean norm with Lipschitz constant $\|\Sigma^{1/2}\|_{2 \rightarrow p} := \sup_{\|u\|_2 \leq 1} \|\Sigma^{1/2}u\|_p$. Thus, by the classical Gaussian concentration inequality for Lipschitz continuous functions, for all $t > 0$,
	\begin{align*}
		\mathrm{P}\left(\|Z\|_p- \mathrm{E}[\|Z\|_p] \geq t \right) \leq \exp\left\{-\frac{t^2}{2\|\Sigma^{1/2}\|_{2 \rightarrow p}^2}\right\}.
	\end{align*}
	In particular,
	\begin{align*}
		\mathrm{P}\left(\|Z\|_p > \mathrm{E}[\|Z\|_p]  + \sqrt{2\log(1/\alpha)} \|\Sigma^{1/2}\|_{2 \rightarrow p}\right) \leq \alpha.
	\end{align*}
	By Chebyshev's inequality
	\begin{align*}
		\mathrm{P}\left(\|Z\|_p- \mathrm{E}[\|Z\|_p] \geq t \right) \leq \frac{\mathrm{Var}(\|Z\|_p)}{t^2},
	\end{align*}
	and therefore
	\begin{align*}
		\mathrm{P}\left(\|Z\|_p > \mathrm{E}[\|Z\|_p]  + \sqrt{\mathrm{Var}(\|Z\|_p)/\alpha} \right) \leq \alpha.
	\end{align*}
	Now, the upper bound follows from the definition of the quantile $c_{n,p}(1-\alpha)$.
	
	To establish the lower bound for $\alpha \in (0,1/2]$, recall the following inequality:
	\begin{align*}
		\left|\mathrm{E}[\|Z\|_p] - c_{n,p}(1/2)\right| \leq \sqrt{\mathrm{Var}(\|Z\|_p)}.
	\end{align*}
	Whence, for all $\alpha \in (0,1/2]$ it follows that
	\begin{align*}
		c_{n,p}(1-\alpha)  \geq c_{n,p}(1/2) \geq \mathrm{E}[\|Z\|_p] - \sqrt{\mathrm{Var}(\|Z\|_p)}.
	\end{align*}
	To conclude, note that by the Gaussian Poincar{\'e} inequality, $\sqrt{\mathrm{Var}(\|Z\|_p)} \leq \|\Sigma^{1/2}\|_{2 \rightarrow p}$. (Notice that this argument yields a tighter lower bound than if we had used the Gaussian concentration inequality for Lipschitz continuous functions.)
\end{proof}

\begin{proof}[\textbf{Proof of Lemma~\ref{lemma:ReverseLyapunovLogConcave}}]
	The result is known with an unspecified constant $C_{s,t} > 0$~\citep[][p. 35]{ledoux2001concentration}. In the following proof we obtain the explicit dependence of $C_{s,t}$ on $s, t$. 
	
	Without loss of generality we can assume that $\left(\mathrm{E}[\|X\|^s]\right)^{1/s} = 1$. The general case follows upon re-scaling of the norm $\|\cdot\|$. Thus, by Markov's inequality we have
	\begin{align*}
		\mathrm{P}\left(\|X\| \geq 4\right) \leq 4^{-s}.
	\end{align*}
	Define $A : = \{ x \in \mathbb{R}^d : \|x\| \leq 4\}$. Then $\mathrm{P}\left( X \in A\right) \geq 1 - 4^s > 1/2$. Hence, by Borell's inequality for log-concave measures~\citep[e.g.][Proposition 2.14]{ledoux2001concentration}, for all $r > 1$,
	\begin{align*}
		\mathrm{P}\left(\|X\|> 4r\right) = 	\mathrm{P}\left(X \in (rA)^c \right) \leq \left(\frac{4^{-s}}{1-4^{-s}}\right)^{r/2} \leq e^{-rs/2},
	\end{align*}
	since $4^s > e^s + 1$ for all $s \geq 1$. Using Tonelli-Fubini, a change of variables, and above inequality we compute
	\begin{align*}
		\mathrm{E}[\|X\|^t] 
		\leq 4^t+ t \int_4^\infty r^{t-1} \mathrm{P}\left(\|X\| \geq r \right) dr
		\leq 4^t + 4^t t \int_1^\infty r^{t-1} e^{-2rs} dr
		\leq 4^t + \left(\frac{2}{s}\right)^t \Gamma(t + 1).
	\end{align*}
	Since $\Gamma(t +1)^{1/t} \leq t$ for all $t \geq 1$ we conclude that for any $1 \leq s \leq t$,
	\begin{align*}
		\left(\mathrm{E}[\|X\|^t]\right)^{1/t} \leq 4 + 2 \frac{t}{s} \leq 6 \frac{t}{s}.
	\end{align*}
	This completes the proof.
\end{proof}

\begin{proof}[\textbf{Proof of Lemma~\ref{lemma:LDPGaussian}}]
	We first prove an asymptotic upper bound. Since
	\begin{align*}
		\limsup_{n \rightarrow \infty} \mathrm{E}[\|Z\|_p]/(\|\Sigma^{1/2}\|_{2 \rightarrow p} t_n) : = c < 1,
	\end{align*}
	for any $\delta \in \left(0,\frac{1-c}{2}\right)$ there exists $N > 0$ such that for all $n \geq N$, 
	\begin{align*}
		\mathrm{E}[\|Z\|_p] < (c + \delta) t_n \|\Sigma^{1/2}\|_{2 \rightarrow p}.
	\end{align*}
	Thus, for all $n \geq N$, 
	\begin{align*}
		&\mathrm{P} \left(\|Z\|_p > t_n\|\Sigma^{1/2}\|_{2 \rightarrow p}\right)\\
		&\quad{}=  \mathrm{P} \left(\|Z\|_p > \mathrm{E}[\|Z\|_p] +  t_n\|\Sigma^{1/2}\|_{2 \rightarrow p} - \mathrm{E}[\|Z\|_p] \right) \\
		& \quad{}\leq  \mathrm{P} \left(\|Z\|_p > \mathrm{E}[\|Z\|_p]+  t_n\|\Sigma^{1/2}\|_{2 \rightarrow p}  - (c + \delta)\sqrt{n} t_n\|\Sigma^{1/2}\|_{2 \rightarrow p} \right) \\
		&\quad{}\leq \mathrm{P} \left(\|Z\|_p > \mathrm{E}[\|Z\|_p]  + \frac{1-c}{2} t_n \|\Sigma^{1/2}\|_{2 \rightarrow p}\right) \\
		&\quad{} \leq e^{-\frac{(1-c)^2}{8}t_n^2},
	\end{align*}
	where the last inequality follows from the Gaussian concentration property of Lipschitz functions. Hence,
	\begin{align*}
		\lim_{n \rightarrow \infty} \frac{1}{t_n^2} \log  \mathrm{P} \left(\|Z\|_p >  t_n\|\Sigma^{1/2}\|_{2 \rightarrow p} \right) \leq  -\frac{(1-c)^2}{8}.
	\end{align*}
	Next, we establish a matching asymptotic lower bound. Let $u^* \in \mathbb{R}^d$ be such that
	\begin{align*}
		\|\Sigma^{1/2}u^*\|_2^2 = \sup_{\|u\|_q = 1}  \|\Sigma^{1/2}u\|_2^2 =\|\Sigma^{1/2}\|_{2 \rightarrow p}^2.
	\end{align*}
	Since $\mathrm{Var}(Z'u^*) =  \|\Sigma^{1/2}u^*\|_2^2$, Gordon's lower bound on Mill's ratio for normal random variables yields
	\begin{align*}
		\mathrm{P} \left(\|Z\|_p > t_n\|\Sigma^{1/2}\|_{2 \rightarrow p}\right) &\geq \mathrm{P} \left(Z'u^* > t_n\|\Sigma^{1/2}\|_{2 \rightarrow p}\right)\\
		&\geq \frac{1}{\sqrt{2\pi}} \frac{\sqrt{\mathrm{Var}(Z'u^*)}}{t_n\|\Sigma^{1/2}\|_{2 \rightarrow p}} \exp \left(- \frac{t_n\|\Sigma^{1/2}\|_{2 \rightarrow p}^2}{2\mathrm{Var}(Z'u^*)} \right)\\
		&=\frac{1}{\sqrt{2\pi}} \frac{e^{-t_n^2/2}}{t_n}.
	\end{align*}
	Hence,
	\begin{align*}
		\lim_{n \rightarrow \infty} \frac{1}{t_n^2} \log  \mathrm{P} \left(\|Z\|_p > t_n\|\Sigma^{1/2}\|_{2 \rightarrow p}\right) \geq - \frac{1}{2}.
	\end{align*}
	This completes the proof. Notice that our proof is so simple and elementary only because we consider large $t_n > 0$ that asymptotically dominate $\mathrm{E}[\|Z\|_p]$.
\end{proof}

\begin{proof}[\textbf{Proof of Lemma~\ref{lemma:LDPStrictly-Log-Concave}}]
	Since
	\begin{align*}
		\limsup_{n \rightarrow \infty}(\sqrt{\lambda_d} \wedge 1) \mathrm{E}[\|X\|_p]/(\|\Sigma^{1/2}\|_{2 \rightarrow p} t_n) : = c < 1,
	\end{align*}
	there exists $N > 0$ such that for all $n \geq N$, 
	\begin{align*}
		\mathrm{E}[\|X\|_p] < t_n  \left(c + \frac{1-c}{2}\right) \frac{\|\Sigma^{1/2}\|_{2 \rightarrow p}}{\sqrt{\lambda_d} \wedge 1} = t_n  \frac{1+c}{2} \frac{\|\Sigma^{1/2}\|_{2 \rightarrow p}}{\sqrt{\lambda_d} \wedge 1}.
	\end{align*}
	Thus, for all $n \geq N$, 
	\begin{align*}
		&\mathrm{P} \left(\|X\|_p > t_n\|\Sigma^{1/2}\|_{2 \rightarrow p}/ (\sqrt{\lambda_d} \wedge 1)\right)\\
		&\quad{}=  \mathrm{P} \left(\|X\|_p > \mathrm{E}[\|X\|_p] +  t_n\|\Sigma^{1/2}\|_{2 \rightarrow p}/ (\sqrt{\lambda_d} \wedge 1) - \mathrm{E}[\|X\|_p] \right) \\
		& \quad{}\leq  \mathrm{P} \left(\|X\|_p > \mathrm{E}[\|X\|_p]+  t_n\|\Sigma^{1/2}\|_{2 \rightarrow p}/ (\sqrt{\lambda_d} \wedge 1) - t_n\frac{1+c}{2} \frac{\|\Sigma^{1/2}\|_{2 \rightarrow p}}{\sqrt{\lambda_d} \wedge 1}\right) \\
		&\quad{}\leq \mathrm{P} \left(\|X\|_p > \mathrm{E}[\|X\|_p]  + t_n \frac{1-c}{2} \frac{\|\Sigma^{1/2}\|_{2 \rightarrow p}}{\sqrt{\lambda_d} \wedge 1}\right) \\
		&\quad{} \leq e^{-\frac{(1-c)^2}{8}t_n^2},
	\end{align*}
	where the last inequality follows from the isoperimetric inequality for log-concave measures~\citep[][Theorem 2.7]{ledoux2001concentration} combined with Proposition 1.2 in~\cite{ledoux2001concentration}. Hence,
	\begin{align*}
		\lim_{n \rightarrow \infty} \frac{1}{t_n^2} \log  \mathrm{P} \left(\|X\|_p > t_n\|\Sigma^{1/2}\|_{2 \rightarrow p} / (\sqrt{\lambda_d} \wedge 1)\right) \leq  -\frac{(1-c)^2}{8}.
	\end{align*}
\end{proof}

\begin{proof}[\textbf{Proof of Lemma~\ref{lemma:LDP-Log-Concave}}]
	We compute
	\begin{align}\label{eq:lemma:LDP-Log-Concave-1}
		\mathrm{Var}(\|X\|_p) \leq \mathrm{E}[\|X\|_p^2] \leq \sup_{\|u\|_q =1} \mathrm{E}[(X'u)^2] = \sup_{\|u\|_q=1}\ |\Sigma^{1/2}u\|_2^2= \|\Sigma^{1/2}\|_{2 \rightarrow p}^2.
	\end{align}
	By Chebychev's inequality,
	\begin{align*}
		\mathrm{P} \left( \|X\|_p \leq \mathrm{E}[\|X\|_p] +  4 \sqrt{\mathrm{Var}(\|X\|_p)} \right) \geq \frac{3}{4},
	\end{align*}
	and thus by Borell's inequality for log-concave measures~\citep[e.g.][Proposition 2.14]{ledoux2001concentration} and eq.~\eqref{eq:lemma:LDP-Log-Concave-1}, for all $t \geq 1$,
	\begin{align}\label{eq:lemma:LDP-Log-Concave-2}
		\mathrm{P} \left( \|X\|_p \geq \mathrm{E}[\|X\|_p] +  4 t \|\Sigma^{1/2}\|_{2 \rightarrow p}\right)  \leq e^{-t/2}.
	\end{align}
	Moreover, since
	\begin{align*}
		\limsup_{n \rightarrow \infty} \mathrm{E}[\|X\|_p]/(\|\Sigma^{1/2}\|_{2 \rightarrow p} t_n) : = c < 1,
	\end{align*}
	there exists $N > 0$ such that for all $n \geq N$, 
	\begin{align}\label{eq:lemma:LDP-Log-Concave-3}
		\mathrm{E}[\|X\|_p] < t_n  \left(c + \frac{1-c}{2}\right) \|\Sigma^{1/2}\|_{2 \rightarrow p} = t_n  \frac{1+c}{2} \|\Sigma^{1/2}\|_{2 \rightarrow p}.
	\end{align}
	Thus, for all $n \geq N$, 
	\begin{align*}
		&\mathrm{P} \left(\|X\|_p > t_n\|\Sigma^{1/2}\|_{2 \rightarrow p} \right)\\
		&\quad{}=  \mathrm{P} \left(\|Z\|_p > \mathrm{E}[\|Z\|_p] +  t_n\|\Sigma^{1/2}\|_{2 \rightarrow p} - \mathrm{E}[\|Z\|_p] \right) \\
		& \quad{}\leq  \mathrm{P} \left(\|Z\|_p > \mathrm{E}[\|Z\|_p]+  t_n\|\Sigma^{1/2}\|_{2 \rightarrow p} - t_n\frac{1+c}{2} \|\Sigma^{1/2}\|_{2 \rightarrow p} \right) \\
		&\quad{}\leq \mathrm{P} \left(\|Z\|_p > \mathrm{E}[\|Z\|_p]  + t_n \frac{1-c}{2}\|\Sigma^{1/2}\|_{2 \rightarrow p}\right) \\
		&\quad{} \leq e^{-\frac{(1-c)}{16}t_n},
	\end{align*}
	where the last inequality follows from eq.~\eqref{eq:lemma:LDP-Log-Concave-2}. Hence,
	\begin{align*}
		\lim_{n \rightarrow \infty} \frac{1}{t_n} \log  \mathrm{P} \left(\|Z\|_p > t_n\|\Sigma^{1/2}\|_{2 \rightarrow p} \right) \leq  -\frac{(1-c)}{16}.
	\end{align*}
\end{proof}

\begin{proof}[\textbf{Proof of Lemma~\ref{lemma:BoundsVarianceUniformRV}}]
	We adapt the proof strategy of Theorem 4.4.3 in~\cite{tanguy2017quelques} (p. 118f) to our setting. Among other things, we slightly improve his arguments and obtain a sharper bound for the variance of the largest coordinate of random vector uniformly distributed over the sphere. This minor improvement can be interpreted as a non-asymptotic analogue to the classical Weibull maximum domain of attraction result for the uniform distribution over the interval $[0,1]$~\cite[e.g.][p. 82]{tanguy2017quelques}.
	
	Throughout the proof, we denote by $\sigma_{n-1}$ the law of $U \sim \mathrm{Unif}(S^{n-1})$ and by $(P_t)_{t \geq 0}$ a Markov semigroup with stationary measure $\sigma_{n-1}$. Define the operator $D_{jk} := x_j \partial_k - x_k \partial_j$ for all $1 \leq j,k \leq n$ and let $\Gamma = [\gamma_1, \ldots, \gamma_n] = [\tilde{\gamma}_1, \ldots, \tilde{\gamma}_d]' \in \mathbb{R}^{d \times n}$. Since $\sigma_{n-1}$ satisfies the log-Sobolev inequality with constant $n-1$, it follows by eq. (16) in~\cite{cordero-erausquin2012hypercontractive} and the discussion prior to their Corollary 4, that, for all $f \in L^2(\sigma_{n-1})$ and for all $T \geq 1/(2(n-1))$,
	\begin{align}\label{eq:lemma:BoundsVarianceUniformRV-1}
		\mathrm{Var}_{\sigma_{n-1}}(f) \leq \frac{e}{n} \sum_{j,k}^n \int_0^T \mathrm{E}_{\sigma_{n-1}}\left[ P_t^2(D_{jk}f)\right] dt.
	\end{align}
	Consider the case $1 \leq p < \infty$.  For $\theta \in \mathbb{R}$ arbitrary, set $f(U) = e^{\theta \|\Gamma U\|_p/2}$. Then,
	\begin{align*}
		D_{jk} f(U) &= \frac{\theta}{2} \left(  \frac{\sum_{\ell=1}^d |\tilde{\gamma}_\ell'U|^{p-1} |\gamma_{\ell k}| U_j }{\|\Gamma U\|_p^{p-1}} -\frac{\sum_{\ell=1}^d |\tilde{\gamma}_\ell'U|^{p-1}|\gamma_{\ell j}|  U_k }{\|\Gamma U\|_p^{p-1}} \right) e^{\frac{\theta}{2}\|\Gamma U\|_p},
	\end{align*}
	and, hence, by eq.~\eqref{eq:lemma:BoundsVarianceUniformRV-1},
	\begin{align}\label{eq:lemma:BoundsVarianceUniformRV-2}
		\mathrm{Var}_{\sigma_{n-1}}\left(e^{\frac{\theta}{2} \|\Gamma U\|_p}\right) \leq \frac{4e}{n}\frac{\theta^2}{4} \sum_{j,k=1}^n \int_0^T \mathrm{E}_{\sigma_{n-1}}\left[ P_t^2\left(\frac{\sum_{\ell=1}^d |\tilde{\gamma}_\ell'U|^{p-1}|\gamma_{\ell k}| U_j }{\|\Gamma U\|_p^{p-1}}e^{\frac{\theta}{2} \|\Gamma U\|_p}\right)\right] dt.
	\end{align}
	By the hypercontractivity of the semi-group $(P_t)_{t \geq 0}$ (implied because $\sigma_{n-1}$ satisifies a log-Sobolev inequality!) it follows that, for all $1 \leq j,k\leq n$,
	\begin{align*}
		\mathrm{E}_{\sigma_{n-1}}\left[ P_t^2\left( \frac{\sum_{\ell=1}^d |\tilde{\gamma}_\ell'U|^{p-1} |\gamma_{\ell k}|  U_j }{\|\Gamma U\|_p^{p-1}}e^{\frac{\theta}{2} \|\Gamma U\|_p} \right)\right] \leq \left( \mathrm{E}_{\sigma_{n-1}}\left[ | U_j|^s \left(\frac{\sum_{\ell=1}^d |\tilde{\gamma}_\ell'U|^{p-1} |\gamma_{\ell k}|  }{\|\Gamma U\|_p^{p-1}} \right)^s e^{\frac{s\theta}{2} \|\Gamma U\|_p}   \right]\right)^{2/s},
	\end{align*}
	where $s \equiv s(t) = 1 + e^{-2(n-1)t}$. Since $1 \leq s \leq 2$, H{\"o}lder's inequality yields
	\begin{align*}
		\left( \mathrm{E}_{\sigma_{n-1}}\left[ | U_j|^s \left(\frac{\sum_{\ell=1}^d |\tilde{\gamma}_\ell'U|^{p-1} |\gamma_{\ell k}|  }{\|\Gamma U\|_p^{p-1}} \right)^s e^{\frac{s\theta}{2} \|\Gamma U\|_p}   \right]\right)^{2/s} \leq \mathrm{E}_{\sigma_{n-1}}\left[ U_j^2 \left(\frac{\sum_{\ell=1}^d |\tilde{\gamma}_\ell'U|^{p-1} |\gamma_{\ell k}|  }{\|\Gamma U\|_p^{p-1}} \right)^2 e^{\theta \|\Gamma U\|_p}\right].
	\end{align*}
	This bound combined with eq.~\eqref{eq:lemma:BoundsVarianceUniformRV-2} and $T = 1/n \geq 1/(2(n-1))$ gives
	\begin{align*}
		\mathrm{Var}_{\sigma_{n-1}}\left(e^{\frac{\theta}{2} \|\Gamma U\|_p}\right) &\leq  \frac{e\theta^2}{n^2} \sum_{j,k=1}^n \mathrm{E}_{\sigma_{n-1}}\left[ U_j^2 \left(\frac{\sum_{\ell=1}^d |\tilde{\gamma}_\ell'U|^{p-1} |\gamma_{\ell k}|  }{\|\Gamma U\|_p^{p-1}} \right)^2 e^{\theta \|\Gamma U\|_p}\right]\\
		&= \frac{e\theta^2}{n^2}\mathrm{E}_{\sigma_{n-1}}\left[\sum_{k=1}^n \left(\frac{\sum_{\ell=1}^d |\tilde{\gamma}_\ell'U|^{p-1} |\gamma_{\ell k}|  }{\|\Gamma U\|_p^{p-1}} \right)^2 e^{\theta \|\Gamma U\|_p}\right]\\
		&\leq \frac{e\theta^2}{n^2} \sup_{\|v\|_q=1}\sum_{k=1}^n (v'\gamma_k)^2 \mathrm{E}_{\sigma_{n-1}}\left[ e^{\theta \|\Gamma U\|_p}\right]\\
		&=\frac{e\theta^2}{n^2}\|\Gamma'\|_{q \mapsto 2}^2 \mathrm{E}_{\sigma_{n-1}}\left[ e^{\theta \|\Gamma U\|_p}\right]\\
		&= \frac{e\theta^2}{n^2} \|\Gamma\|_{2 \mapsto p}^2 \mathrm{E}_{\sigma_{n-1}}\left[ e^{\theta \|\Gamma U\|_p}\right],
	\end{align*}
	where $q \geq 1$ is the conjugate exponent to $p \geq 1$, i.e.  $1/p + 1/q = 1$. Since $\theta \in \mathbb{R}$ is arbitrary, Corollary 3.2 in~\cite{ledoux2001concentration} implies that there exists an absolute constant $c > 0$ such that for all $t \geq 0$,
	\begin{align*}
		\mathrm{P} \left( \left|\|\Gamma U\|_p - \mathrm{E}[\|\Gamma U\|_p] \right| \geq t \right) \leq 6 e^{-ct n / \|\Gamma\|_{2\rightarrow p} }.
	\end{align*}
	Hence, by Fubini, for all $1 \leq p < \infty$,
	\begin{align*}
		\mathrm{Var}(\|\Gamma U\|_2) \leq  \frac{12}{c} \frac{\|\Gamma\|_{2 \rightarrow p}^2}{n^2}.
	\end{align*}
	Next, consider the case $p = \infty$. The strategy is the same as in the previous case. For $\theta \in \mathbb{R}$ arbitrary, we now set $f(U) = e^{\theta \|\Gamma U\|_\infty/2}$. Also, note that $\mathbf{1}\{ | \sum_{i=1}^n \gamma_{i\ell} U_i| = \|\Gamma U\|_\infty \} = \mathbf{1}\{\ell = \arg\max_{1 \leq k \leq d} \left| \sum_{i=1}^n \gamma_{ik} U_i\right| \}$. Then, 
	\begin{align*} 
		D_{jk} f(U) &= \frac{\theta}{2} U_j \gamma_{k\ell}  \mathbf{1}\left\{ \left| \sum_{i=1}^n \gamma_{i\ell} U_i\right|  = \|\Gamma U\|_\infty  \right\} \mathrm{sign} \left( \sum_{i=1}^n \gamma_{i\ell} U_i \right) e^{\frac{\theta}{2} \|\Gamma U\|_\infty}\\
		&\quad - \frac{\theta}{2} U_k \gamma_{j\ell}  \mathbf{1}\left\{ \left| \sum_{i=1}^n \gamma_{i\ell} U_i\right| = \|\Gamma U\|_\infty \right\} \mathrm{sign} \left( \sum_{i=1}^n \gamma_{i\ell} U_i \right) e^{\frac{\theta}{2} \|\Gamma U\|_\infty},
	\end{align*}	
	and, hence, by eq.~\eqref{eq:lemma:BoundsVarianceUniformRV-1},
	\begin{align}\label{eq:lemma:BoundsVarianceUniformRV-3}
		\begin{split}
			&\mathrm{Var}_{\sigma_{n-1}}\left(e^{\frac{\theta}{2} \|\Gamma U\|_\infty}\right)\\
			&\quad{}\leq \frac{4e}{n}\frac{\theta^2}{4} \sum_{j,k=1}^n \int_0^T \mathrm{E}_{\sigma_{n-1}}\left[ P_t^2\left(U_j \gamma_{k\ell}  \mathbf{1}\left\{ \left| \sum_{i=1}^n \gamma_{i\ell} U_i\right|  = \|\Gamma U\|_\infty  \right\} \mathrm{sign} \left( \sum_{i=1}^n \gamma_{i\ell} U_i \right) e^{\frac{\theta}{2} \|\Gamma U\|_\infty}\right)\right] dt.
		\end{split}
	\end{align}
	As in the proof of statement (i), hypercontractivity of the semi-group $(P_t)_{t \geq 0}$ combined with H{\"o}lder's inequality yields, for all $1 \leq j,k \leq n$,
	\begin{align*}
		&\mathrm{E}_{\sigma_{n-1}}\left[ P_t^2\left(U_j \gamma_{k\ell}  \mathbf{1}\left\{ \left| \sum_{i=1}^n \gamma_{i\ell} U_i\right|  = \|\Gamma U\|_\infty  \right\} \mathrm{sign} \left( \sum_{i=1}^n \gamma_{i\ell} U_i \right) e^{\frac{\theta}{2} \|\Gamma U\|_\infty}\right)\right]\\
		& \quad{} \leq \mathrm{E}_{\sigma_{n-1}}\left[ U_j^2 \gamma_{k\ell}^2  \mathbf{1}\left\{ \left| \sum_{i=1}^n \gamma_{i\ell} U_i\right|  = \|\Gamma U\|_\infty  \right\} e^{\theta \|\Gamma U\|_\infty}\right].
	\end{align*}
	This bound combined with eq.~\eqref{eq:lemma:BoundsVarianceUniformRV-3} and $T = 1/n \geq 1/(2(n-1))$ gives 
	\begin{align*}
		\mathrm{Var}_{\sigma_{n-1}}\left(e^{\frac{\theta}{2} \|\Gamma U\|_\infty}\right) &\leq  \frac{eT \theta^2}{n} \sum_{j,k=1}^n  \mathrm{E}_{\sigma_{n-1}}\left[ U_j^2 \gamma_{k\ell}^2  \mathbf{1}\left\{ \left| \sum_{i=1}^n \gamma_{i\ell} U_i\right|  = \|\Gamma U\|_\infty  \right\} e^{\theta \|\Gamma U\|_\infty}\right]\\
		&=\frac{eT \theta^2}{n} \mathrm{E}_{\sigma_{n-1}}\left[\left( \sum_{k=1}^n \gamma_{k \ell}^2 \right) \mathbf{1}\left\{ \left| \sum_{i=1}^n \gamma_{i\ell} U_i\right|  = \|\Gamma U\|_\infty  \right\} e^{\theta \|\Gamma U\|_\infty}\right]\\
		&\leq \max_{1 \leq \ell \leq d} \left( \sum_{k=1}^n \gamma_{k \ell}^2  \right) \frac{eT \theta^2}{n} \mathrm{E}_{\sigma_{n-1}}\left[ e^{\theta \|\Gamma U\|_\infty}\right]\\
		&= \|\Gamma\|_{2 \rightarrow \infty}^2 \frac{eT \theta^2}{n} \mathrm{E}_{\sigma_{n-1}}\left[ e^{\theta \|\Gamma U\|_\infty}\right]\\
		&\leq \|\Gamma\|_{2 \rightarrow \infty}^2 \frac{e \theta^2}{n^2} \mathrm{E}_{\sigma_{n-1}}\left[ e^{\theta \|\Gamma U\|_\infty}\right].
	\end{align*}
	Since $\theta \in \mathbb{R}$ is arbitrary, we can invoke Corollary 3.2 in~\cite{ledoux2001concentration} and integrate out the tail probability to obtain
	\begin{align*}
		\mathrm{Var}(\|\Gamma U\|_\infty) \leq  \frac{12}{c} \frac{\|\Gamma\|_{2\rightarrow \infty}^2}{n^2}.
	\end{align*}
\end{proof}

\end{document}